\documentclass[10pt,reqno]{amsart}
\usepackage{amscd,amsmath,amsopn,amssymb,amsthm,multicol}
\usepackage{hyperref}
\DeclareMathOperator{\Diff}{Diff}

\DeclareMathOperator{\ad}{ad}
\DeclareMathOperator{\Id}{Id}
\DeclareMathOperator{\hor}{hor}

\DeclareMathOperator{\diag}{diag}
\DeclareMathOperator{\Rad}{Rad}
\DeclareMathOperator{\Ad}{Ad}

\DeclareMathOperator{\Lin}{Lin}
\DeclareMathOperator{\Aut}{Aut}

\DeclareMathOperator{\Isom}{Isom}
\DeclareMathOperator{\trace}{trace}
\DeclareMathOperator{\rk}{rk}
\DeclareMathOperator{\Radinj}{Radinj}

\renewenvironment{proof}[1][Proof]{\textbf{#1.} }
{\ \rule{0.5em}{0.5em}}

\newtheorem{theorem}{Theorem}
\newtheorem{pred}{Proposition}
\newtheorem{lemma}{Lemma}
\newtheorem{corollary}{Corollary}

\newtheorem{vopros}{Question}

\theoremstyle{definition}

\newtheorem{definition}{Definition}
\newtheorem{remark}{Remark}
\newtheorem{example}{Example}

\begin{document}

\title
[On $\delta$-homogeneous Riemannian manifolds]
{On $\delta$-homogeneous Riemannian manifolds}
\author{V.N.~Berestovski\u\i\ and Yu.G.~Nikonorov}

\address{Berestovski\u\i\  Valeri\u\i\  Nikolaevich \newline
Omsk Branch of Sobolev Institute of Mathematics SD RAS \newline
644099, Omsk, ul. Pevtsova, 13, Russia}

\email{berestov@ofim.oscsbras.ru}

\address{Nikonorov\ Yuri\u\i\  Gennadievich\newline
Rubtsovsk Industrial Institute \newline
of Altai State Technical University after I.I.~Polzunov \newline
658207, Rubtsovsk, Traktornaya, 2/6, Russia}

\email{nik@inst.rubtsovsk.ru}

\thanks
{The first author is supported in part by RFBR (grant N
04-01-00315-a). The second author is supported in part by RFBR
(grant N 05-01-00611-a) and by the Council on grants of the
President of Russian Federation for supporting of young russian
scientists and leading scientific schools of Russian Federation
(grants NSH-8526.2006.1 and MD-5179.2006.1)}

\begin{abstract}

We study in this paper previously defined by V.N.~Berestovskii and
C.P.~Plaut $\delta$-homogeneous spaces in the case of Riemannian
manifolds. Every such manifold has non-negative sectional
curvature. The universal covering of any $\delta$-homogeneous
Riemannian manifolds is itself $\delta$-homogeneous. In turn,
every simply connected Riemannian $\delta$-homogeneous manifold is
a direct metric product of an Euclidean space and compact simply
connected indecomposable homogeneous manifolds; all factors in
this product are itself $\delta$-homogeneous. We find different
characterizations of $\delta$-homogeneous Riemannian spaces, which
imply that any such space is geodesic orbit (g.o.) and every
normal homogeneous Riemannian manifold is $\delta$-homogeneous.
The g.o. property and the $\delta$-homogeneity property are
inherited by closed totally geodesic submanifolds. Then we find
all possible candidates for compact simply connected
indecomposable Riemannian $\delta$-homogeneous non-normal
manifolds of positive Euler characteristic and a priori
inequalities for parameters of the corresponding family of
Riemannian $\delta$-homogeneous metrics on them (necessarily
two-parametric). We prove that there are only two families of
possible candidates: non-normal (generalized) flag manifolds
$SO(2l+1)/U(l)$ and $Sp(l)/U(1)\cdot Sp(l-1)$, $l\geq 2$,
investigated earlier by W.~Ziller, H.~Tamaru, D.V.~Alekseevsky and
A.~Arvanitoyeorgos. At the end we prove that the corresponding
two-parametric family of Riemannian metrics on
$SO(5)/U(2)=Sp(2)/U(1)\cdot Sp(1)$ satisfying the above mentioned
(strict!) inequalities, really generates $\delta$-homogeneous
spaces, which are not normal and are not naturally reductive with
respect to any isometry group.

\vspace{2mm}

\vspace{2mm} \noindent Key words and phrases: homogeneous spaces,
homogeneous spaces of positive Euler characteristic, geodesic
orbit spaces, Clifford-Wolf translations, geodesics, normal
homogeneous Riemannian spaces, submetry, Riemannian submersion.

\end{abstract}

\maketitle

\tableofcontents

\section{Introduction}

Historically, the assembly of homogeneous Riemannian manifolds
under considerations has been gradually extended and about thirty
years ago it included all such manifolds. Nevertheless, the
division of them into particular classes is very important.

We shall mention here classes, which can be characterized by some
properties of their isometry groups with connection to their
geodesics. B.~Riemann separated all manifolds of constant sectional
curvature, which are characterized by the property of free
movability of figures. Later E.~Cartan introduced and classified
all \textit{symmetric} Riemannian manifolds. Then K.~Nomizu
introduced and studied \textit{naturally reductive} homogeneous
Riemannian manifolds which include as special cases symmetric
spaces and \textit{normal homogeneous} manifolds. The latter have
non-negative sectional curvatures and include all symmetric spaces
with nonnegative sectional curvatures. A little later A.~Selberg
introduced another generalization of symmetric spaces, namely,
\textit{weakly symmetric} Riemannian spaces. At last,
\textit{geodesic orbit (g.o.)} homogeneous Riemannian manifolds have been
discovered. This class includes properly as special all
previously mentioned classes, see \cite{KV}, \cite{Zi96}.
Every simply connected Riemannian g.o. manifold
of dimension $\leq 5$ is naturally reductive \cite{KV}.
Riemannian g.o. manifolds
of dimension $6$, which are not naturally reductive,
are classified in \cite{KV}; for recent results in dimension $7$
we refer to \cite{DKN}.
A.~Selberg proved that every weakly symmetric Riemannian manifold
$M$ is \textit{commutative}, i.e. it admits a transitive motion Lie group
$G$ such that the Lie algebra of
$G$-invariant differential operators on $M$ is commutative \cite{Selb}.
If $G$ is connected and $M=G/H$, the latter is equivalent to the property that
the functional space $L^1(H\backslash G/H)$ is commutative, i.e
$(G,H)$ is a \textit{Gelfand pair}, or the property that
for every unitary irreducible representation of $G$, the dimension
of $H$-fixed set is $\leq 1$, i.e.
$(G,H)$~ is a \textit{spherical pair}.
J.~Lauret obtained an example of commutative
non weakly symmetric Riemannian manifold \cite{Laur}.
On the other hand, if $(G,H)$ is a spherical pair with compact
simple Lie group $G$ and its closed subgroup $H$,
then $G/H$ is weakly symmetric \cite{Ng}. A classification of such pairs
is known from \cite{Kr}, \cite{Ng}, \cite{AkhVinb}.
Let us
remark that besides symmetric spaces, there is no complete classification of
manifolds in other classes, mentioned above, although normal
homogeneous manifolds do not require in some sense such
classification.

We prove in this paper that the previously defined in \cite{BerP}
$\delta$-homogeneous spaces constitute in the case of Riemannian
manifolds  a new class of homogeneous manifolds situated
between normal homogeneous and g.o. manifolds. These manifolds,
unlike all previously mentioned classes, have very simple, purely
metric definition, which really can be applied to any metric
space. Namely, an arbitrary metric space $(M,\rho)$ is called
\textit{$\delta$-homogeneous}, if for every two points $x,y\in M$
there is an isometry $f$ of $(M,\rho)$ onto itself, which moves
$x$ to $y$ and has the maximal displacement at the point $x$, i.e.
$f(x)=y$ and $\rho(x,f(x))\geq \rho(z,f(z))$ for all $z\in M$. If
we can always take such a motion $f$ from an isometry group $G$ of
$(M,\rho)$, then $(M,\rho)$ is called
\textit{$G$-$\delta$-homogeneous}. In the Riemannian case we shall take
as $G$ only connected transitive Lie groups.

The consideration and methods in this paper go from general to
more and more specific.

In Section 2 we bring main definitions, earlier results, and
simple examples. In particular, any Lie group with a bi-invariant
Riemannian metric or any direct metric product of
$\delta$-homogeneous spaces is $\delta$-homogeneous. Every
$\delta$-homogeneous Riemannian manifold has non-negative sectional
curvature (Proposition \ref{pos}).

In Section 3 we get general results on $\delta$-homogeneous
Riemannian manifolds $(M,\mu)$. If $(M,\mu)$ is $G$-$\delta$-homogeneous
and $G$ normalizes a closed subgroup $H$ of the full
isometry group $I(M)$ of $(M,\mu)$, then the quotient (orbit)
space $H\backslash M$ with the quotient Riemannian metric is
$\delta$-homogeneous (Theorem \ref{delta1}). As a corollary
(\ref{NormCol}), we get that every normal homogeneous Riemannian
manifold is $\delta$-homogeneous. Then we prove that the universal
locally isometric covering of $(M,\mu)$ is $\delta$-homogeneous
(Corollary \ref{covr}); $(M,\mu)$ is either compact or it is
isometric to a direct metric product of an Euclidean space and
some compact $\delta$-homogeneous Riemannian manifold (Theorem
\ref{product}). Since any homogeneous Riemannian manifold
$(M,\mu)$ is an orbit space of its universal covering
$(\tilde{M},\tilde{\mu})$ by central discrete subgroup $\Gamma$ in
the (unit) connected component of $I(\tilde{M})$, where $\Gamma$
is isomorphic to $\pi_1(M)$, then previous results imply that the
study of $\delta$-homogeneous Riemannian manifolds entirely
reduces to the simply connected compact case. Then we get four
useful necessary and sufficient conditions for a (homogeneous)
connected Riemannian manifolds $(M,\mu)$ to be
$\delta$-homogeneous; we will mention here two of them. First:
$(M=G/H,\mu)$ with the corresponding inner metric $\rho$ is
$G$-$\delta$-homogeneous if and only if it is $G$-normal in the
generalized sense (Corollary \ref{general1}). The latter means
that there is a bi-invariant Finsler (inner) metric $F$ on $G$
such that the natural projection $p:(G,F)\rightarrow (G/H,\rho)$
is submetry (see Definition \ref{su}). Notice that in the case
when $F$ is Riemannian (inner), $(M,\mu)$ would be $G$-normal.
Second: $(M=G/H,\mu)$ is $G$-$\delta$-homogeneous if and only if
every geodesic $\gamma$ in $(M,\mu)$ is an orbit of a 1-parameter
motion group of $(M,\mu)$ in $G$, generated by a Killing vector
field, attaining a maximal value of its length at $\gamma$
(Theorem \ref{killi}). As a corollary (\ref{goc}), every
$\delta$-homogeneous Riemannian manifold is geodesic orbit
(g.o.). At the same time, g.o. Lobachevski's space of constant
negative curvature cannot be $\delta$-homogeneous by Proposition
\ref{pos}.

In Section 4 we prove that every closed totally geodesic
submanifold of a $\delta$-homogeneous (respectively, g.o.)
Riemannian manifold is $\delta$-homogeneous (respectively, g.o.)
itself, see Theorem \ref{gon2} (respectively, \ref{gon3}). As a
corollary (Theorem \ref{Ram}), every factor of a
$\delta$-homogeneous or g.o. direct metric product has the same
property. By all previous results, the study of all
$\delta$-homogeneous Riemannian spaces reduces to the compact
simply connected indecomposable case; we can separate further the
cases of zero or positive Euler characteristic. In the second half
of the section we find some algebraic properties of geodesic
vectors on a homogeneous Riemannian space $(G/H,\mu)$, i.e.
vectors in the Lie algebra $\mathfrak{g}$ of the Lie group $G$,
that tangent to a 1-parameter subgroup in $G$ with geodesic orbit
through the point $H\in (G/H,\mu)$. Besides later applications, we
use them to prove that if $(G/H,\mu)$ is $G$-$\delta$-homogeneous
and $L$ is a Lie subgroup of $G$ such that $H\subset L\subset G$,
then $L/H$ with the metric, induced by $\mu$, is
$\delta$-homogeneous.

In Section 5 we find additional isometries of $\delta$-homogeneous
Riemannian manifolds and find some applications of them.

Similarly to geodesic vectors, we define a $\delta$-vector on
$(G/H,\mu)$ as a vector in $\mathfrak{g}$ that tangent to (the
unique) right-invariant vector field $Y$ on $G$ with the property
that the Killing vector field $X=dp(Y)$ on $(G/H,\mu)$ has maximal
value of its length at the point $H\in (G/H,\mu)$. Remark that
every $\delta$-vector is a geodesic vector. In Section 6 we find
general properties of $\delta$-vectors. We use essentially these
properties later. The space $(G/H,\mu)$ is
$G$-$\delta$-homogeneous if and only if every vector $v\in
T_{H}(G/H)$ can be represented in the form $v=dp(w)$ for some
$\delta$-vector $w$.

In Section 7 we give (mainly known) results on compact simply
connected homogeneous spaces $M=G/H$, in particular, Hopf-Samelson
Theorem \ref{On}, which implies that $\chi(M)\geq 0$ and
characterizes the case, when $\chi(M)>0$, by the condition
$\rk(G)=\rk(H)$. Let us mark Theorem \ref{pose}, which states that
every proper Lie subalgebra $\mathfrak{h}$ of the Lie algebra
$\mathfrak{g}$ of a simple compact connected Lie group $G$,
containing the Lie algebra $\mathfrak{t}$ of a maximal torus
$T\subset G$, is the Lie algebra of the unique closed connected
Lie subgroup $H\subset G$. Moreover, $M=G/H$ is a simply connected
compact connected homogeneous space of positive Euler
characteristic. This gives an algebraic description of all simply
connected compact indecomposable homogeneous Riemannian
manifolds $(M,\mu)$ with $\chi(M)>0$ by the Kostant's Theorem
\ref{KostN}.

In Section 8 we give  known results on compact simply connected
homogeneous spaces of positive Euler characteristic. Also we prove
Theorem \ref{Kon57} which implies that every naturally reductive
homogeneous Riemannian manifolds of positive Euler characteristic
is normal (hence, $\delta$-homogeneous).

In Section 9 algebraic corollaries of $\delta$-homogeneity of the
first and the second order are found in Theorem \ref{ncdo}. Note
that the first order condition is simply the condition for
geodesic vectors.

In the next sections we consider only compact simply connected
homogeneous spaces of positive Euler characteristic.

In Section 10 we find some algebraic identities and inequalities
for $\delta$-homogeneous manifolds of one special type. As we
shall prove in Section 13, any (compact simply connected)
indecomposable $\delta$-homogeneous non-normal Riemannian manifold
$(M,\mu)$ with $\chi(M)>0$ must have such type. Especially
important are Propositions \ref{t31.4}, \ref{t31.5}, and
\ref{t31.9n}.

In Section 11 is given necessary information on roots and
structural constants of compact simple Lie algebras with respect
to their Killing forms and Cartan subalgebras. Mark especially the
identity (\ref{length}).

We show in Section 12 that all $G_2$-$\delta$-homogeneous
Riemannian metrics on homogeneous spaces with positive Euler
characteristic are normal.

Observe that roots of (simple compact) Lie algebras
$A_l$, $D_l$, $e_6$, $e_7$ and $e_8$ have one and the same length,
while the roots of Lie algebra $g_2$, $B_l$, $C_l$ and $f_4$ have
two different lengths. One knows also that the Weyl group of
any simple Lie algebra acts transitively on the set of roots with
equal lengths. With the help of these facts and the identity
(\ref{length}) we prove in Section 13 that the set of
$G$-$\delta$-homogeneous Riemannian metrics on $G/H$ with
$\chi(G/H)>0$ and compact simple Lie group $G$ is one- or
two-parametric; we have necessarily the first case (that is, only
$G$-normal metrics), if the Lie algebra $\mathfrak{g}$ has roots
of equal length, as for Lie groups $G=SU(l+1)$, $SO(2l)$, $E_6$,
$E_7$, $E_8$ (Proposition \ref{two} and Corollary \ref{onel}). We
shall have only $G$-normal metrics also in the case, when $H=T$
(Proposition \ref{rflag}). In the case of two-parametric family we
get with the help of Proposition \ref{t31.9n} a priori
inequalities (\ref{three}) for these parameters in Proposition
\ref{ineq}.

Further investigations of two-parametric case in Section 14 shows
that possible candidates one can find only among flag manifolds $SO(2l+1)/U(l)$
and $Sp(l)/U(1)\cdot Sp(l-1)$,
where $l\geq 2$.
All invariant metrics on these manifolds are weakly symmetric (hence, g.o.)
\cite {Zi96}. Moreover, among (generalized) flag manifolds only
$SO(2l+1)/U(l)$ and $Sp(l)/U(1)\cdot Sp(l-1)$, $l\geq 2$,
admit non-normal invariant g.o. metrics \cite{AA}.

Using Proposition \ref{t31.4} and spectra of matrices, we prove in
Section 15 that two-parametric family of Riemannian metrics on
$SO(5)/U(2)$, which satisfies the inequalities (\ref{three}),
really give us $SO(5)$-$\delta$-homogeneous spaces. The limiting
cases of this inequalities represent $SO(5)$-normal and
$SO(6)$-normal spaces respectively, while all other metrics are
$SO(5)$-$\delta$-homogeneous and non-normal (Theorem \ref{main}).
We are planning to investigate all other possible cases, mentioned
in the previous paragraph, separately.

Some unsolved questions are posed in different places of the text.

The first author is very obliged to Mathematics Department of
University of Tennessee, Knoxville, USA, for hospitality and
visiting position while a part of this paper have been prepared.

\section{Preliminaries}

\begin{definition}
Let $(X,d)$ be a metric space and $x\in X$. An isometry $f:X
\rightarrow X$ is called {\it a $\delta(x)$-translation} ({\it a
Clifford-Wolf translation}), if $x$ is a point of maximal
displacement of $f$, i.e. for every $y \in X$ the relation
$d(y,f(y))\leq d(x,f(x))$ holds (respectively, $f$ displaces all
points of $(X,d)$ the same distance, i.e. $d(y,f(y))=d(x,f(x))$
for every $y \in X$).
\end{definition}

\begin{definition}\label{Gener}
A metric space $(X,d)$ is called {\it (G)-$\delta$-homogeneous}
(respectively, {\it (G)-Clifford-Wolf homogeneous}), if for every
$x,y \in X$ there exists a $\delta(x)$-translation (respectively,
{\it Clifford-Wolf translation}) of $(X,d)$ (from an isometry
group $G$), moving $x$ to $y$.
\end{definition}

It is clear that any Clifford-Wolf translation is a
$\delta(x)$-translation for every point $x\in X$, any
($G$)-Clifford-Wolf homogeneous space is
($G$)-$\delta$-homogeneous, and the latter one is
($G$)-homogeneous.

\begin{example}
Every Lie group with a bi-invariant inner metric $(G,r)$ and every
odd-dimen\-sional Euclidean sphere (of the unit radius)
$S^{2n+1}\subset \mathbb{E\,}^{2(n+1)}$ with the induced inner
(Riemannian) metric is Clifford-Wolf homogeneous space. In the
first case it is enough to use left translations on some fixed
element of the group. The second statement is proved essentially
by Clifford himself, which explains the term in Definition
\ref{Gener}.
\end{example}

\begin{example}
One can easily see that a direct metric product of
$\delta$-(respectively, Clifford-Wolf) homogeneous spaces is again
$\delta$-(respectively, Clifford-Wolf) homogeneous.
\end{example}

In the paper \cite{BerP} the following results are obtained.

\begin{theorem}[Berestovskii-Plaut \cite{BerP}]\label{BerP0}
Every locally compact
$\delta$-homogeneous space of curvature bounded below
in the sense of Alexandrov
has non-negative curvature.
\end{theorem}

\begin{theorem}[Berestovskii-Plaut \cite{BerP}]\label{BerP}
Every non-compact locally compact homogeneous inner metric space
of non-negative curvature in the sense of Alexandrov is isometric
to a direct metric product of finite dimensional Euclidean space
and a compact homogeneous inner metric space of non-negative
curvature.
\end{theorem}

\begin{remark} In the Riemannian case, the last theorem easily follows
from Toponogov's theorem in \cite{T}, which states that every
complete Riemannian manifold $(M,\mu)$ with nonnegative sectional
curvature, containing a metric line, is isometric to a direct
Riemannian product $(N,\nu)\times \mathbb{R}$. Later J.~Cheeger and
D.~Gromoll in \cite{CG} generalized Toponogov's theorem to complete
Riemannian manifolds of nonnegative Ricci curvature.
\end{remark}

If $(M,\mu)$ is a Riemannian manifold with inner metric $\rho$,
then Theorem \ref{BerP0} implies

\begin{pred}\label{pos}
Every $\delta$-homogeneous Riemannian manifold $(M,\rho)$
has non-negative sectional curvature.
\end{pred}

\begin{definition}\label{su}
A map of metric spaces $f:(M,r)\rightarrow (N,q)$ is called a
\textit{submetry}, if it maps every closed ball $B(x,s)\subset
(M,r)$ with the radius $s$ and the center $x$ onto the closed ball
$B(f(x),s)\subset (N,q)$ with the radius $s$ and the center
$f(x)$, \cite{BG}.
\end{definition}

Note that a smooth map of complete Riemannian spaces is submetry
if and only if it is a Riemannian submersion \cite{BG}.

\begin{definition}\label{nor}
A locally compact inner metric (respectively, Riemannian) space
$(M=G/H,\rho)$ with a transitive locally compact topological
(respectively, Lie) group $G$ and a stabilizer subgroup $H$ at a
point $x\in M$ is called {\it $G$-normal in generalized
(respectively, usual) sense}, if $G$ admits a bi-invariant
(respectively, Riemannian bi-invariant) inner metric $r$ such that
the natural projection $(G,r)\rightarrow (G/H,\rho)$ is a
submetry.
\end{definition}

\section{General properties of $\delta$-homogeneous spaces}\label{prostr}

\begin{definition}\label{R}
An inner metric space $(M,\rho)$ is called {\it restrictively
($G$)-$\delta$-homogeneous} (respectively, {\it restrictively
($G$)-Clifford-Wolf homogeneous}) if for every $x\in M$ there
exists a number $r(x)>0$ such that for every two points $y,z$ in
the open ball $U(x,r(x))$ there exists a $\delta(y)$-translation
(respectively, a Clifford-Wolf translation) of the space $(M,
\rho)$ (from the isometry group $G$), moving $y$ to $z$. The
supremum $R(x)$ of all such numbers $r(x)$ is called {\it the
($G$)-$\delta$-homogeneity radius} (respectively, {\it the
($G$)-Clifford-Wolf homogeneity radius}) of the space $(M,\rho)$
at the point $x$.
\end{definition}

\begin{pred}\label{pr}
Every restrictively ($G$)-$\delta$-homogeneous locally compact
complete inner metric space is ($G$)-$\delta$-homogeneous.
\end{pred}

\begin{proof}
It is clear that (in the notation of Definition \ref{R}) the
function $R(x)$, $x\in M$, is equal identically to $+\infty$, i.e.
the space $(M,\rho)$ is ($G$)-$\delta$-homogeneous, or it
satisfies the inequality $|R(x_{1})-R(x_{2})|\leq
\rho(x_{1},x_{2})$. In the last case the function $R(x)$, $x\in
M$, is positive and continuous.

Let us consider arbitrary points $x,y$ of a metric space
$(M,\rho)$, and suppose that this space satisfies the above-stated
condition. Then one can join the points $x$ and $y$ by some
shortest $[x,y]$. According to the above discussion, one can
divide sequentially this shortest by points $x_{0}=x, x_{1},\dots
,x_{m}=y$ such that for every $l$, where $0\leq l\leq m-1$, there
exists a $\delta(x_{l})$-translation $f_{l}$ of the space
$(M,\rho)$ (from the group $G$), moving the point $x_{l}$ to the
point $x_{l+1}$. Now the triangle inequality implies that the
composition $f:=f_{m-1}\circ \dots \circ f_{0}$ is a
$\delta(x)$-translation of the space $(M,\rho)$ (from the group
$G$), moving the point $x$ to the point $y$.
\end{proof}

\begin{theorem}\label{delta1}
Let $(M,r)$ be a locally compact inner metric space which is
$G$-$\delta$-homogeneous. Suppose that the group $G$ normalizes
some closed subgroup $H$ of the full isometry group $\Isom(M)$ of
$M$ (supplied by the compact-open topology). Then the quotient
(orbit) space $H \backslash M$ with the quotient metric $\rho$ is
a ($G$)-$\delta$-homogeneous (locally compact inner metric) space.
\end{theorem}

\begin{proof}
According to  S.E.~Cohn-Vossen theorem \cite{KF}, every complete
locally compact inner metric space is finitely-compact, i.e.,
every its closed bounded subset is compact. It is proved in the
paper \cite{B} that any closed subgroup of the full isometry group
(with the compact-open topology) of arbitrary finitely-compact
space has closed orbits. This implies that the group $H$ has
closed orbits in $M$.

On the ground of this fact it is easy to prove that the canonical
projection $p:(M,r)\rightarrow (H \backslash M,\rho)$ is a
submetry. This is equivalent to the following two properties:

1) the map $p$ does not increase distances;

2) for every three points $x,y\in H \backslash M$, $\xi \in p^{-1}(x)$,
there exists a point $\eta\in p^{-1}(y)$ such that
$r(\xi,\eta)=\rho(x,y)$.

Now let us consider arbitrary points $x,y\in H \backslash M$ and
the corresponding points $\xi, \eta$ from Property 2). By
condition there is a $\delta(\xi)$-translation $F$ of the space
$(M,r)$ from the group $G$ such that $F(\xi)=\eta$. Since the
group $G$ normalizes the group $H$, there is an isometry $f$ of
the space $(H \backslash M,\rho)$, induced by the isometry $F$.
Moreover, $f(x)=p(F(\xi))=p(\eta)=y$. Now for any point
$z=p(\zeta)\in H \backslash M$ Properties 1) and 2) imply the
relations

$$
\rho(x,f(x))=\rho(x,y)=r(\xi,\eta)= r(\xi,F(\xi))\geq
$$
$$
r(\zeta, F(\zeta)) \geq \rho(p(\zeta),p(F(\zeta)))=\rho(z,f(z)),
$$
i.e. $f$ is a $\delta(x)$-translation of the space $(H \backslash
M,\rho)$ moving the point $x$ to the point $y$. Therefore, the
space $(H \backslash M,\rho)$ is $G$-$\delta$-homogeneous.
\end{proof}

\begin{corollary}\label{NormCol}
Every ($G$)-normal in the generalized sense homogeneous locally
compact inner metric space is ($G$)-$\delta$-homogeneous. As a
corollary, any ($G$)-normal (maybe, in the generalized sense)
homogeneous Riemannian manifold is ($G$)-$\delta$-homogeneous.
\end{corollary}

\begin{proof}
Let a ($G$)-normal (in the generalized sense) homogeneous space
under consideration be a (metric) quotient space $(G/H,\rho)$ of a
locally compact topological group $(G,r)$ with a bi-invariant
inner metric $r$ by its compact subgroup $H$. Then the group of
left translations of the group $(G,r)$ is a transitive group of
Clifford-Wolf translations, and it commutes with the group of
right translations by elements of the subgroup $H$ which consists
of some isometries of the space $(G,r)$. Now it is enough to use
Theorem \ref{delta1}.
\end{proof}

\begin{pred}\label{cover}
The universal locally isometric covering of a $\delta$-homogeneous
(respectively, a restrictively Clifford-Wolf homogeneous)
Busemann's $G$-space is a $\delta$-homogeneous (respectively, a
restrictively Clifford-Wolf homogeneous) Busemann's $G$-space.
\end{pred}

\begin{proof}
Busemann's $G$-spaces are defined in his book \cite{Bus}.

Let $p:(\tilde{M},\tilde{\rho})\rightarrow (M,\rho)$ be the
universal locally isometric covering map for a
$\delta$-homogeneous (respectively, a restrictively Clifford-Wolf
homogeneous) Busemann's $G$-space $(M,\rho)$. It is clear that
$(\tilde{M},\tilde{\rho})$ is a Busemann's $G$-space. By Theorem
28.10 in \cite{Bus}, the group $G$ of all motions of the space
$(\tilde{M},\tilde{\rho})$, which cover motions of the space
$(M,\rho)$, is transitive on $\tilde{M}$, and the group $\Gamma$
of deck transformations of the covering $p$ is a normal subgroup
of the group $G$. Therefore, there is a number $r>0$ such that the
map $p$ is isometry on every open ball $U(x,r)\subset
(\tilde{M},\tilde{\rho})$ .

According to Proposition \ref{pr}, it is enough to show that the
space $(\tilde{M},\tilde{\rho})$ is restrictively
$\delta$-homogeneous (respectively, restrictively Clifford-Wolf
homogeneous). Consider arbitrary points  $x,y$ in
$(\tilde{M},\tilde{\rho})$ with the condition $\tilde{\rho}(x,y)<
r$. Since $(M,\rho)$ is $\delta$-homogeneous (respectively,
restrictively Clifford-Wolf homogeneous), there is a
$\delta(p(x))$-translation (respectively, a Clifford-Wolf
translation) $f$ of the space $(M,\rho)$ such that $f(p(x))=p(y)$.
From the above discussion we get that there is the unique map $F$
of the space $(\tilde{M},\tilde{\rho})$ onto itself covering the
map $f$ such that $F(x)=y$. It is clear that $F$ is an isometry of
the space $(\tilde{M},\tilde{\rho})$ and also a
$\delta(x)$-translation (respectively, a Clifford-Wolf
translation). This means that the space $(\tilde{M},\tilde{\rho})$
is restrictively $\delta$-homogeneous (respectively, restrictively
Clifford-Wolf homogeneous).
\end{proof}

\begin{corollary}\label{covr}
The universal Riemannian covering of a $\delta$-homogeneous
(respectively, a restrictively Clifford-Wolf homogeneous)
Riemannian manifold is $\delta$-homogeneous (respectively,
restrictively Clifford-Wolf homogeneous).
\end{corollary}

\begin{lemma}\label{eucl}
Suppose that the Riemannian manifold $(M,\mu)$ is isometric to the
direct metric product $(K,\mu_{1})\times
(\mathbb{E\,}^{m},\mu_{2})$, where $(K,\mu_{1})$ is a compact
homogeneous Riemannian manifold, and $(\mathbb{E\,}^{m},\mu_{2})$
is a finite dimensional Euclidean space. Then every isometry $f$
of the space $(M,\mu)$ has the form $f=f_{1}\times f_{2}$, where
$f_{1}$ (respectively, $f_{2}$) is an isometry of the space
$(K,\mu_{1})$ (respectively, $(\mathbb{E\,}^{m},\mu_{2})$).
\end{lemma}

\begin{proof}
It is easy to see that a geodesic in $(M,\mu)$ is a metric line if
and only if it is situated in some Euclidean subspace $\{k\}\times
\mathbb{E\,}^{m}$. Therefore, any isometry $f$ of the space
$(M,\mu)$ transposes such subspaces. Since $f$ keeps the
orthogonality, $f$ must transpose also all fibers of the form
$K\times \{e\}$. This proves Lemma.
\end{proof}

\begin{lemma}\label{dp}
If $M=M_1\times M_2$ is a direct product of Riemannian manifolds,
then every its isometry of the form $f=f_1\times f_2$ is a
$\delta(x)$-translation for the point $x=(x_1,x_2)\in M$ if and
only if both isometries $f_1: M_1\rightarrow M_1$ and
$f_2:M_2\rightarrow M_2$ are $\delta$-translations at the points
$x_1\in M_1$ and $x_2\in M_2$ respectively.
\end{lemma}

\begin{proof}
Let us remind that
$$
\rho((x_1,x_2),(y_1,y_2))=\sqrt{\rho_1^2(x_1,y_1)+\rho_2^2(x_2,y_2)},
$$
where $\rho$, $\rho_1$, $\rho_2$ are inner metrics of spaces $M$,
$M_1$, $M_2$ respectively. This easily implies the sufficiency.
Suppose that $f=f_1\times f_2$ is a $\delta$-translation of the
space $M$ at the point $x=(x_1,x_2)$, but, for instance, $f_1$  is
not a $\delta$-translation at the point $x_1$. Then there is a
point $x'_1$ such that
$\rho_1(x'_1,f_1(x'_1))>\rho_1(x_1,f_1(x_1))$. Therefore,
$$
\rho((x_1,x_2),f(x_1,x_2))=\sqrt{\rho_1^2(x_1,f_1(x_1))+
\rho_2^2(x_2,f_2(x_2))}<
$$
$$
\sqrt{\rho_1^2(x'_1,f_1(x'_1))+\rho_2^2(x_2,f_2(x_2))}=
\rho((x'_1,x_2),f(x'_1,x_2)),
$$
which contradicts to assumptions of Lemma.
\end{proof}

\begin{theorem}\label{product}
Any $\delta$-homogeneous Riemannian manifold $(M,g)$ is either
compact, or it is isometric to the direct metric product of an
Euclidean space and some compact $\delta$-homogeneous Riemannian
manifold.
\end{theorem}

\begin{proof}
This theorem immediately follows from Proposition \ref{pos},
Theorem \ref{BerP}, Lemma \ref{eucl} and Lemma \ref{dp}.
\end{proof}

From Theorem \ref{product} and Proposition \ref{cover} we
immediately obtain

\begin{corollary}\label{pi}
The universal Riemannian covering $(\tilde{M},\tilde{\mu})$ of a
$\delta$-homogeneous compact Riemannian manifold $(M,\mu)$ is
compact if and only if $\pi_{1}(M)$ is finite. In the opposite
case $(\tilde{M},\tilde{\mu})$ is isometric to a nontrivial direct
metric product of a compact simply connected $\delta$-homogeneous
Riemannian space and an Euclidean space.
\end{corollary}

\begin{theorem}\label{pr1}
A homogeneous space $M=G/H$ of a connected Lie group $G$ by its
compact subgroup $H$ admits an invariant Riemannian
$\delta$-homogeneous metric if and only if $G/H$ admits an
invariant Riemannian metric of non-negative sectional curvature.
\end{theorem}

\begin{proof}
The necessity follows from Proposition \ref{pos}.

Let us prove the sufficiency. Suppose that $M=G/H$ admits an
invariant Riemannian metric $\mu$ of non-negative sectional
curvature.

If $M$ is compact, then the Lie group $G$ is compact and it admits
a bi-invariant Riemannian metric $\gamma$. Then there is an unique
Riemannian metric $\nu$ on $M$ such that the canonical projecture
$p:(G,\gamma)\rightarrow (M,\nu)$ is a Riemannian submersion.
Moreover, $\nu$ is invariant on $G/H$, and $(G/H,\nu)$ is a
$G$-normal homogeneous Riemannian manifold. According to Corollary
\ref{NormCol}, $(G/H,\nu)$ is a $\delta$-homogeneous space.

Suppose, that $M$ is noncompact. Then by Theorem \ref{BerP}, all
assumptions of Lemma \ref{eucl} are fulfilled, moreover,
$(K,\mu_{1})$ has non-negative sectional curvature. Therefore, the
Lemma \ref{eucl} is valid. Obviously, the set of all isometries of
the type $\{f_{1}| f=(f_1,f_2)\in G\}$ forms a precompact
transitive isometry group $G_{1}$ of the compact space
$(K,\mu_{1})$ (relatively to the compact-open topology) with the
closure $\Gamma_{1}:=\overline{G_{1}}$, which is a compact
effective transitive isometry Lie group of the space
$(K,\mu_{1})$. Consequently, the manifold $K$ admits a
$\Gamma_{1}$-invariant Riemannian metric $\gamma_{1}$ such that
$(K,\gamma_{1})$ is a normal homogeneous space of the Lie group
$\Gamma_{1}$. According to Corollary \ref{NormCol},
$(K,\gamma_{1})$ is a $\delta$-homogeneous space. The last
reasonings imply that the Riemannian metric
$g_{0}=\gamma_{1}\times \mu_{2}$ on $M$ is invariant under the
action of the group $G$. In this case the Riemannian manifold
$(M,\mu_{0})=(K,\gamma_{1})\times (\mathbb{E\,}^{m},\mu_{2})$ is a
$\delta$-homogeneous space as a direct metric product of
$\delta$-homogeneous spaces.
\end{proof}

\begin{theorem}\label{D}
Let $(M,\mu)$ be a smooth connected compact Riemannian manifold
with inner metric $\rho$, and $G$ be the identity component of the
full isometry group of $(M,\mu)$. Then the function $d:G\times
G\rightarrow \mathbb{R}$ defined by the formula
\begin{equation}
\label{d} d(g,h)=\max_{x\in M}\rho(g(x),h(x)),
\end{equation}
determines a  bi-invariant metric on $G$ compatible with its
compact-open topology. In this case $(G,d)$ is locally isometric
to $(G,D)$ for some bi-invariant inner metric $D$ on $G$. Under
identification of the Lie algebra $G_{e}$ of the group $G$ with
the Lie algebra of Killing vector fields on $(M,\mu)$, $D$
coincides with the bi-invariant Finsler metric on $G$, determined
by the $\Ad(G)$-invariant norm $||\cdot||$ on $G_{e},$ defined by
the formula
\begin{equation}
\label{no} ||X||=\max_{x\in M}\sqrt{\mu(X(x),X(x))}.
\end{equation}
\end{theorem}

\begin{proof}
One can check directly the bi-invariance of the metric $d$. The
compactness of $(M,\mu)$ implies the compactness of the Lie group
$G$. Then, since $G$ is connected, the exponential map of the Lie
algebra $G_{e}$ to $G$ is surjective.

Let $g\neq e$ be arbitrary element in $G$. Then $g=\exp(X)$ for
some suitable Killing vector field $X$ on $(M,\mu)$. Let
$$
||X||=\max_{x\in M}\sqrt{\mu(X(x),X(x))}=\sqrt{\mu(X(y),X(y))}.
$$
According to Proposition 5.7 of Chapter VI in \cite{KN}, the curve
$\gamma(t)=\exp(tX)(y)$, $0\leq t\leq 1$, is a segment of a
geodesic in $(M,\mu)$ with the length $||X||$. It is known that
for any other point $x\in M$ the curve $\exp(tX)(x)$, $0\leq t\leq
1$, is parameterized proportionally to the arc-length with the
coefficient of proportionality $\sqrt{\mu(X(x),X(x))}$, which does
not exceed $||X||$. Therefore, the length of any arc of the second
curve does not exceed the length of the corresponding arc of the
geodesic $\gamma$.

The injectivity radius of the compact smooth manifold $(M,\mu)$ is
bounded below by some number $r>0.$ If $0\leq s||X||\leq r$;
$t,s\in [0,1]$, then it implies that for $g(s)=\exp(sX)$,
$g(t)=\exp(tX)$,  the point $\gamma(t)$ is the point of maximal
displacement on $(M,\rho)$ for the motion $g(s)$, since
$\rho(g(s)(\gamma(t)),\gamma(t))=s||X||$ according to equalities
$$
g(s)(\gamma(t))=g(s)(g(t)(y))=g(s+t)(y)=\gamma(s+t).
$$
Hence, $d(g(t),g(t+s))=s||X||$, the length of the curve $g(t)$,
$0\leq t\leq 1$, in $(G,d)$ equals to $||X||$. Therefore, one can
join any two point in $(G,d)$ by a curve of finite length (with
respect to the metric $d$). Let $D$ be the inner metric
corresponding to $d$.

There exists a positive number $s_0$ such that $\exp: g\rightarrow
G$ is a homeomorphism of some open subset $V$ of $g$, containing
the zero, onto the open ball $U(e,s_0)$ with the radius $s_0$ in
$(G,d)$. Then the above reasonings imply that the curve $g(t)$,
$0\leq t\leq 1$, is a geodesic in $(G,D)$, and $D(g,h)=d(g,h)$, if
$d(g,h)< \min(r,s)$. Also, $d\leq D$.

From the above calculations of the length of the geodesic
$g(t)=\exp(tX)$, $0\leq t\leq 1$, in $(G,D)$, it is clear that $D$
is the bi-invariant Finsler (inner) metric on $G$ determined by
the $\Ad(G)$-invariant norm $||\cdot||$ on $G_{e}$, which defined
by the formula (\ref{no}). It is easy to check that this formula
defines some norm on $G_{e}$.
\end{proof}

\begin{vopros}
Whether the metrics $d$ and $D$ coincide on $G$?
\end{vopros}

\begin{theorem}\label{killi}
Let $(M,\mu)$ be a compact homogeneous Riemannian manifold. Then
there exists a positive number $s>0$ such that for arbitrary
motion $f$ of the space $(M,g)$ with maximal displacement
$\delta$, which is less than $s$, there is unique Killing vector
field $X$ on $(M,g)$ such that $\max_{x\in
M}\sqrt{\mu(X(x),X(x))}=1$ and $\gamma_{X}(\delta)=f$, where
$\gamma_{X}(t)$, $t\in \mathbb{R}$ is the one-parameter motion
group in $(M,g)$ generated by the field $X$. If also $f$ is a
Clifford-Wolf translation, then the Killing field $X$ has constant
unit length on $(M,\mu)$.
\end{theorem}

\begin{proof}
Let us supply the identity component $G$ of the full isometry
group of $(M,g)$ with the bi-invariant metric $d$ as in Theorem
\ref{D}. There is sufficiently small number $s>0$ (which we can
suppose smaller than the injectivity radius $r$ of the manifold
$(M,\mu)$) such that the exponential map $\exp: g\rightarrow G$ is
a homeomorphism of some neighborhood $V$ of the zero in $g$ onto
an open ball $U(e,s)$ in $(G,d)$. Then for every motion $f$ of the
space $(M,\mu)$ with the condition $d(f,e)=\delta <s$ there exists
the unique vector $Y\in V$ such that $\exp(Y)=f$. It was shown in
the proof of Theorem \ref{D} that for all such motions $f$ we have
$D(f,e)=d(f,e)$. This common value is equal also to the length of
the path $\exp(\tau Y)$, $0\leq \tau\leq 1$, which joins elements
$e$ and $f$, with respect to the bi-invariant norm $||\cdot||$ on
$TG$ from Theorem \ref{D}, and to the length $||Y||$. By the
definition, $||Y||=\max_{x\in M}\sqrt{\mu(Y(x),Y(x))}$. Now it is
clear that $X=(1/{\delta}) Y$ is an desired vector. The uniqueness
of $X$ follows from the above arguments.

Let us suppose also that $f$ is a Clifford-Wolf translation. By
the above construction we have
\begin{equation}\label{1}
||X||=1=\max_{x\in
M}\sqrt{\mu(X(x),X(x))}=\sqrt{\mu(X(x_{1}),X(x_{1}))}
\end{equation}
for some point $x_{1}\in M$. We state that
$$
\sqrt{\mu(X(x),X(x))}\equiv 1.
$$
Indeed, in the opposite case there would be a point $x_{0}\in M$
such that $\sqrt{\mu(X(x_{0}),X(x_{0}))}=\varepsilon< 1$. Then the
path $c(t)=\exp(tX)(x_{0})$, $0\leq t\leq \delta$, joins the point
$x_{0}$ with the point $f(x_{0})$ and has the length $\delta
\varepsilon$.  Therefore,
$$
\rho(x_{0},f(x_{0}))\leq \delta \varepsilon < \delta=
\rho(x_{1},f(x_{1})),
$$
because, according to the condition (\ref{1}), the orbit of the
point $x_{1}$ under the action of the one-parameter group
$\exp(tX)$, $t\in \mathbb{R}$, is a geodesic \cite{KN}, and
$\delta< r$. But this contradicts to the fact that $f$ is a
Clifford-Wolf translation.
\end{proof}

\begin{theorem}\label{subm}
Let $(M,\mu)$ be a compact connected ($G$)-$\delta$-homogeneous
Riemannian manifold with inner metric $\rho$, and let $G$ be a
closed connected (Lie) subgroup of the full isometry group of
$(M,\mu)$, supplied by the bi-invariant inner metric $D$ as in
Theorem \ref{D} (more exactly, by it's restriction to $G$). Then
$D$ is an inner bi-invariant metric on $G$. Let us fix a point
$x_{0}\in M$ and define a projection $p:G\rightarrow M$ by the
formula $p(g)=g(x_{0})$ such that under usual identification of
$M$ with $G/H$, where $H$ is the stabilizer of $G$ at the point
$x_{0}$, $p$ coincides with the canonical projection
$p:G\rightarrow G/H$. Then the map $p:(G,D)\rightarrow (M,\rho)$
is a submetry.
\end{theorem}

\begin{proof}
The first statement easily follows from arguments in the last two
paragraphs in the proof of Theorem \ref{D}, applied to $G$.

Now it is enough to check the properties 1) and 2) from the proof
of Theorem \ref{delta1}.

1) Let $g,h\in G$. Then
$$\rho(p(g),p(h))=\rho (g(x_{0}),h(x_{0}))\leq
\max_{x\in M}\rho(g(x),h(x))=d(g,h)\leq D(g,h),
$$
i.e. $p$ does not increase distances.

2) Consider any points $x,y$ in $M$ and put $\rho(x,y)=a$. Let us
choose arbitrary shortest $K$ in $(M,\rho)$ joining points $x$ and
$y$; consider a geodesic $\gamma(s)$, $s\in \mathbb{R}$, in
$(M,\mu)$ parameterized by the arc-length such that $\gamma(0)=x$,
$\gamma(a)=y$ and $\gamma(s)\in K$, $0\leq s\leq a$. Since
$(M,\rho)$ is $G$-$\delta$-homogeneous, there is
$\delta(x)$-translation $g_t\in G$ of $(M,\rho)$, moving the point
$x$ to the point $\gamma(t), 0< t\leq a.$ Now if $t$ is small
enough, then by Theorems \ref{D} and \ref{killi}, there is an
one-parameter group of motions $g(s)=\gamma_X(s)\in G$, $s\in
\mathbb{R}$, such that $g(t)=g_t$ and $\max_{y\in
M}\sqrt{\mu(X(y),X(y))}=\sqrt{\mu(X(x),X(x))}.$ Then
$g(s)(x)=\gamma(s), s\in \mathbb{R}.$

Therefore, $D(e=g(0),g(s))=d(e,g(s))=s$ for $0\leq s \leq a$.
Suppose that $p(h)=h(x_{0})=x$ for some element $h\in G$. Then
$$
y=\gamma(a)=g(a)(x)=g(a)(h(x_{0}))=p(g(a)h)
$$
and
$$
D(h,g(a)h)=D(e,g(a))=a=\rho(x,y).
$$
\end{proof}

On the ground of Corollary \ref{NormCol} and Theorem \ref{subm}
we obtain

\begin{corollary}\label{general1}
A compact connected Riemannian manifold is
($G$)-$\delta$-homogeneous if and only if it is ($G$)-normal in
the generalized sense.
\end{corollary}

Let us consider a compact Riemannian homogeneous manifold
$(G/H,\mu)$, some $\Ad(G)$-invariant inner product
$\langle \cdot,\cdot \rangle$ on the Lie algebra $\mathfrak{g}$
of the group $G$, the corresponding
$\langle \cdot,\cdot \rangle$-orthogonal direct sum decomposition
$\mathfrak{g}=\mathfrak{h}\oplus \mathfrak{p}$
($\mathfrak{h}$ is the Lie algebra of $H$),
and $\Ad(H)$-invariant inner product $(\cdot,\cdot)$
on $\mathfrak{p}$ which defines the Riemannian metric $\mu$.
Then we can state the previous corollary as follows:

\begin{theorem}\label{body}
A compact Riemannian manifold $(G/H,\mu)$ is
$G$-$\delta$-homogeneous for Lie group $G$ if and only if there
exists an $\Ad(G)$-invariant centrally symmetric (relative to
zero) convex body $B$ in $\mathfrak{g}$ such that
$$
P(B)=\{v\in \mathfrak{p}\,|\, (v,v)\leq 1\},
$$
where $P:\mathfrak{g}\rightarrow \mathfrak{p}$ is $\langle \cdot,\cdot
\rangle$-orthogonal projection. One can take $C=\{w\in \mathfrak{g}\,|\,
||w||\leq 1\}$ as $B$.
\end{theorem}

\begin{corollary}\label{a}
The vector space $\mathfrak{p}$ and the inner product $(\cdot,\cdot)$ are
invariant under $\Ad(N_G(H_0))$, where $N_G(H_0)$ is the normalizer
of the connected unit component $H_0$ of $H$ in $G$.
\end{corollary}

\begin{proof} Evidently, $\mathfrak{h}$ is $\Ad (N_G(H_0))$-invariant.
Then $\mathfrak{p}$
is also $\Ad (N_G(H_0))$-invariant, because $\langle \cdot,\cdot \rangle$
is $\Ad (G)$-invariant.
Now the $\Ad (N_G(H_0))$-invariance of $(\cdot,\cdot)$ follows from Theorem
\ref{body}.
\end{proof}

\begin{remark}
It follows from Theorems \ref{subm} and \ref{main} that in general
case the metric $D$ on $G$ is not Riemannian even in the case when
$(M,\mu)$ is a $\delta$-homogeneous Riemannian manifold. This is
the reason for the words "in the generalized sense" in the
statement of Theorem \ref{general1}.
\end{remark}

\begin{theorem}\label{go}
A Riemannian manifold $(M,\mu)$ is ($G$)-$\delta$-homogeneous if
and only if any of two following conditions are satisfied:

1) For every tangent vector $v\in M_x,$ where $x$ is any point in
$M,$ there is a Killing vector field $X$ (in the Lie algebra $RG$
of right-invariant vector fields on the Lie group $G$) on $M$ such
that $X(x)=v$ and $\mu(X(x),X(x))=\max_{y\in M}\mu(X(y),X(y)).$

2) Every geodesic $\gamma$ in $M$ is an orbit of a 1-parameter
motion group of $M$ (in $G$) generated by a Killing vector field,
attaining a maximal value of its length on $\gamma.$
\end{theorem}

\begin{proof}
Let us remark at first that we can suggest that the vector $v$ in
the condition 1) is non-zero; then the condition 2) implies
condition 1), while the condition 2) follows from the condition 1)
and Proposition 5.7 of the chapter VI in \cite{KN}, which states
that an integral trajectory of a Killing vector field $X$ on $M$,
going through a point $x\in M,$ is a geodesic, if $x$ is a
critical value of the function $\mu(X,X)$ and $X(x)\neq 0.$

Let suppose that $(M,\mu)$ is $\delta$-homogeneous. Then Theorems
\ref{killi} and \ref{product} immediately imply the condition 2).

Sufficiency of 2). It's clear that the condition 2) implies that
$M$ is ($G$)-homogeneous. Then there is a constant $r>0$ such that
$\Radinj (M)>r.$ Let $x,y \in M$ and $\rho(x,y)=t<r.$ Then there is
unique geodesic $\gamma(s), s\in \mathbb{R},$ parameterized by arc
length such that $\gamma(0)=x, \gamma(t)=y.$ By the condition,
$\gamma(s)=g(s)(x),$ where $g(s), s\in \mathbb{R},$ is a
1-parameter motion group of $M$ (in $G$), generated by a Killing
vector field $X,$ such that $\mu(X(x),X(x))=\max_{z\in
M}\mu(X(z),X(z)).$ Then it is clear that for every $z\in X,$
$\rho(x,y)=\rho(x,g(t)(x))\geq \rho(z,g(t)(z)).$ We proved that
$M$ is restrictively ($G$)-$\delta$-homogeneous. Hence $M$ is
($G$)-$\delta$-homogeneous by Proposition \ref{pr}.
\end{proof}

\begin{definition}
A Riemannian manifold $(M,\mu)$ is called
\textit{($G$)-geodesic orbit} (($G$)-g.o.), if every
geodesic in $M$ is an orbit of a one-parameter isometry subgroup
(in $G$).
\end{definition}

More extensive information on geodesic orbit manifolds
(or {\it geodesic orbit spaces} by another terminology)
one can find e.g. in \cite{AA, KV, tam, tam1, Zi96}.

\begin{corollary}
\label{goc} Every ($G$)-$\delta$-homogeneous Riemannian manifold
is ($G$)-geodesic orbit (($G$)-g.o.) manifold.
\end{corollary}

\section{Totally geodesic submanifolds}
\label{totally}

In this section we investigate some totally geodesic submanifolds
of $\delta$-homogeneous and g.o. Riemannian manifolds.

\begin{pred}[Theorem 8.9 of Chapter VII in \cite{KN}]\label{gon1}
Let $M$ be a Riemannian manifold, $N$ is its totally geodesic
submanifold, $X$ is a Killing field on $M$.
Consider a smooth vector field $\widetilde{X}$ on
$N$, with is tangent (with respect to
$N$) component of the field $X$. Then $\widetilde{X}$
is a Killing field on the Riemannian manifold $N$.
\end{pred}

In \cite{KN} this proposition is used to prove that every closed
totally geodesic submanifold of a homogeneous Riemannian manifold
is homogeneous itself (Corollary 8.10 of Chapter VII in
\cite{KN}). Here we give some refinement of this classical result.

\begin{theorem}\label{gon2}
Every closed totally geodesic submanifold of a
$\delta$-homogeneous Riemannian manifold is
$\delta$-homogeneous itself.
\end{theorem}

\begin{proof}
Let $N$ be a closed totally geodesic submanifold of a
$\delta$-homogeneous Riemannian manifold $M$. Since $M$ is
homogeneous, it is complete. Since $N$ is closed submanifold of
$M$, it is complete too. Let $U\neq 0$ be a tangent vector at some
point $x\in N$. By Theorem \ref{go} to prove the
$\delta$-homogeneity of $N$ it is enough to show that there is a
Killing field $Y$ on $N$, whose value at the point $x$ is $U$, and
the maximal value of the length of $Y$ is attained at the point
$x$.

Since $M$ is $\delta$-homogeneous Riemannian manifold, there is a
Killing field $X$ on $M$ such that its value at the point $x$ is
$U$, and the maximal value of its length is attained at the point
$x$. Now as a required Killing field $Y$ we can take
$\widetilde{X}$, the tangent component of the field $X$ to $N$.
According to Proposition \ref{gon1}, this field is Killing on $N$
and $\widetilde{X}(x)=X(x)$ obviously. Since at the point $x$ the
length of the field $X$ is maximal among all points $y\in M$, then
$x$ is a point of maximal value for the length of the field
$\widetilde{X}$ (the length of the field $\widetilde{X}$ does not
exceed the length of the field $X$ at all points of the manifold
$N$). Theorem is proved.
\end{proof}

\begin{corollary}\label{gon4}
Every closed totally geodesic submanifold
of a normal homogeneous Riemannian manifold is
$\delta$-homogeneous.
\end{corollary}

\begin{remark}
Let $M$ be a Riemannian manifold, $F$ is some set of its
isometries. Then every connected component of the set of points of
$M$, which are fixed under every isometry in $F$, is a closed
totally geodesic submanifold of $M$. By the same manner, if $K$ is
some set of Killing fields on $M$, then every connected component
of the set of points of $M$, which are zeros for every Killing
field in $K$, is a closed totally geodesic submanifold of $M$
\cite{KN}.
\end{remark}

\begin{theorem}\label{gon3}
Every closed totally geodesic submanifold
of geodesic orbit (g.o.) Riemannian manifold
is geodesic orbit itself.
\end{theorem}

\begin{proof}
Let $N$ be a closed totally geodesic submanifold of a geodesic
orbit Riemannian manifold $M$. It is clear that $M$ and $N$ are
complete. Let $U\neq 0$ be a tangent vector at some point $x\in
\widetilde{M}$. It is enough to prove that there is a Killing
field $Y$ on $N$ with the following properties:

1) the value $Y$ at the point $x$ is $U$;

2) $x$ is a critical point of the length of the field $Y$ on
$N$.

Indeed, in this case a geodesic passing through $x$ in the
direction $U$ is an orbit of an one-dimensional motion group
generated by the Killing field $Y$ (this one-parameter group is
correctly defined because of the completeness of $N$).

Since $M$ is a geodesic orbit Riemannian manifold, there is a
Killing field $X$ on $M$, whose value at the point $x$ is $U$, and
such that $x$ is a critical point of the length of the field $X$.
Now as a required Killing field $Y$ one can consider
$\widetilde{X},$ the tangent component of the field $X$ to $N$.
According to Proposition \ref{gon1}, it is a Killing field on $N$,
and, moreover, $\widetilde{X}(x)=X(x)$.

Now we need to prove only that $x$ is a critical point of the
length of the field $\widetilde{X}$ on $N$.
Let
$Z=X-\widetilde{X}$ be the normal component of the field $X$ on
the manifold $N$, and let $g$ be the metric tensor on
$M$. It is clear that
$$
g(\widetilde{X},\widetilde{X})=g(X,X)-g(Z,Z).
$$
The point $x$ is a zero point for $g(Z,Z)$, therefore, $x$ is a
point of the minimal value of $g(Z,Z)$ on $N$. Consequently, $x$
is a critical point both to the function $g(X,X)$ and to the
function $g(Z,Z)$ on the manifold $N$. But in this case $x$ is a
critical point for the function $g(\widetilde{X},\widetilde{X})$
also. Therefore, $x$ is a critical point of the length of the
field $\widetilde{X}$ (since $\widetilde{X}(x)=U\neq 0)$. Theorem
is proved.
\end{proof}

According to Lemma \ref{dp}, the metric product of
$\delta$-homogeneous spaces is
$\delta$-homogeneous itself.
In the Riemannian case we have the conversion to this statement:

\begin{theorem}\label{Ram}
Let $M=M_{0}\times M_{1}\times \dots \times M_{k}$ be a direct
metric decomposition of a $\delta$-homogeneous (respectively,
g.o.) Riemannian manifold $M$ with the maximal Euclidean factor
$M_{0}$. Then all factors of this product are $\delta$-homogeneous
(respectively, g.o.). If $M$ is $\delta$-homogeneous, then $M_{i}$
are compact for $i\neq 0$. Besides, an isometry $f=f_0\times \dots
\times f_k$ of the manifold $M$, which is a product of
$\delta$-translations, is a $\delta$-translation itself.
\end{theorem}

\begin{proof}
Since every fiber of the product under consideration is a complete
totally geodesic submanifold, then according to Theorem \ref{gon2}
(Theorem \ref{gon3}), all factors are $\delta$-homogeneous
(respectively, g.o.), which proves the first statement. The second
statement follows from the maximality of the Euclidean factor
$M_{0}$, Proposition \ref{pos} and Theorem \ref{BerP}. The last
statement of Theorem follows from Lemma \ref{dp}.
\end{proof}

Since every g.o. (in particular, every $\delta$-homogeneous)
Rimannian manifold is homogeneous, it is useful to remind an
algebraic description of homogeneous Riemannian manifolds. Let
$(M,\mu)$ be a homogeneous Riemannian manifold with a closed
connected transitive isometry group $G$, and $H$ is its isotropy
subgroup at a given point $x\in M$. Then $M$ is naturally
identified with the coset space $G/H$. Consider the Lie algebras
$\mathfrak{h}$ and $\mathfrak{g}$, $\mathfrak{h} \subset
\mathfrak{g}$, of the groups $G$ and $H$. It is possible to choose
some $\Ad(H)$-invariant complement $\mathfrak{p}$ to
$\mathfrak{h}$ in $\mathfrak{g}$, which could be identified with
the tangent space $M_x$ of $(M,\mu)$ at the point $x$. In this
case the homogeneous Riemannian metric $\mu$ is identified with
some $\Ad(H)$-invariant inner product $(\cdot,\cdot)$ on
$\mathfrak{p}$, whereas $\mathfrak{g}$ is identified with the Lie
algebra of Killing vector fields on $(M,\mu)$ (see details in
\cite{Bes}, Chapter VII).

\begin{remark}\label{biinv}
If $M$ is compact, then $G$ is compact too, therefore, there exists
some $\Ad(G)$-invariant inner product
$\langle \cdot,\cdot \rangle$ on the Lie algebra $\mathfrak{g}$ of
the group $G$. In this case as $\mathfrak{p}$ we can consider
a $\langle \cdot,\cdot \rangle$-orthogonal complement to $\mathfrak{h}$ in
$\mathfrak{g}$.
Note also that restrictions of
$(\cdot,\cdot)$ and $\langle \cdot,\cdot \rangle$ to any
$\Ad(H)$-invariant and
$\Ad(H)$-irreducible submodule $\mathfrak{q}\subset \mathfrak{p}$
are proportional one to another.
\end{remark}

Let $(M=G/H,\mu)$ be a homogeneous Riemannian manifolds with a
closed connected transitive isometry group $G$, which is generated
by some $\Ad(H)$-invariant inner product $(\cdot,\cdot)$ on
$\mathfrak{p}$ in the above notation. For Killing fields $X,Y\in
\mathfrak{p}$ we have the following equality:
\begin{equation}\label{connec}
\nabla_XY(x)=-\frac{1}{2}[X,Y]_{\mathfrak{p}}+U(X,Y),
\end{equation}
where the (bilinear symmetric) map
$U:\mathfrak{p}\times \mathfrak{p} \rightarrow \mathfrak{p}$ is
defined by the formula
\begin{equation}\label{connec1}
2(U(X,Y),Z)=([Z,X]_{\mathfrak{p}},Y)+(X,[Z,Y]_{\mathfrak{p}})
\end{equation}
for any $Z\in \mathfrak{p}$ \cite{Bes}.
In \cite{Al68} it is proved the following (compare with \cite{Wallach72}, Theorem 4.1)

\begin{pred}[\cite{Al68}]\label{NWAL}
Let $(M=G/H, \mu)$ be any homogeneous Riemannian manifold and $T$
be any torus in $H$, $C(T)$ is its centralizer in $G$. Then
the orbit $M_T=C(T)(x)$ is a totally geodesic submanifold of $(M, \mu)$.
\end{pred}

\begin{proof}
It is easy to get that the Lie subalgebra $\mathfrak{l}$ of $C(T)$ in
$\mathfrak{g}$ has
the form $\mathfrak{l}=\mathfrak{k}\oplus \mathfrak{q}$,
where $\mathfrak{q}=\{X\in \mathfrak{p}\,|\,[X,\mathfrak{t}]=0\}$,
$\mathfrak{k}=\{X\in \mathfrak{h}\,|\,[X,\mathfrak{t}]=0\}$,
$\mathfrak{t}\subset \mathfrak{k}$ is the Lie algebra of $T$.

According to (\ref{connec}), to prove Proposition we need to show
that $U(X,Y)\in \mathfrak{q}$ for any $X,Y\in \mathfrak{q}$.
Let $W\in \mathfrak{p}$, $Z\in \mathfrak{t}$,
then
$$
2([Z,U(X,Y)],W)=-2(U(X,Y),[Z,W])=
-([[Z,W],X]_{\mathfrak{p}},Y)-(X,[[Z,W],Y]_{\mathfrak{p}})=
$$
$$
([[W,X],Z]_{\mathfrak{p}},Y)+(X,[[W,Y],Z]_{\mathfrak{p}})=
([W,X]_{\mathfrak{p}},[Z,Y])+([Z,X],[W,Y]_{\mathfrak{p}})=0.
$$
Since $W\in \mathfrak{p}$ may be chosen arbitrary, we have $[Z,U(X,Y)]=0$
for any $Z\in \mathfrak{t}$. This means that $U(X,Y)\in \mathfrak{q}$.
\end{proof}

\begin{remark}
If $T$ is a maximal torus in $H$, then subalgebra $\mathfrak{k}=\mathfrak{t}$
is a part
of the center of Lie algebra $\mathfrak{l}$. Therefore, in this case
$\mathfrak{q}$ is
the Lie algebra of some subgroup $Q\subset G$. Moreover, we can
consider $M_T$ as an orbit of $Q$ through the point $x\in M$.
\end{remark}

Now we consider some properties of g.o. manifolds. If we represent
a homogeneous Riemannian metric $\mu$ on $M=G/H$ as a suitable
$\Ad(H)$-invariant inner product $(\cdot,\cdot)$ on $\mathfrak{p}$
in the above notation, we can consider a useful notion of {\it
geodesic vectors} on $(M,\mu)$. A vector $X+Y$, where $Y\in
\mathfrak{p}$ and $Y\in \mathfrak{h}$, is called {\it geodesic},
if the orbit of one-parameter group generated by the Killing field
$X+Y$ is a geodesic of $(M,g)$, passing through the point $x\in M$
with stabilizer group $H$ in the direction $X$. It is clear that a
homogeneous Riemannian manifold $(G/H=M,\mu)$ is G-g.o. manifold
if and only if for any $X\in \mathfrak{p}$ there is $Y\in
\mathfrak{h}$ such that the vector $X+Y$ is geodesic. It is well
known the following criterion for geodesic vectors (see e.g.
\cite{KV}).

\begin{pred}\label{geovec}
A vector $X+Y$, where $X\in \mathfrak{p}$ and $Y\in \mathfrak{h}$,
is geodesic if and only if for every $V\in \mathfrak{p}$ the
equality $([X+Y,V]_{\mathfrak{p}},X)=0$ holds.
\end{pred}

\begin{pred}\label{simplen}
Let $(G/H,\mu)$ be a $G$-g.o.-space. For any $X\in \mathfrak{p}$
and $Y\in \mathfrak{h}$ such that $X+Y$ is geodesic vector we have
the equality $U(X,X)=[X,Y]$, where $U$ is defined by
(\ref{connec1}).
\end{pred}

\begin{proof}
For the geodesic vector $X+Y$ we have the equality
$$
0=(X,[V,X+Y]_{\mathfrak{p}})=(X,[V,X]_{\mathfrak{p}})+(X,[V,Y])=
(X,[V,X]_{\mathfrak{p}})+([Y,X],V)=
$$
$$
(U(X,X)+[Y,X],V)
$$
for every $V\in \mathfrak{p}$. Therefore, $U(X,X)=[X,Y]$.
\end{proof}

\begin{definition}\label{natred}
A homogeneous Riemannian manifold $(M,\mu)$ is called {\it
(G)-naturally reductive}, if there exist a connected Lie subgroup
$G\subset \Isom (M)$, acting transitively and effectively on $M$
and a $\Ad(H)$-invariant decomposition
$\mathfrak{g}=\mathfrak{h}\oplus \mathfrak{p}$, where $\mathfrak{h}$ is
the Lie algebra of the isotropy subgroup $H\subset G$ at some
point in $x\in M$, such that one of the following equivalent
statements holds:

(1) every geodesic in $M$ through the point $x$ is an orbit of a
one-parameter subgroup in $G$, generated by some $X\in \mathfrak{p}$;

(2) $\mu([Z,X]_{\mathfrak{p}},Y)+\mu(X,[Z,Y]_{\mathfrak{p}})=0$ for all
$X,Y,Z \in \mathfrak{p}$ (in other words, $U\equiv 0$).
\end{definition}

We obviously get from Proposition \ref{simplen}

\begin{corollary}\label{simple1}
Let $(G/H,\mu)$ be a $G$-g.o.-space. If for any $X\in
\mathfrak{p}$ there is some $Z\in \mathfrak{h}$ such that $X+Z$ is
geodesic vector and $[Z,X]=0$, then $(G/H,\mu)$ is $G$-naturally
reductive.
\end{corollary}

Now we get some simple general remarks.

\begin{pred}\label{simple3}
Let $(G/H,\mu)$ be a $G$-g.o.-space. Consider any
$\Ad(H)$-invariant submodule $\mathfrak{q}\subset \mathfrak{p}$.
Then for every $X,Y\in
\mathfrak{q}$ we have $U(X,Y)\in \mathfrak{q}$.
\end{pred}

\begin{proof}
Consider some geodesic vectors $X+Z_1$, $Y+Z_2$, $X+Y+Z_3$, where
$X,Y \in \mathfrak{q}$ and $Z_i\in \mathfrak{h}$. We get from
Proposition \ref{simplen}  that $U(X,X)=[X,Z_1]\subset
\mathfrak{q}$, $U(Y,Y)=[X,Z_2]\subset \mathfrak{q}$,
$U(X+Y,X+Y)=[X+Y,Z_3]\subset \mathfrak{q}$. Therefore,
$2U(X,Y)=U(X+Y,X+Y)-U(X,X)-U(Y,Y)\in \mathfrak{q}$.
\end{proof}

\begin{pred}\label{gotg}
Let $(G/H, \mu)$ be a $G$-g.o. manifold ($G$-$\delta$-homogeneous
manifold), and $L$ is a Lie subgroup of $G$ such that $H\subset L
\subset G$. Then the orbit of the group $L$ through the point $x$
in $G/H$ is a totally geodesic submanifold of $(G/H, \mu)$. In
particular, $L/H$ with the metric, induced by $\mu$, is g.o. space
(respectively, $\delta$-homogeneous space).
\end{pred}

\begin{proof}
Let $\mathfrak{l}$ be a Lie algebra of $L$. Consider the
decomposition $\mathfrak{l}=\mathfrak{h}\oplus \mathfrak{q}$,
where $\mathfrak{q}=\mathfrak{p}\cap \mathfrak{l}$. Then the
module $\mathfrak{q} \subset \mathfrak{p}$ is $\Ad(H)$-invariant.
According to Proposition \ref{simple3} we have $U(X,Y)\in
\mathfrak{q}$ for every $X,Y\in \mathfrak{q}$. On the other hand,
for every $X,Y\in \mathfrak{q}$ we have $[X,Y]\in
\mathfrak{l}=\mathfrak{h}\oplus \mathfrak{q}$. Therefore, by
(\ref{connec}) we get $\nabla_XY(x)\subset \mathfrak{q}$ for any
$X,Y \in \mathfrak{q}$. This means that the homogeneous
submanifold  $L/H$ (with the induced metric) is totally geodesic
in $(G/H, \mu)$. The last statement follows from Theorem
\ref{gon3} (respectively, \ref{gon2}).
\end{proof}

At the end of this section we note one special property of compact
homogeneous Riemannian manifolds.

\begin{pred}\label{simple2}
Let $(M=G/H,\mu)$ be a compact homogeneous Riemannian manifold. Consider
any $\Ad(H)$-invariant and $\Ad(H)$-irreducible submodule
$\mathfrak{q}\subset \mathfrak{p}$, where $\mathfrak{p}$ is a
$\langle \cdot, \cdot \rangle$-orthogonal
complement to $\mathfrak{h}$, and $\langle \cdot, \cdot \rangle$ is some
$\Ad(G)$-invariant inner product on $\mathfrak{g}$ (see Remark \ref{biinv}).
Then for every $X,Y\in \mathfrak{q}$ we have $U(X,Y)=0$.
\end{pred}

\begin{proof}
Since the module $\mathfrak{q}$ is $\Ad(H)$-invariant and
$\Ad(H)$-irreducible,
$(\cdot,\cdot)|_{\mathfrak{q}}=
\alpha \langle \cdot,\cdot \rangle|_{\mathfrak{q}}$ for some $\alpha >0$.
Therefore, for any $Z\in \mathfrak{p}$ we have
$$
2(U(X,Y),Z)=([Z,X]_{\mathfrak{p}},Y)+(X,[Z,Y]_{\mathfrak{p}})=
\alpha \langle [Z,X]_{\mathfrak{p}},Y
\rangle+ \alpha \langle X,[Z,Y]_{\mathfrak{p}} \rangle= 0,
$$
since $\langle \cdot,\cdot \rangle$ is $\Ad(G)$-invariant.
\end{proof}

\section{Additional symmetries of $\delta$-homogeneous metrics}

Remind that the group $G$ acts on the homogeneous space $G/H$ by
the transformation $L_b:G/H\rightarrow G/H$ ($b\in G$), where
$$
L_b(cH)=bcH.
$$
Let $N_G(H)$ be the normalizer of $H$ in the group $G$. For every
$a\in N_G(H)$ one can correctly define a $G$-equivariant
diffeomorphism $R_a:G/H\rightarrow G/H$ acting by the following
rule:
$$
R_a(cH)=cHa^{-1}=ca^{-1}H.
$$

\begin{theorem}\label{the0.5}
Let $(G/H,\rho)$ be a compact $G$-$\delta$-homogeneous Riemannian
manifold with a connected transitive isometry group $G$, $N_G(H)$
is the normalizer of the subgroup $H$ in the group $G$. Then for
every $a\in N_G(H)$ the diffeomorphism $R_a:G/H\rightarrow G/H$ is
a Clifford-Wolf translation on the Riemannian manifold
$(G/H,\rho)$.
\end{theorem}

\begin{proof}
It is clear that the isometricity of the map $R_a$ is equivalent
to that that for all elements $c\in G$ and $a\in N_G(H)$, the
differential $dr_{a^{-1}}(c)$ preserves the length of every vector
$u\in {\hor}_c\subset G_c$, where ${\hor}_c$ means the horizontal
subspace of the corresponding Riemannian submersion
$pr:(G,\nu)\rightarrow (G/H,\mu)$ in $G_c$ and
$dr_{a^{-1}}(hor_c)=hor_{ca^{-1}}$. Here $r$, $l$ denote the
operations of right and left translations in $G$. We have the
evident equality
$$
r_{a^{-1}}=l_c\circ l_{a^{-1}}\circ (l_a\circ r_{a^{-1}})\circ
l_{c^{-1}},
$$
and the corresponding composition of their differentials. Now it
is clear that $l_{c^{-1}}(c)=e$, $dl_{c^{-1}}(hor_c)={\hor}_e=p$,
and $d(l_a\circ r_{a^{-1}})(e)=\Ad(a)$. But the last map preserves
the space $p$ and the scalar product $(\cdot,\cdot)$ by Corollary
\ref{a} and evident inclusion $N_{G}H\subset N_{G}(H_0)$. All
differentials of left translations preserve the horizontal
distribution and length of horizontal vectors. So, the map $R_a$
is an isometry. It is a Clifford-Wolf translation, because it is
generated by the right translation $r_a$ of $G$, commuting with
all left translations of $G$, which generate a transitive isometry
group of $(G/H,\rho)$.
\end{proof}

\begin{lemma}\label{lem1}
The transformation $R_a$ of (effective) homogeneous space $G/H$
for $a\in N_G(H)$ coincides with a transformation $L_b$ for some
$b\in G$ if and only if $a$ is the product of some central element
of the group $G$ and some element of the group $H$.
\end{lemma}

\begin{proof}
Suppose that $R_a=L_b$ for some $b\in G$. Since $R_a$ obviously
commutes with every transformation $L_d$, $d\in G$, we obtain that
$b$ is in the center of $G$. Further, the condition $R_a=L_b$ is
equivalent to the next one: $ca^{-1}H=bcH=cbH$ for any $c\in G$.
Therefore, $a= \tilde{b}d$, where $\tilde{b}=b^{-1}$ is a central
element of $G$, and $d$ is some element of the group $H$. The
converse is obvious.
\end{proof}

\begin{theorem}\label{the1}
Let $(G/H,\rho)$ be a compact $\delta$-homogeneous Riemannian
manifold with a closed connected transitive semisimple isometry
Lie group $G$. Then the group $N_G(H)/H$ is finite.
\end{theorem}

\begin{proof}
According to Theorem \ref{the0.5}, for every $a\in N_G(H)$ the
diffeomorphism $R_a:=G/H\rightarrow G/H$, acting by the rule
$R_a(cH)=cHa^{-1}=ca^{-1}H$, is an isometry of $(G/H,\rho)$.

If $\dim(N_G(H))>\dim(H)$, then one can choose a continuous family
of isometries of the form $R_a$, which are not in the group $G$.
Really, let us consider a vector $U$, which is in the Lie algebra
of the group $N_G(H)$, but not in $h$. Consider $a=\exp(tU)\in
N_G(H)$ for some real number $t$. Then the transformation $R_a$ is
an isometry of $(G/H,\rho)$. Since the center of the group $G$ is
discrete, with using of Lemma \ref{lem1} we get that for some open
set $O\subset \mathbb{R}$ all the transformations $R_a$ for $a\in
O$ are not in the group $G$. But this contradicts to the fact that
$G$ is the full connected isometry group
 of the Riemannian manifold $(G/H,\rho)$.

Therefore, we conclude that $\dim(N_G(H))=\dim(H)$, and the group
$N_G(H)/H$ is finite, since it is compact.
\end{proof}

\begin{example}\label{DNORM2}
Let $G$ be a connected compact semisimple Lie group, and $\mu$ is
some left-invariant Riemannian metric on $G,$ so $G$ is a closed
connected transitive isometry group of the Riemannian manifold
$(G,\mu)$. Then $(G,\mu)$ is not $G$-$\delta$-homogeneous. Really,
if $(G,\mu)$ is $G$-$\delta$-homogeneous, then according to
Theorem \ref{the1}, the group $N_G(H)/H$ is finite. But in our
case $H=\{e\}$ is trivial, and $N_G(H)/H=G$ is not discrete.
\end{example}

According to the previous example we need to discuss
$\delta$-homogeneous left-invariant metrics on compact Lie groups.
It is clear that any bi-invariant metric $\rho$ on a compact Lie
group $G$ is $G$-$\delta$-homogeneous. But there exist
$\delta$-homogeneous left-invariant metrics on $G$ which are not
bi-invariant. One can show this as follows.

Let $G$ be a compact connected semisimple Lie group, and let $K$
be a connected subgroup of $G$. Among all left-invariant metrics
on $G$ we consider a subclass $\mathcal{M}_{G,K}$ of metrics which
are right-invariant with respect to $K$. It is easy to see that
the subclass $\mathcal{M}_{G,K}$ consists of $(G\times K)$-invariant
metrics on the homogeneous space $M=(G\times K)/\diag(K)$ (we use
the natural inclusion $K\subset G$). Indeed, every metric from
$\mathcal{M}_{G,K}$ has $G\times K$ as a transitive motion group
with the isotropy subgroup $\diag(K)$ at the unit $e\in G$. On the
other hand, it is clear that $G$ is transitive on the space
$M=(G\times K) /\diag(K)$.

Now let us consider a $(G\times K)$-normal homogeneous metric
$\rho$ on $M$. Then the Riemannian homogeneous space $(M,\rho)$ is
$(G\times K)$-$\delta$-homogeneous (Corollary \ref{NormCol}). But
the above discussion implies that $(M,\rho)$ is isometric to the
Lie group $G$ with some left-invariant metric ${\rho}_1$. This
metric could be bi-invariant, but it is easy to see that the set
of $(G\times K)$-normal homogeneous metric $\rho$ on $M$ is more
extensive than the set of bi-invariant metrics on $G$ (for more
details see \cite{DZ}). Therefore, we obtain $\delta$-homogeneous
left-invariant metrics on $G$ which are not bi-invariant.

\begin{example}\label{DNORM2n}
Let $F$ be a connected compact simple Lie group, $G=F^k$, $k\geq
2$, $H=\diag(F)\subset G$. Let us consider the space
$G/H=F^k/\diag(F)$ supplied with the metrics ${\rho}$ generated by
the Killing form of $F^k$. Then the homogeneous Riemannian
manifold $(G/H,{\rho})$ is $\delta$-homogeneous. On the other hand
it is isometric to the Lie group $F^{k-1}$ with some
left-invariant metric ${\rho}_1$. If $k\geq 3$, the metric
${\rho}_1$ is not bi-invariant.
\end{example}

\begin{remark}\label{rem2}
It is obvious that for a compact $G$-$\delta$-homogeneous Riemannian
manifold $(G/H,\rho)$ with positive Euler characteristic all
conditions of Theorem \ref{the1} are fulfilled. Really, any
connected one-dimensional central subgroup of $G$ induces on $G/H$ a
non-vanishing vector field, but this implies that $\chi(G/H)=0$. On
the other hand, in the case of positive Euler characteristic the
statement of Theorem \ref{the1} is well known, since the groups $H$
and $G$ have one and the same rank.
\end{remark}

\section{$\delta$-vectors}\label{general}

Let suppose that $M=(G/H,\mu)$ be a compact homogeneous connected
Riemannian manifold with connected (compact) Lie group $G$. Let
$\mathfrak{g}=\mathfrak{h}\oplus \mathfrak{p}$, $\langle
\cdot,\cdot \rangle$, and $(\cdot,\cdot)$ be the same as in the
section \ref{totally}. We identify the Lie algebra of Killing
vector fields on $M$ with the Lie algebra $\mathfrak{g}$ of right
invariant vector fields on $G$ and use $\Ad(G)$-invariant
(Chebyshev's ) norm $||\cdot||$  on $\mathfrak{g}$ and
corresponding bi-invariant inner metric $D$ on $G$ from Theorem
\ref{D}.

From Section \ref{prostr} we get the following

\begin{pred}\label{sub}
The map $p:(G,D)\rightarrow (G/H,\mu)$ does not increase
distances. It is a submetry if and only if $M$ is
$G$-$\delta$-homogeneous.
\end{pred}

\begin{definition}\label{defdv1}
A vector $w\in \mathfrak{g}$ is called $\delta$-vector on the Riemannian
homogeneous manifold $(M=G/H,\mu)$ if
$|P(w)|:=\sqrt{(P(w),P(w))}=||w||$,
where $P$ is as in Theorem \ref{body}.
(This is equivalent to the
condition that for any $a\in G$, $(w_{\mathfrak{p}},w_{\mathfrak{p}})\geq
(\Ad(a)(w)|_{\mathfrak{p}},\Ad(a)(w)|_{\mathfrak{p}})$.)
\end{definition}

\begin{pred}\label{unique}
Let suppose that for a vector $v\in \mathfrak{p}$, the set $W(v)$
of all $\delta$-vectors of the form $w=v+u$, $u\in \mathfrak{h}$
(such that $||w||=\sqrt{(v,v)}$) is non empty. Then $W(v)$ is
compact and convex. Moreover, there is a unique vector $w=w(v)\in
W(v)$ with the smallest distance $\sqrt{<w-v,w-v>}.$
\end{pred}

\begin{proof}
We can suggest that $\sqrt{(v,v)}=1$. Since $p$ in Proposition
\ref{sub} doesn't increase distances, then $P$ in Theorem
\ref{body} has the same property, and really $||w||=1$. Let suppose that $w_1, w_2\in
W(v)$, $0\leq t\leq 1$, and $w=tw_1+(1-t)w_2$. Then by the
triangle inequality,
$$
||w||=||tw_1+(1-t)w_2||\leq t||w_1||+(1-t)||w_2||=t+(1-t)=1.
$$
Since $P$ is a linear map, then
$$
P(w)=P(tw_1+(1-t)w_2)=tP(w_1)+(1-t)P(w_2)=tv+(1-t)v=v.
$$
One more, because $P$ doesn't increase distances, it follows from
the last two relations that $||w||=1$ and $w\in W(v)$. So, the set
$W(v)$ is convex. Evidently, it is compact, and we proved the
first statement.

It follows from compactness of $W(v)$ the existence of a vector
$w\in W(v)$ with the smallest $|w-v|_1=\sqrt{<w-v,w-v>}$. If we have another such a
vector $w^{\prime} \neq w$, then by the previous statement,
$w^{\prime \prime}:=\frac{1}{2}(w+w^{\prime})\in W(v)$ and
$$
2|w^{\prime \prime}-v|_1=|(w-v)+(w^{\prime}-v)|_1< |(w-v)|_1+|(w^{\prime}-v)|_1=2|w-v|_1,
$$
a contradiction.
\end{proof}

\begin{remark}\label{notation}
Let $v\in \mathfrak{p}$ with $W(v)\neq \emptyset $. According to Proposition \ref{unique},
there is  a unique vector $w \in W(v)$ with the smallest distance $\sqrt{<w-v,w-v>}$.
Later on we shall use
a notation $w(v)$ for this vector and a notation $u(v)$ for the vector
$w(v)-v \in \mathfrak{h}$.
\end{remark}

\begin{pred}\label{four}
Consider any vector $v\in \mathfrak{p}$ with
$W(v)\neq \emptyset$. The following four statements are equivalent
(see Remark \ref{notation}): $w(v)=v$, $u(v)=0$,
$||v||=|v|$, and the corresponding vector field $X(v)$ on $M$ is
infinitesimal $\delta(x_0)$-translation for the point $x_0=p(e)$.
\end{pred}

\begin{pred}\label{fix}
If $W(v)\neq \emptyset,$ then the inequalities $u(v)\neq 0$ and
$||v||>|v|$ are equivalent. In this case the following statements
are satisfied: for every element $g\in G$, such that
$\Ad(g)(\mathfrak{h})=\mathfrak{h}$, the equality $\Ad(g)(v)=v$
(respectively, $\Ad(g)(v)=-v$) implies that $\Ad(g)(u(v))=u(v)$
(respectively, $\Ad(g)(u(v))=-u(v)$).
\end{pred}

\begin{proof}
This follows easily from Propositions \ref{unique}, \ref{sub} and
the fact that $||\cdot||$, ${\langle \cdot,\cdot \rangle}_e$ are
$\Ad(G)$-invariant and invariant under central symmetry.
\end{proof}

From Theorem \ref{body} we get the following

\begin{pred}\label{delram}
A homogeneous Riemannian manifold $(G/H,\mu)$ with connected Lie
group $G$ is $G$-$\delta$-homogeneous if and only if for every
vector $v\in \mathfrak{p}$ there exists a vector $u\in \mathfrak{h}$
such that the vector $v+u$ is a $\delta$-vector.
\end{pred}

\section{On the topology of compact homogeneous spaces}

In general case a Cartan subalgebra $\mathfrak{k}$ of a Lie algebra
$\mathfrak{g}$ is defined as a nilpotent Lie subalgebra in
$\mathfrak{g}$, which coincides with its normalizer in
$\mathfrak{g}$. If a Lie algebra $\mathfrak{g}$ is compact, i.e.
is the Lie algebra of some compact Lie group $G$, then
$\mathfrak{k}$ is a maximal commutative subalgebra in
$\mathfrak{g}$, hence, is the Lie algebra of a maximal torus $T$
in $G$.

\begin{theorem}[\cite{Ad}]
\label{tor} Any two maximal tori in a compact (connected)
Lie group $G$ are conjugate by an inner automorphism of the Lie
group $G.$
\end{theorem}

Thus, the \textit{rank} $\rk(G)$ of a compact Lie group $G$ is
(correctly) defined as the dimension of a Cartan subalgebra $\mathfrak{k}$ in
$\mathfrak{g}$, or, what is equivalent, the dimension of a maximal torus in
$G$.

\begin{theorem}[\cite{HS}, \cite{On}]\label{On}
Let $M=G/H$ be a homogeneous space, where
$G,H$ are connected compact Lie groups. Then $\chi(M)\geq 0.$ The
following statements are equivalent:
(i) $\chi(M)> 0$;
(ii) $\rk(G)=\rk(H)$.
If $\chi(M)>0$, then the manifold $M$ is formal
and $\chi(M)=\frac{|W_{G}|}{|W_{H}|}$, where $|W_{G}|$
(respectively, $|W_{H}|$) is the order of the Weyl group $W_{G}$
(respectively, $W_{H}$)of the Lie group $G$ (respectively, $H$).
\end{theorem}

\begin{theorem}\label{kill}
Let $M=(G/H,\mu)$ be a compact simply connected
homogeneous Riemannian  manifold with compact Lie groups $G$ and
$H$, and $G$ is connected. Then the following conditions are
equivalent:

1) $\chi(M)=0$;

2) $\rk G > \rk H$;

3) There is a right-invariant vector field on $G$, projecting under
canonical map $p:G \rightarrow M$ to nowhere vanishing Killing
field on $M$;

4) All characteristic numbers of the Riemannian  manifold $M$,
defined for principal bundle $\pi: SO(M)\rightarrow M$ of
orthonormal  oriented bases on $M$, are equal to zero.
\end{theorem}

\begin{proof}
It follows from homotopic sequence of the bundle $p:G \rightarrow
G/H$, connectedness of $G$, and simply connectedness of $G/H$ that
the group $H$ is connected. So, all conditions of Theorem \ref{On}
are satisfied. Then the conditions 1) are 2) equivalent.

It is clear that the condition 3) implies the condition 1).

We will show that the condition 2) implies the statement 3).

Let us consider $U\in \mathfrak{g}$ such that the dimension of the closure in
$G$ of one-parameter group $\exp(tU)$ coincides with $\rk(G)$,
which, in turn, is strongly greater than $\rk(H)$. We state that
$\Ad(s)(U) \not \in \mathfrak{h}$ for all $s\in G$. Actually, let suppose
that $V:=\Ad(s)(U) \in \mathfrak{h}$. Since $\Ad(s)$ is an inner automorphism
of Lie algebra $\mathfrak{g}$, then the dimension of the closure in $G$ of
one-parameter group $\exp(tV)$ also coincides with $\rk(G)$. On the
other hand, this closure is a torus in $H$, because $H$ is a
closed Lie subgroup of Lie group $G$. This contradicts to the
inequality $\rk(H)< \rk(G)$. Clearly, a right-invariant vector field
$W$ on $G$ with the condition $W(e)=U$ projects under the map $p$
to a Killing vector field on $M$ without zeroes.

Since any two maximal tori in a compact Lie group are conjugate,
then one can easily prove that the condition 3) implies the
condition 2), because the equality $\rk(G)= \rk(H)$ implies that
every maximal torus is conjugate by an inner automorphism of Lie
group $G$ to a subgroup in $H$. Thus every right-invariant vector
field on $G$ projects to a Killing field on $M$, which necessarily
vanishes at some points.

Characteristic numbers from the condition 4) are defined only for
even-dimensional Riemannian manifold $M$. In this case also Euler
characteristic is a characteristic number (corresponding to the
characteristic Euler class) by Gauss-Bonnet theorem. Then in this
case the condition 1) follows from the condition 4); The statement
4) follows from the condition 3) (even from the more weaker
existence condition of nowhere vanishing Killing vector field on
arbitrary compact smooth oriented Riemannian manifold of even
dimension) by Bott's theorem \cite{Bot} (a proof is also given in
Theorem 6.1 of Chapter 2 in \cite{K}).

In odd-dimensional case $\chi(M)=0$ and the condition 1) is
satisfied, hence 2) and 3), as we said before. If we suggest that
characteristic numbers of odd-dimensional (compact Riemannian)
manifold are equal zero by definition, then the condition 4) is
automatically satisfied. Thus, in this case all 4 conditions are
equivalent and always satisfied.
\end{proof}

\begin{pred}[\cite{Wallach72}]\label{w}
Every even-dimensional homogeneous
Riemannian manifold $M$ of positive sectional curvature has
positive Euler characteristic.
\end{pred}

\begin{proof}
According to Berger's theorem \cite{Berger66}, any Killing field
on an even-dimensional Riemannian manifold of positive sectional
curvature must vanish at some point. If $M=G/H$ would have zero
Euler characteristic, then by Theorem \ref{kill} $M$ would admit
nowhere vanishing Killing vector field. Thus $\chi(M)>0$ by
Hopf-Samelson theorem.
\end{proof}

\begin{remark}
Example of a flat even-dimensional torus, which has zero Euler
characteristic, shows that the statement of Proposition \ref{w} is
not true under the condition of non-negativeness of sectional
curvature. Notice that by Poincar\'e duality, any compact
odd-dimensional triangulated (in particular, smooth) manifold has
zero Euler characteristic.
\end{remark}

\begin{corollary}
All CROSS'es, besides odd-dimensional, i.e. besides $S^{2k+1}$ and
$\mathbb{R}P^{2k+1}$, have positive Euler characteristic.
\end{corollary}

\begin{theorem}
\label{simple} Any simply connected compact homogeneous Riemannian
manifold $(M,g)$ admits a semi-simple compact transitive isometry
group. If moreover the connected component of the group of all
isometries of the space $(M,g)$ is not semi-simple, then
$\chi(M)=0$ and $(M,g)$ is a total space of a Riemannian
submersion, which is a non-trivial principal bundle with simply
connected homogeneous Riemannian base $(M_{1},g_{1})$ and
pair-wise isometric totally geodesic flat tori as fibers. Under
this the connected component of the group of all motions of the
space $(M_{1},g_{1})$ is semi-simple. If $(M,g)$ is
$\delta$-homogeneous, then $(M_{1},g_{1})$ is also
$\delta$-homogeneous.
\end{theorem}

\begin{proof}
The proof follows the line of the paper \cite{B2}. The first
statement of theorem we get on the ground of Corollary 4 of the
section 3 in the chapter 2 in \cite{GO}.

Under this the connected component $G$ of the full isometry group
of the space $(M,g)$ is not semi-simple if and only if $G$ has
non-trivial connected component $C$ of it's center. Then the group
$C$ acts as a non-trivial connected group of Clifford-Wolf
translations on $(M,g)$. Thus $\chi(M)=0$.

It is clear that the orbits of one-parameter subgroups of the
group $C$ in $(M,g)$ are geodesic (see also \cite{B2}). Thus the
orbits of the group $C$ are pair-wise isometric flat totally
geodesic tori in $(M,g)$.

The simply connectedness of $M$ and connectedness of fibers of
Riemannian submersion $p:(M,g)\rightarrow (M_{1},g_{1})$ imply the
non-triviality of the bundle $p$ and simply connectedness of the
space $M_{1}$.

On the ground of Theorem \ref{delta1}, the metric quotient (orbit)
space $(C \backslash M,g_{1})$ $:=(M_{1},g_{1})$ is
$\delta$-homogeneous Riemannian manifold, if $(M,g)$ is
$\delta$-homogeneous Riemannian manifold.
\end{proof}

\begin{remark}
If $(M,g)$ is a homogeneous compact Riemannian manifold and
$\chi(M)> 0,$ then by Theorem \ref{simple}, the connected
component (of effective) full isometry group of the manifold
$(M,g)$ is semi-simple. The opposite statement is not true: the
connected component of full isometry group of Euclidean sphere
$S^{2l-1},l\geq 3,$ is simple Lie group $SO(2l)$ and semi-simple
Lie group $SO(4)$ with Lie algebra $so(4)=so(3)\oplus so(3)$ in
the case of the sphere $S^3.$
\end{remark}

\begin{remark}
The well-known example of Berger spheres $S^{2n+1}=U(n+1)/U(n)$
shows that in general case the connected component $G$ of the unit
for full isometry Lie group of the space $(M,\mu)$ is not
semi-simple, (even if $(M,g)$ is normal); in this case the
universal covering Lie group of $G$ is non-compact. One needs to
note also that for Berger spheres $U(n+1)/U(n)$ (with normal
metrics) the Lie algebra of isotropy group $U(n)$ is not
orthogonal to the center of Lie algebra $u(n+1)$ with respect to
corresponding $\Ad(U(n+1))$-invariant scalar product.
\end{remark}

It follows from Proposition \ref{pos} that every
$\delta$-homogeneous Riemannian manifold has non-negative
sectional curvature.

\begin{vopros}
Whether every compact $\delta$-homogeneous Riemannian
manifold with a finite fundamental group has positive Ricci
curvature?
\end{vopros}

\begin{pred} Let $\mathfrak{h}$  is a Lie subalgebra of a  Lie algebra
$\mathfrak{g}$
of a connected Lie group $G$ and $N_{\mathfrak{g}}(\mathfrak{h})=\mathfrak{h}$,
where $N_{\mathfrak{g}}(\mathfrak{h})$ is
the normalizer of $\mathfrak{h}$ in $\mathfrak{g}$.
Then $\mathfrak{h}$ is a Lie algebra of a
unique closed connected Lie subgroup $H$ in $G$.
\end{pred}

\begin{proof} Let $H_1=\{g\in G: \Ad(g)(\mathfrak{h})\subset \mathfrak{h}\}$.
Then $H_1$ is
closed subgroup in $G$. Hence its connected component $H$ is
closed. By Cartan theorem, $H$ is a Lie subgroup of $G$.
Evidently, Lie algebra of $H$ is equal to $N_{\mathfrak{g}}(\mathfrak{h})$,
which is by condition is equal to $\mathfrak{h}$,
so $H$ is required Lie subgroup.
\end{proof}

One can easily deduce from this the following statements.

\begin{pred} If $\mathfrak{h}$ is a reductive Lie subalgebra of $\mathfrak{g}$,
containing a maximal commutative subalgebra $\mathfrak{t}$ in $\mathfrak{g}$,
then $N_{\mathfrak{g}}(\mathfrak{h})=\mathfrak{h}$.
\end{pred}

\begin{theorem}\label{pose}
Let $G$ be a simple compact connected Lie group and
$\mathfrak{t}$ be Lie algebra of a maximal torus $T\subset G$. Then every
proper Lie subalgebra $\mathfrak{h}\subset \mathfrak{g}$, such that
$\mathfrak{t}\subset \mathfrak{h}$, is a
Lie algebra of the unique closed connected Lie subgroup $H\subset
G$. Moreover, $G/H$ is a simply connected compact connected
homogeneous space of positive Euler characteristic.
\end{theorem}

\section{Homogeneous spaces of positive Euler characteristic}
\label{HomManPECH}

Here we recall some properties of homogeneous spaces with positive
Euler characteristic.

\begin{theorem}[\cite{Shchet2}]\label{Shchet.1}
If $M$ and $M'$ are homogeneous spaces of connected compact Lie
groups, $\chi(M)>0, \chi(M')>0$ and $M$ is homotopically
equivalent to $M'$ then $M$ and $M'$ are diffeomorphic.
\end{theorem}

Now we outline some structure results about homogeneous spaces of
positive Euler characteristic (see \cite{On}, 19.5). Let $G/H$ be
an almost effective compact homogeneous space of positive Euler
characteristic with connected group $G$. From Theorem \ref{kill}
we know that the center of $G$ is discrete (hence, $G$ is
semi-simple), and that there is a maximal torus $T\subset G$ such
that $T\subset H$. Since the center of $G$ is contained in every
maximal torus of $G$, we get the following

\begin{pred}[\cite{Wang49}]\label{effect}
If a compact connected Lie group $G$ acts effectively on the space
$M=G/H$ of positive Euler characteristic, then the center of $G$ is trivial.
\end{pred}

\begin{theorem}[\cite{Kost57}]\label{KostN}
Let $(G/H,\mu)$ be a simply connected compact almost effective
homogeneous Riemannian manifold of positive Euler characteristic.
Then $(G/H,\mu)$ is indecomposable if and only if $G$ is
simple. In particular, a simple and a non-simple compact Lie
groups can not both act transitively and effectively as a group of
motions on a compact Riemannian manifold $M$ with positive Euler
characteristic.
\end{theorem}

A.~Borel and J.~de~Siebenthal obtained in \cite{BorSieb} the
classification of subgroups with maximal rank of compact Lie
groups (see also Section 8.10 in \cite{W}). This classification
give us a description of compact homogeneous spaces with positive
Euler characteristic. A complete description of homogeneous spaces
of classical Lie groups with positive Euler characteristic have
been obtained also by H.C.~Wang in \cite{Wang49}.

We will concern later with special cases of compact homogeneous
manifolds of positive Euler characteristic, namely, the
(generalized) flag manifolds. They can be described as orbits $M$
of a compact connected Lie group $G$ by the adjoint
representation. In other words, $M=G/H,$ where $H=Z_G(S)$ is the
centralizer of a non-trivial torus $S\subset G;$ the Lie group $H$
is always connected. Under this orbits of regular elements in $g$
are called \textit{(full) flag manifolds}.

The chapter 8 in \cite{Bes} contains the following statements:
simply connected compact homogeneous K\"ahler manifolds are
exactly (generalized) flag manifolds. Any latter manifold
(admitting a canonical K\"ahler-Einstein structure, unique in a
sense) is a rational complex algebraic (hence complex projective)
manifold. In a special case $G=Sp(l)$, the stabilizer sub-groups,
whose center is 1-dimensional, are sub-groups $U(l-m)\times
Sp(m).$ Among the corresponding orbits $M^{Sp(l)}_{l-m},$ the only
ones for which the normal metric is K\"ahler (hence
K\"ahler-symmetric) are $M^{Sp(l)}_{1},$ that is
$\mathbb{C}$$P^{2l-1}=Sp(l)/U(1)\times Sp(l-1),$ and
$M^{Sp(l)}_{l},$ isomorphic to $Sp(l)/U(l),$ which is the manifold
of totally isotropic complex $l$-subspaces of $\mathbb{C}^{2l}.$
The space $M^{SO(2l+1)}_{l}=SO(2l+1)/U(l)$ is the manifold of
complex flags of type $l.$

Using chapter 15 in \cite{On}, we can add more. Any (compact
generalized) flag manifold $M,$ supplied with the above mentioned
canonical K\"ahler-Einstein structure, is isomorphic to
$\mathbb{G}/\mathbb{H},$ where $\mathbb{G}$ is a complex connected
Lie group and $\mathbb{H}$ is a closed complex parabolic Lie
subgroup in $\mathbb{G}.$ We recall that a connected complex Lie
subgroup of $\mathbb{G}$ is called \textit{parabolic}, if it
contains a Borel subgroup of $\mathbb{G}.$ A \textit{Borel
subgroup} in $\mathbb{G}$ is any its maximal connected solvable
complex Lie subgroup. Thus $M$ is a so-called \textit{flag
homogeneous space}. Under this, the corresponding complex
structure on $M$ is induced by complex structure on $\mathbb{G}$.
Any parabolic subgroup of $\mathbb{G}$ contains $\Rad(\mathbb{G})$,
a normal subgroup in $\mathbb{G}$. Hence $M$ is a flag homogeneous
space of semi-simple complex Lie group
$\mathbb{G}_0:=\mathbb{G}/\Rad(\mathbb{G})$. Under this
$M=G_0/H_0,$ where $G_0$ is any compact real form of
$\mathbb{G}_0$ and $H_0=G_0\cap \mathbb{H}_0$ for
$\mathbb{H}_0=\mathbb{H}/\Rad(\mathbb{G}).$

It is proved in Corollary 7.12, p. 301 in \cite{MT} that a maximal
connected Lie subgroup $H$ of maximal rank in a compact connected
Lie group $G$ is a connected component of the normalizer (=of the
centralizer) of some element $g\in G.$ On the ground of this
Corollary and connected results, the Table 5.1 in \cite{MT} is
given of all maximal connected compact subgroups $H$ of maximal
rank (more exactly, their Lie subalgebras) in a compact connected
simple Lie groups $G.$ In particular, $G/H$ is an orbit of the
above mentioned element $g\in G$ with respect to the action of the
group $I(G)$ of all inner automorphisms of the Lie group $G.$ A
(generalized) flag manifolds also can be considered as such
orbits, when $g\in G$ is taken in a diffeomorphic image
$\exp_G(U),$ where $U$ is an open ball with the center $0\in
\mathfrak{g}$ with respect to an $\Ad(G)$-invariant Euclidean
metric on $\mathfrak{g}$.

\begin{theorem}[\cite{Shchet1}]\label{Shchet}
Let $G$ and $G'$ be connected compact Lie groups, $H\subset G$ and
$H'\subset G'$ their connected Lie subgroups of maximal rank,
provided that the natural action of $G$ and $G'$ on $M=G/H$ and
$M'=G'/H'$ are locally effective. Suppose that the graded rings
$H(M,Z)$ and $H(M',Z)$ are isomorphic. Then

(i) If $M=M_1\times \dots \times M_s$  and  $M'=M_1'\times \dots
\times M_t'$ are the canonical decompositions of $M$ and $M'$,
then $s=t$ and $M_k$ is diffeomorphic to $M_k'$ after an
appropriate permutation of the factors.

(ii) If $G$ and $G'$ are simple then either the pairs $(G,H)$ and
$(G',H')$ are locally isomorphic or (up to transposition) they are
locally isomorphic to the pairs of the following list:

$$
G=SU(2n)\,(n\geq 2), H=S(U(1)\times U(2n-1)); \quad
$$
$$
G'=Sp(n), H'=U(1)\cdot Sp(n-1);\quad M=M'=\mathbb{C}P^{2n-1}.
$$

$$
G=SO(7), H=SO(6); \quad G'=G_2, H'=SU(3); \quad M=M'=S^6.
$$

$$
G=SO(7), H=SO(5)\times SO(2); \quad G'=G_2, H'=SU(2)\cdot SO(2);
\quad M=M'=Gr^{+}_{7,2}.
$$

$$
G=SO(2n)\, (n\geq 4), H=U(n);  \quad G'=SO(2n-1), H'=U(n-1); \quad
M=M'=I^0Gr^{C}_{2n,n}.
$$
\end{theorem}

Theorem \ref{Shchet} implies easily the classification of
transitive actions of connected compact Lie groups on simply
connected homogeneous spaces of positive Euler characteristic.

Moreover, from results of \cite{On92} and \cite{Shchet1} we have

\begin{theorem}\label{On-Shch}
Let $(G/H,\mu)$ be a simply connected Riemannian homogeneous
manifold of positive Euler characteristic, and $G$ is a simple
connected Lie group. Then the full connected isometry group of
$(G/H,\mu)$ is $G/C$ ($C$ is the center of $G$), excepting the
cases when $(G/H, \mu)$ is one of the following manifolds:

1) $G/H=Sp(n)/U(1)\cdot Sp(n-1)$ ($n\geq 2$), $\mu$ -- symmetric
(Fubini) metric on $\mathbb{C}P^{2n-1}=SU(2n)/S(U(1)\times
U(2n-1))$;

2) $G/H=SO(2n-1)/U(n-1)$ ($n\geq 4$), $\mu$ -- symmetric metric on
$I^0Gr^{C}_{2n,n}=SO(2n)/U(n)$;

3) $G/H=G_2/SU(2)\cdot SO(2)$, $\mu$ -- symmetric metric on
$Gr^{+}_{7,2}=SO(7)/SO(5)\times SO(2)$;

4) $G/H=G_2/SU(3)$ (strongly isotropy irreducible), $\mu$ --
arbitrary $G$-invariant metric.

In the first three cases the metric $\mu$ is not $G$-normal, in the last case
$\mu$ is metric of constant curvature on $S^6=SO(7)/SO(6)$.
\end{theorem}

\begin{proof}
Using Proposition \ref{effect} and Theorem \ref{Shchet}, we easily
get the main statements. We need only to show that in Cases 1),
2), and 3) the metric $\mu$ is not $G$-normal. It follows from
results of \cite{On92}. Really, in that paper the author proved
that the full connected isometry group of a simply connected
$G$-normal homogeneous space $M = G/H$ of a connected simple
compact Lie group $G$, is $G\cdot {\Aut}_G(M)^0$ (a locally direct
product), where
$$
{\Aut}_G(M) = \{f\in \Diff (M)\,|\, f(gx) = gf(x), g \in G, x \in M
\},
$$
excepting the following cases: $G_2/SU(3)=S^6$, $Spin(7)/G_2=S^7$,
$Spin(8)/G_2=S^7\times S^7$. Only one of these spaces (namely,
$G_2/SU(3)=S^6$) has positive Euler characteristic. Moreover, it
is strongly isotropy irreducible. We need to note also that
${\Aut}_G(M)^0$ is trivial for spaces $M = G/H$ of positive Euler
characteristic (it is easy to see from Theorem \ref{KostN}).
\end{proof}

Now we describe the sets of $G$-invariant metrics on the spaces
$G/H$ from items 1), 2), 3) of Theorem \ref{On-Shch}. Note, that
each of these spaces is a (generalized) flag manifold. Note also,
that $G$-invariant metrics on the space $G/H=G_2/SU(3)$
constitutes a one-dimensional family of pairwise homothetic
metrics.

\begin{example}\label{1pec}
It is known (see e.g. \cite{Zi1}) that the set of $G$-invariant
metrics on $G/H=Sp(n)/U(1)\cdot Sp(n-1)$ ($n\geq 2$) is
two-parametric. More exactly, let $\langle \cdot ,\cdot \rangle$
be an $\Ad(Sp(n))$-invariant inner product on the Lie algebra
$\mathfrak{g}=sp(n)$. In this case
$\mathfrak{h}=u(1)\oplus sp(n-1)  \subset \mathfrak{k}:=
sp(1)\oplus sp(n-1)\subset \mathfrak{g}$. Let us consider an $\langle \cdot
,\cdot \rangle$-orthogonal decomposition
$$
sp(n)=\mathfrak{g}=\mathfrak{h}\oplus \mathfrak{p}= \mathfrak{h}\oplus
\mathfrak{p}_1 \oplus \mathfrak{p}_2,
$$
where $\mathfrak{h}\oplus \mathfrak{p}_2=\mathfrak{k}=sp(1)\oplus sp(n-1)$.
Then the modules $\mathfrak{p}_1$
and $\mathfrak{p}_2$ are $\Ad(H)$-invariant, $\Ad(H)$-irreducible, and
pairwise inequivalent with respect to $\Ad(H)$. Therefore, any
$Sp(n)$-invariant metric on $G/H=Sp(n)/U(1)\cdot Sp(n-1)$ is
generated by one of inner products on $\mathfrak{p}$ of the form
$$
(\cdot,\cdot)=x_1 \langle \cdot ,\cdot \rangle|_{\mathfrak{p}_1}+
x_2 \langle \cdot ,\cdot \rangle|_{\mathfrak{p}_2}
$$
for some positive $x_1$ and $x_2$. Note, that the subset of
$SU(2n)$-invariant (symmetric) metrics on $G/H$ consists of the
metrics with the relation $x_2=2x_1$. In this case the full
connected isometry group is a quotient-group of $SU(2n)$ by its
center, and the metric $\mu$ is $SU(2n)$-normal, and
$(Sp(n)/U(1)\cdot Sp(n-1),\mu)$ is isometric to the complex
projective space $\mathbb{C}P^{2n-1}=SU(2n)/U(1)\cdot S(U(2n-1))$)
with the Fubini metric. Note also, that any $Sp(n)$-invariant
metric on $Sp(n)/U(1)\cdot Sp(n-1)$ is weakly symmetric and,
hence, g.o.-metric \cite{Zi96}.
\end{example}

\begin{example}\label{2pec}
The set of $G$-invariant metrics on $G/H=SO(2n-1)/U(n-1)$ ($n\geq
3$) is two-parametric also. More exactly, let $\langle \cdot
,\cdot \rangle$ be an $\Ad(SO(2n-1))$-invariant inner product on
the Lie algebra $\mathfrak{g}=so(2n-1)$. In this case
$\mathfrak{h}=u(n-1) \subset
\mathfrak{k}:=so(2n-2)\subset g=so(2n-1)$. Let us consider an $\langle \cdot
,\cdot \rangle$-orthogonal decomposition
$$
so(2n-1)=\mathfrak{g}=\mathfrak{h}\oplus \mathfrak{p}= \mathfrak{h}
\oplus \mathfrak{p}_1 \oplus \mathfrak{p}_2,
$$
where $\mathfrak{h}\oplus \mathfrak{p}_2=\mathfrak{k}=so(2n-2)$.
Then the modules $\mathfrak{p}_1$ and $\mathfrak{p}_2$ are
$\Ad(H)$-invariant, $\Ad(H)$-irreducible, and pairwise
inequivalent with respect to $\Ad(H)$. Therefore, any
$SO(2n-1)$-invariant metric on $G/H=SO(2n-1)/U(n-1)$ is generated
by one of inner products on $\mathfrak{p}$ of the form
$$
(\cdot,\cdot)=x_1 \langle \cdot ,\cdot \rangle|_{\mathfrak{p}_1}+
x_2 \langle \cdot ,\cdot \rangle|_{\mathfrak{p}_2}
$$
for some $x_1>0$ and $x_2>0$. Note, that the subset of
$SO(2n)$-invariant (symmetric) metrics on $G/H$ consists of the
metrics with the relation $x_2=2x_1$ \cite{Ker}. As in the
previous case, every $SO(2n-1)$-invariant metric on
$SO(2n-1)/U(n-1)$ is weakly symmetric and, hence, g.o.-metric
\cite{Zi96}. Note also that $SO(5)/U(2)$ coincides with
$Sp(2)/U(1)\cdot Sp(1)$ as a homogeneous space.
\end{example}

\begin{example}\label{3pec}
Let us consider now the space $G/H=G_2/SU(2)\cdot SO(2)$, where
$H=SU(2)\cdot SO(2)\subset SU(3)$, and $G_2/SU(3)$ is strongly
isotropy irreducible ($G/H=SO(7)/SO(5)\times SO(2)=Gr^{+}_{7,2}$).
It is easy to see  that there is a subgroup $SO(4)\subset G_2$
such that $SU(2)\cdot SO(2)=SU(3)\cap SO(4)$. Therefore, we have
$\langle \cdot,\cdot \rangle$-orthogonal decomposition
$$
g_2=\mathfrak{h}\oplus \mathfrak{p}=\mathfrak{h}\oplus \mathfrak{p}_1\oplus
\mathfrak{p}_2 \oplus \mathfrak{p}_3,
$$
where $\langle \cdot,\cdot \rangle$ is some $\Ad(G_2)$-invariant
inner product on $g_2$, $su(3)=\mathfrak{h}\oplus \mathfrak{p}_3$,
$so(4)=\mathfrak{h}\oplus \mathfrak{p}_2$,
$\dim(\mathfrak{p}_2)=2$,
$\dim(\mathfrak{p}_1)=\dim(\mathfrak{p}_3)=4$, and every module
$\mathfrak{p}_i$ is $\Ad(G_2)$-invariant and
$\Ad(G_2)$-irreducible. Moreover, the modules $\mathfrak{p}_1$,
$\mathfrak{p}_2$, and $\mathfrak{p}_3$ are pairwise inequivalent
with respect to $\Ad(H)$ \cite{Ker}. Therefore, we have
$3$-parametric family of $G_2$-invariant metrics on $G/H$, every
of each is generated by some inner product
$$
(\cdot,\cdot)= x_1 \langle \cdot,\cdot \rangle|_{\mathfrak{p}_1}+
x_2\langle \cdot,\cdot \rangle|_{\mathfrak{p}_2}
+x_3\langle \cdot,\cdot \rangle|_{\mathfrak{p}_3}
$$
on $\mathfrak{p}$ for some positive $x_i$, $i=1,2,3$.
From \cite{Ker} we know that $SO(7)$-invariant (symmetric) metrics on $G/H$
are exactly metrics with the following relations:
$$
x_2=2x_1,\quad x_3=3x_1.
$$
Remind, that $G_2/SU(2)\cdot SO(2)$ is a flag manifold. The results of
the paper \cite{AA} implies that any $G_2$-invariant g.o.-metric $\mu$
on the space $G_2/SU(2)\cdot SO(2)$ is either $G_2$-normal or $SO(7)$-normal
(symmetric).
\end{example}

Now we shall give a simple description of naturally reductive
homogeneous manifolds of positive Euler characteristic, which
follows from Theorem \ref{KostN}.

\begin{theorem}\label{Kon57}
Let $M$ be a compact naturally reductive homogeneous Riemannian
manifold of positive Euler characteristic. Then $M$ is
$G_1$-normal homogeneous for some (transitive on $M$) semi-simple
Lie subgroup $G_1 \subset G$, where $G$ is the full connected
isometry group of $M$.
\end{theorem}

\begin{proof}
The group $G$ is semisimple, since $\chi(M)>0$ (see Theorem
\ref{simple}).

In the proof of the statement of Theorem  we can  assume without
loss of generality that $M$ is simply connected. Really, the
universal Riemannian covering $\widetilde{M}$ of $M$ has a
semisimple transitive group of motion $\widetilde{G}$, which is a
covering of $G$. Since $G$ and $\widetilde{G}$ have one and the
same Lie algebra, $\widetilde{G}$ is compact, therefore,
$\widetilde{M}$ is compact too. If $\widetilde{M}$ is normal homogeneous with
respect to some semisimple subgroup $\widetilde{G}_1 \subset \widetilde{G}$,
then $M$ is $G_1$-normal homogeneous,
where $G_1\subset G$ is the image of $\widetilde{G}_1$ under the
natural covering epimorphism $\pi:\widetilde{G} \rightarrow G$.

Moreover, we can assume in addition
that $M$ is indecomposable. Really, if $M=M_1\times \cdots \times M_s$ is
the de~Rham decomposition of $M$ then every $M_i$ is naturally reductive
homogeneous manifold
(\cite{Kost56}, Corollary 7; see also \cite{KN}, Chapter X, theorem 5.2).
If we prove that every $M_i$ is normal homogeneous (with respect to
some transitive subgroup of its full connected isometry group),
then $M$ is normal homogeneous too.

Let $M$ be a compact simply connected indecomposable naturally
reductive homogeneous manifold with $\chi(M)> 0$, and $G$ is its
(semisimple) connected isometry group.  From Kostant theorem
(Theorem 4 in \cite{Kost56}) we get that there is a subgroup $G_1
\subset G,$ transitive on $M$, with the following property: there
is an $\Ad(G_1)$-invariant non-degenerate quadratic form $Q$ on
the Lie algebra $\mathfrak{g}_1$ of the group $G_1$ such that the
Riemannian metric of $M$ is generated by by the restriction of $Q$
to $Q$-orthogonal compliment $\mathfrak{p}$ to $\mathfrak{h}_1$ in
$\mathfrak{g}_1$ ($H_1$ is the stabilizer group of some point of
$M$ with respect to the action of $G_1$, and $\mathfrak{h}_1$ is
the corresponding subalgebra of $\mathfrak{g}_1$).

Note that the group $G_1$ is simple according to Theorem
\ref{KostN}. But since $G_1$ is simple, $Q$ is a multiple of the
Cartan-Killing form of $\mathfrak{g}_1$, therefore, $Q$ is
positive definite on $\mathfrak{g}_1$, and $M$ is $G_1$-normal.
Theorem is proved.
\end{proof}

We obviously get from Theorem \ref{Kon57} and Corollary
\ref{NormCol}

\begin{corollary}\label{Kon57.1}
Every compact naturally reductive homogeneous Riemannian manifold
with positive Euler characteristic is $\delta$-homogeneous.
\end{corollary}

According to Corollary \ref{Kon57.1}, a compact naturally
reductive homogeneous Riemannian manifolds $M$, which is not
$\delta$-homogeneous, satisfies the condition $\chi(M)=0$. In Section 15
we obtain examples of $\delta$-homogeneous Riemannian manifolds with
positive Euler characteristic, which are not normal homogeneous
(consequently, are not naturally reductive).

\section{On algebraic corollaries of the $\delta$-homogeneity}

Let $(G/H,\mu)$ be a compact $G$-$\delta$-homogeneous Riemannian
manifold with a connected Lie group $G$, and let $\langle
\cdot,\cdot \rangle$ be an $\Ad(G)$-invariant inner product on the
Lie algebra $\mathfrak{g}$ of the group $G$. Denote by
$\mathfrak{h}$ the Lie algebra of the group $H$, and consider some
$\Ad(H)$-invariant complement $\mathfrak{p}$ to $\mathfrak{h}$ in
$\mathfrak{g}$ (e.g., we can take $\mathfrak{p}$ from the $\langle
\cdot,\cdot \rangle$-orthogonal decomposition
$\mathfrak{g}=\mathfrak{h}\oplus \mathfrak{p}$). It is well know
that the metric $\mu$ is generated by some $\Ad(H)$-invariant
inner product $(\cdot,\cdot)$ on $\mathfrak{p}$, and there is the
equality
\begin{equation}\label{rav1}
(\cdot,\cdot)=x_1 \langle \cdot,\cdot \rangle|_{\mathfrak{p}_1}+ x_2 \langle
\cdot,\cdot \rangle|_{\mathfrak{p}_2}+ \cdots + x_s \langle \cdot,\cdot
\rangle|_{\mathfrak{p}_s}
\end{equation}
for some $\Ad(H)$-invariant pairwise orthogonal (with respect to
both inner products) submodules $\mathfrak{p}_i$ ($1\leq i \leq s$) of the
$\Ad(H)$-module $\mathfrak{p}$ and for some positive numbers $x_i$ ($1\leq i
\leq s$) such that $x_1<x_2<\cdots < x_s$. Note that the modules
$\mathfrak{p}_i$ need not to be $\Ad(H)$-irreducible.

For a vector
$Z\in \mathfrak{g}$ let us denote by $Z_{\mathfrak{p}}$ and $Z_{\mathfrak{h}}$
its projections to subspaces $\mathfrak{p}$ and $\mathfrak{h}$ respectively,
and for a
vector $U\in \mathfrak{p}$ we will denote by $U_i$ its projection to
$\mathfrak{p}_i$,
$1\leq i \leq s$. The symbol $|\cdot|$ denotes the norm on $\mathfrak{p}$,
generated by the scalar product $(\cdot,\cdot)$.

We will give at first another simple proof of the fact that every
($G$)-normal homogeneous Riemannian manifold $(G/H,\mu)$ is
($G$)-$\delta$-homogeneous. Let us consider for this the
decomposition (\ref{rav1}), where $s=1$ and $x_1=1$. Choose any
$X\in \mathfrak{p}$ and show that the vector $X$ is
$\delta$-vector, see Definition \ref{defdv1}. Let $a\in G$, then
by $\Ad(G)$-invariance of the scalar product $\langle \cdot,\cdot
\rangle$ we get
$$
\langle \Ad(a)(X), \Ad(a)(X)\rangle =\langle X, X\rangle,
$$
thus
$$
|\Ad(a)(X)_{\mathfrak{p}}|^2= \langle \Ad(a)(X)_{\mathfrak{p}},
\Ad(a)(X)_{\mathfrak{p}}\rangle \leq
\langle \Ad(a)(X), \Ad(a)(X)\rangle =\langle X, X\rangle =|X|^2.
$$
Proposition \ref{delram} implies that $(G/H,\mu)$ is
($G$)-$\delta$-homogeneous.

Now we derive some corollaries from $\delta$-homogeneity of
Riemannian manifolds in terms of Lie algebras.

Let us consider in a $G$-$\delta$-homogeneous Riemannian manifold
$(G/H,\mu)$ (with a closed connected transitive isometry group
$G$) a geodesic $\gamma$, passing through the point $eH$ in the
direction $V$, $V\in \mathfrak{p}-\{0\}$. Suppose, that the
Killing field $V+U$, $U\in \mathfrak{h}$ admits the maximum of its
length on $\gamma$, and that this field generates an one-parameter
motion group, one of whose orbit is $\gamma$ (Theorem \ref{go}).

\begin{pred}
In the above condition the function $\varphi:G\rightarrow
\mathbb{R},$ defined by the formula
$\varphi(g)=|(\Ad(g)(V+U))_{\mathfrak{p}}|$, where $g \in G$, has the
absolute maximum at the point $g=e$.
\end{pred}

\begin{corollary}
In the above condition one has the following:
\begin{equation}\label{rav2}
(V,[X,V+U]_{\mathfrak{p}})=0 \mbox{ for all } X\in \mathfrak{g},
\end{equation}
\begin{equation}\label{rav3}
(V,[X,[X,V+U]]_{\mathfrak{p}})+|[X,V+U]_{\mathfrak{p}}|^2\leq 0
\mbox{ for all } X\in \mathfrak{g}.
\end{equation}
\end{corollary}

\begin{proof}
Let us consider arbitrary $X\in \mathfrak{g}$. Then the function
$f(t)=|(\Ad(e^{tX})(V+U))_{\mathfrak{p}}|^2$ has its absolute maximum at the
point $t=0$. Now the statement of Corollary follows from the
following:
$$
f(t)=|V|^2+2(V,[X,V+U]_{\mathfrak{p}})t+
\left(|[X,V+U]_{\mathfrak{p}}|^2+(V,[X,[X,V+U]]_{\mathfrak{p}}) \right)
t^2+o(t^2) \mbox{
when } t\rightarrow 0 \,.
$$
\end{proof}

\begin{remark}
Note that for $X\in \mathfrak{h}$ the relations
(\ref{rav2}) and (\ref{rav3}) are fulfilled for any invariant metric.
\end{remark}

\begin{remark}
The equation $(V,[X,V+U]_{\mathfrak{p}})=0$ in the previous
corollary is a well known criterion for geodesic vectors
(see Proposition \ref{geovec}).
\end{remark}

Now we easily obtain

\begin{theorem}\label{ncdo}
Let $(G/H, \mu)$ be a $G$-$\delta$-homogeneous Riemannian manifold
with connected Lie group $G$. Then for every $V\in \mathfrak{p}$
there is $U\in \mathfrak{h}$ such that for every $X\in
\mathfrak{g}$ the following conditions fulfilled:
$$
(V,[X,V+U]_{\mathfrak{p}})=0, \quad
(V,[X,[X,V+U]]_{\mathfrak{p}})+|[X,V+U]_{\mathfrak{p}}|^2\leq 0.
$$
\end{theorem}

\section{On $\delta$-homogeneous manifold of one special type}

Let $G$ be a compact connected Lie group, $H\subset K\subset G$
are its closed subgroup. Fix some $\Ad(G)$-invariant inner product
$\langle \cdot,\cdot \rangle$ on the Lie algebra $\mathfrak{g}$ of the group
$G$. Consider $\langle \cdot,\cdot \rangle$-orthogonal
decomposition
$$
\mathfrak{g}=\mathfrak{h}\oplus \mathfrak{p}=\mathfrak{h}\oplus
\mathfrak{p}_1 \oplus \mathfrak{p}_2,
$$
where $\mathfrak{k}=\mathfrak{h}\oplus \mathfrak{p}_2$ is a Lie algebra of the group $K$.
Obviously,
$[\mathfrak{p}_2,\mathfrak{p}_1]\subset \mathfrak{p}_1$. Let $\mu$
be a $G$-invariant Riemannian metric on $G/H,$ generated by the
inner product
$$
(\cdot,\cdot)=x_1\langle \cdot,\cdot \rangle|_{\mathfrak{p}_1}+
x_2\langle \cdot,\cdot \rangle|_{\mathfrak{p}_2}
$$
on $\mathfrak{p}$ for some $x_1>0, x_2>0$ with $x_1 \neq x_2$.

For any vector $V\in \mathfrak{g}$ we denote by $V_\mathfrak{p}$
and $V_\mathfrak{h}$ its ($\langle \cdot,\cdot
\rangle$-orthogonal) projection to $\mathfrak{h}$ and
$\mathfrak{p}$ respectively.

\begin{pred}[\cite{tam1}]\label{t31.4}
Let $W=X+Y+Z$ be a geodesic vector on $(G/H, \mu)$,
 where $X\in \mathfrak{p}_1$, $Y\in \mathfrak{p}_2$, $Z\in \mathfrak{h}$. Then
we have the following equalities:
\begin{equation}\label{th31.dop1}
[Z,Y]=0,\quad \quad [X,Y]=\frac{x_1}{x_2-x_1}[X,Z].
\end{equation}
\end{pred}

\begin{proof}
By Theorem \ref{ncdo}, for any $U\in \mathfrak{g}$ the equality $(X+Y,[U,
X+Y+Z]_{\mathfrak{p}})=0$ holds. Therefore, we have
$$
(X+Y,[U, X+Y+Z]_{\mathfrak{p}})=x_1\langle X, [U,X+Y+Z] \rangle+
x_2\langle Y, [U,X+Y+Z] \rangle=
$$
$$
x_1\langle [X+Y+Z,X], U \rangle+x_2\langle [X+Y+Z,Y], U \rangle=
$$
$$
\langle (x_2-x_1)[X,Y]+x_1[Z,X]+x_2[Z,Y],U\rangle=0
$$
for any $U\in \mathfrak{g}$. Since $[Z,Y]\in \mathfrak{p}_2$ and
$[X,Y],[Z,X]\in \mathfrak{p}_1$,
this proves Proposition.
\end{proof}

\begin{pred}\label{t31.5}
Let $W=X+Y+Z$ be a $\delta$-vector on $(G/H, \mu)$,
where $X\in \mathfrak{p}_1$, $Y\in \mathfrak{p}_2$, $Z\in \mathfrak{h}$.
Then for any $U\in \mathfrak{p}_1$ the following
inequality holds:
$$
-x_1 \langle [U,X]_\mathfrak{h},[U,X]_\mathfrak{h} \rangle+ (x_2-x_1) \langle
[U,X]_{\mathfrak{p}_2},[U,X]_{\mathfrak{p}_2} \rangle+
(x_1-x_2) \langle [U,Y],[U,X]
\rangle+
$$
\begin{equation}\label{th31.dop2}
(x_1-x_2) \langle [U,Y],[U,Y] \rangle+
x_1 \langle [U,X],[U,Z] \rangle+
\end{equation}
$$
(2x_1-x_2) \langle [U,Y],[U,Z] \rangle+
x_1 \langle [U,Z],[U,Z] \rangle \leq 0.
$$
\end{pred}

\begin{proof}
According to Theorem \ref{ncdo} we get the inequality
$$(X+Y,[U,[U,X+Y+Z]]_{\mathfrak{p}})+([U,X+Y+Z]_
{\mathfrak{p}},[U,X+Y+Z]_{\mathfrak{p}})\leq 0.$$
It is clear that $[Z,X],[Z,U],[Y,X], [Y,U] \in \mathfrak{p}_1$,
$[Z,Y]\in \mathfrak{p}_2$. Therefore, using $\Ad(G)$-invariance of
$\langle \cdot,\cdot \rangle$, we obtain
$$
0\geq (X+Y,[U,[U,X+Y+Z]]_{\mathfrak{p}})+([U,X+Y+Z]_
{\mathfrak{p}},[U,X+Y+Z]_{\mathfrak{p}})=
$$
$$
-x_1\langle [U,X], [U, X+Y+Z]\rangle-x_2\langle [U,Y], [U, X+Y+Z]\rangle+
$$
$$
x_1\langle [U,X]_{\mathfrak{p}_1}+[U,Y+Z], [U,X]_{\mathfrak{p}_1}+[U,Y+Z]\rangle+
x_2\langle [U,X]_{\mathfrak{p}_2}, [U,X]_{\mathfrak{p}_2}\rangle=
$$
$$
-x_1\langle [U,X], [U,X]\rangle
-x_1\langle [U,X], [U,Y]\rangle
-x_1\langle [U,X], [U,Z]\rangle
-x_2\langle [U,Y], [U,X]\rangle
-x_2\langle [U,Y], [U,Y]\rangle
$$
$$
-x_2\langle [U,Y], [U,Z]\rangle
+x_1\langle [U,X]_{\mathfrak{p}_1}, [U,X]_{\mathfrak{p}_1}\rangle
+x_1\langle [U,Y], [U,Y]\rangle
+x_1\langle [U,Z], [U,Z]\rangle
$$
$$
+2x_1\langle [U,Y], [U,X]\rangle
+2x_1\langle [U,X], [U,Z]\rangle
+2x_1\langle [U,Y], [U,Z]\rangle
+x_2\langle [U,X]_{\mathfrak{p}_2}, [U,X]_{\mathfrak{p}_2}\rangle=
$$
$$
-x_1 \langle [U,X]_\mathfrak{h},[U,X]_\mathfrak{h} \rangle+ (x_2-x_1) \langle
[U,X]_{\mathfrak{p}_2},[U,X]_{\mathfrak{p}_2} \rangle+
(x_1-x_2) \langle [U,Y],[U,X]
\rangle+
$$
$$
(x_1-x_2) \langle [U,Y],[U,Y] \rangle+
x_1 \langle [U,X],[U,Z] \rangle+
(2x_1-x_2) \langle [U,Y],[U,Z] \rangle+
x_1 \langle [U,Z],[U,Z] \rangle,
$$
which proves Proposition.
\end{proof}

\begin{corollary}\label{nado1}
If in conditions of Proposition \ref{t31.5} $X=0$, then for any
$U\in \mathfrak{p}_1$ we have
\begin{equation}\label{th31.dop4}
(x_1-x_2) \langle [U,Y],[U,Y] \rangle+
(2x_1-x_2) \langle [U,Y],[U,Z] \rangle+
x_1 \langle [U,Z],[U,Z] \rangle \leq 0.
\end{equation}
\end{corollary}

\begin{pred}\label{t31.6}
For any $\delta$-vector $X+Y+Z$ on $(G/H, \mu)$ the vector $Y+Z$ is a
$\delta$-vector on $K/H$ (with the induced metric). In particular, if $(G/H, \mu)$ is
$G$-$\delta$-homogeneous, then $K/H$ with the induced metric is
$K$-$\delta$-homogeneous.
\end{pred}

\begin{proof}
For any $\Ad(a)$, where $a\in K$,  we have
$\Ad(a)(\mathfrak{p}_1)=\mathfrak{p}_1$.
Moreover, $\Ad(a)|_{\mathfrak{p}_1}$ is orthogonal transformation.
Since
$$
(X,X)+(Y,Y)=(X+Y,X+Y)\geq (\Ad(a)(X+Y+Z)|_{\mathfrak{p}},
\Ad(a)(X+Y+Z)|_{\mathfrak{p}})=
$$
$$
(X,X)+(\Ad(a)(Y+Z)|_{\mathfrak{p}},\Ad(a)(Y+Z)|_{\mathfrak{p}})
$$
for any $a\in K$, the vector $Y+Z$ is $\delta$-vector for $K/H$.
Remark that really the Riemannian subspace $K/H$ of $(G/H, \mu)$
is $K$-normal, because $\mathfrak{k}=\mathfrak{h}\oplus
\mathfrak{p}_2$.
\end{proof}

\begin{pred}\label{t31.7}
For any geodesic vector
$X+Y+Z$ on $(G/H,\mu)$ the vector $Y+Z$ is geodesic vector on
$K/H$ (with the induced metric).
\end{pred}

\begin{proof}
By Proposition \ref{geovec}, $X+Y+Z$ is geodesic if and only if
for any $U\in g$ we have $(X+Y,[U,X+Y+Z]_{\mathfrak{p}})=0$. Let  $U\in
\mathfrak{p}_2\oplus \mathfrak{h}$, then $[U,X+Y+Z]_{\mathfrak{p}_1}=[U,X]$,
$[U,X+Y+Z]_{\mathfrak{p}_2}=[U,Y+Z]_{\mathfrak{p}}$.
Therefore, we have $(Y,[U,Y+Z]_{\mathfrak{p}})=0$,
since $(X,[U,X])=0$. Since $U\in \mathfrak{h}\oplus \mathfrak{p}_2$
may be arbitrary, we
get that the vector $Y+Z$ is a geodesic vector on $K/H$.
\end{proof}

\begin{pred}\label{t31.8}
If vectors $\widetilde{X}+Y+Z$ and
$X+Y+Z$ both are $\delta$-vectors on $(G/H,\mu)$, then
$$
x_1 \langle [\widetilde{X},X]_\mathfrak{h},[\widetilde{X},X]_\mathfrak{h}
\rangle \geq
(x_2-x_1) \langle [\widetilde{X},X]_{\mathfrak{p}_2},
[\widetilde{X},X]_{\mathfrak{p}_2}
\rangle.
$$
\end{pred}

\begin{proof}
From Proposition \ref{t31.4}
we have the equality
$[\widetilde{X},Y]=x_1/(x_2-x_1)[\widetilde{X},Z]$.
Putting $U=\widetilde{X}$ in the inequality
(\ref{th31.dop2}) and using the above equality, we prove Proposition.
\end{proof}

\begin{pred}\label{t31.9}
Suppose that $(G/H, \mu)$ is $G$-$\delta$-homogeneous.
Let $X \in \mathfrak{p}_1$, $Y \in
\mathfrak{p}_2$, $a=\exp(tY)$ for some $t\in \mathbb{R}$,
$\widetilde{X}=\Ad(a)(X)$. Then the following inequality holds:
$$
x_1 \langle [\widetilde{X},X]_\mathfrak{h},[\widetilde{X},X]_\mathfrak{h}
\rangle \geq
(x_2-x_1) \langle [\widetilde{X},X]_{\mathfrak{p}_2},
[\widetilde{X},X]_{\mathfrak{p}_2}
\rangle.
$$
\end{pred}

\begin{proof}
Let $Z\in \mathfrak{h}$ be such a vector that $X+Y+Z$ is $\delta$-vector.
From Proposition \ref{t31.4} we have $[Z,Y]=0$. This implies that
$\Ad(a)(Z)=Z$. Besides this, $\Ad(a)(Y)=Y$, and
$(X,X)=(\widetilde{X},\widetilde{X})$, since $\Ad(a)|_{\mathfrak{p}_1}$ is
$(\cdot,\cdot)$-orthogonal. Therefore, the vector
$\widetilde{X}+Y+Z=\Ad(a)(X+Y+Z)$ is $\delta$-vector too. Now we
can apply Proposition \ref{t31.8}.
\end{proof}

Since for $a=\exp(tY)$ we have $\Ad(a)(X)=X+[Y,X]t+o(t)$ when
$t\to 0$, we get the following infinitesimal version of
Proposition \ref{t31.9}.

\begin{pred}\label{t31.9n}
Suppose that $(G/H, \mu)$ is $G$-$\delta$-homogeneous. Let $X \in
\mathfrak{p}_1$, $Y \in \mathfrak{p}_2$, then the following inequality holds:
$$
x_1 \langle [[Y,X],X]_\mathfrak{h},[[Y,X],X]_\mathfrak{h} \rangle
\geq (x_2-x_1) \langle
[[Y,X],X]_{\mathfrak{p}_2},[[Y,X],X]_{\mathfrak{p}_2} \rangle.
$$
\end{pred}

\section{Root systems of compact simple Lie algebras}
\label{sv}

We give here some information about root systems of a compact
simple Lie algebra $(\mathfrak{g},\langle \cdot,\cdot \rangle=-B)$ with the
Killing form $B$, which can be find in books \cite{Hel,Burb4}.

The Lie algebra $\mathfrak{g}$ admits a direct
$\langle \cdot,\cdot \rangle$-orthogonal decomposition
$\mathfrak{t}\oplus \Lin \{\cup_{\alpha \in \Delta}V_{\alpha}\}$ into (non-zero) vector
subspaces , where $\alpha \in \mathfrak{t}^{\ast}$ is some
(non-zero) real-valued linear form on the Cartan subalgebra
$\mathfrak{t}$ of Lie algebra $\mathfrak{g}$,
$V_{\alpha}=V_{-\alpha}$ is some 2-dimensional
$\ad(\mathfrak{t})$-invariant vector subspace, and $\Lin$ means a linear span. Using the
restriction (of non-degenerate) inner product $\langle \cdot,\cdot \rangle$ to
$\mathfrak{t}$, we will naturally identify  $\alpha$ with vector
in $\mathfrak{t}$. All such forms (vectors) $\alpha$ are called
\textit{roots} of Lie algebra $(\mathfrak{g},\langle \cdot,\cdot \rangle)$, and
the set $\Delta$ of all such roots $\alpha$ is called \textit{root
system} of Lie algebra $(\mathfrak{g},\langle \cdot,\cdot \rangle)$. It is easy to see that
$[V_{\alpha}, V_{\alpha}]$ is one-dimensional subalgebra of $\mathfrak{t}$ spanned on
the root $\alpha$,
and $[V_{\alpha}, V_{\alpha}]\oplus V_{\alpha}$ is a Lie algebra isomorphic to $su(2)$.
This implies that
vector subspaces $V_{\alpha},\alpha \in \Delta,$ admit bases
$\{u_{\alpha},v_{\alpha}\}$ with the following commutator
relations
\begin{equation}
\label{N} [h,u_{\alpha}]=-\langle \alpha,h \rangle v_{\alpha},\quad
[h,v_{\alpha}]=\langle \alpha,h \rangle u_{\alpha}, \quad h\in \mathfrak{t}, \quad
[u_{\alpha},v_{\alpha}]=-\frac{4}{\langle \alpha, \alpha \rangle}{\alpha}.
\end{equation}

Moreover, for $\alpha\neq \pm\beta,$

\begin{equation}
\label{N1}
[u_{\alpha},u_{\beta}]=N_{\alpha,\beta}u_{\alpha+\beta}+
N_{\alpha,-\beta}u_{\alpha-\beta},
\end{equation}

\begin{equation}
\label{N2}
[v_{\alpha},v_{\beta}]=-N_{\alpha,\beta}u_{\alpha+\beta}+
N_{\alpha,-\beta}u_{\alpha-\beta},
\end{equation}

\begin{equation}
\label{N3}
[u_{\alpha},v_{\beta}]=N_{\alpha,\beta}v_{\alpha+\beta}-
N_{\alpha,-\beta}v_{\alpha-\beta},
\end{equation}

\begin{equation}
\label{N4}
[v_{\alpha},u_{\beta}]=N_{\alpha,\beta}v_{\alpha+\beta}+
N_{\alpha,-\beta}v_{\alpha-\beta},
\end{equation}
where (integer) numbers $N_{\alpha,\beta}$ are defined as follows:
$$
N_{-\alpha,-\beta}=N_{\alpha,\beta},\quad
N_{\alpha,\beta}=\pm(q+1)
$$
for $\alpha, \beta, \alpha+\beta \in \Delta$, where $q$ is the
greatest integer number $j$ such that $\beta-j\alpha \in \Delta$.
We suggest in these formulas that $N_{\gamma,\delta}=0,$ if
$\gamma+\delta$ is not a root. From (\ref{N}) and the invariance of
$\langle \cdot,\cdot \rangle$ with respect to automorphisms of $\mathfrak{g}$,
it is easy to obtain
\begin{equation}
\label{length}
\langle u_{\alpha},u_{\alpha}\rangle=\langle v_{\alpha},v_{\alpha} \rangle=
\frac{4}{\langle\alpha,\alpha \rangle}.
\end{equation}

The formulas above imply

\begin{lemma}
\label{nontr}
$$[V_{\alpha},V_{\beta}]=V_{\alpha+\beta}+V_{\alpha-\beta}.$$
\end{lemma}

The root system $\Delta$ is invariant relative to the Weyl group
$W=W(T)$. Besides this:

(i) For every root $\alpha \in \Delta \subset \mathfrak{t}$
the Weyl group $W$ contains the reflection $\varphi_{\alpha}$ in the plane $P_{\alpha},$
which is orthogonal to the root $\alpha$ (with respect to $\langle \cdot,\cdot \rangle$).

(ii) Reflections from (i) generate $W$.

We list below the root systems of that simple compact Lie groups
which we shall need later:

$$
A_l: e_i-e_j, \quad i\neq j, \quad i,j=0,1,\dots,l.
$$
$$
B_l: \pm e_i, \quad i=1,2,\dots,l; \quad \pm e_i\pm e_j,\quad i<j,
\quad i,j=1,2,\dots, l.
$$
$$
C_l:  \pm 2e_i, i=1,2,\dots,l;\quad  \pm e_i\pm e_j, \quad i<j,
\quad i,j=1,2,\dots, l.
$$
$$
D_l:  \pm e_i\pm e_j, \quad i<j,  \quad i,j=1,2,\dots, l.
$$
$$
g_2:  e_i- e_j;\quad  \pm \left(\sum_{i=1}^3e_i-3e_j\right),\quad
i,j=1,2,3.
$$
$$
f_4: \pm e_i, \quad \pm e_i\pm e_j, \quad \frac{1}{2} (\pm e_1 \pm e_2 \pm e_3
\pm e_4), \quad i,j=1,2,3,4.
$$

Here $A_{l-1}=su(l), B_l=so(2l+1), C_l=sp(l), D_l=so(2l)$.
Let us
remark that all roots of any Lie algebra $A_l,D_l,e_6,e_7, e_8$
have one and the same lengths. The roots of any other simple Lie
algebra have two different lengths, so we have the systems
$\Delta_l\subset \Delta$ and $\Delta_s\subset \Delta$ of all long
and short roots respectively. If $\alpha \in \Delta_l, \beta \in
\Delta_s$ for $B_l,C_l,f_4$ (respectively $g_2$), then
$|\alpha|=\sqrt{2}|\beta| $ (respectively
$|\alpha|=\sqrt{3}|\beta|$), where $|X|=\sqrt{\langle X,X \rangle}$.
In all cases two roots of equal length
may constitute the angles $\frac{\pi}{3}, \frac{\pi}{2},
\frac{2\pi}{3}.$ The roots of different length for $B_l,C_l,f_4$
(respectively, $g_2$) may constitute the angles $\frac{\pi}{4},
\frac{\pi}{2}, \frac{3\pi}{4}$ (respectively $\frac{\pi}{6},
\frac{\pi}{2}, \frac{5\pi}{6}$).

By theorem \ref{pose}, all simply connected homogeneous spaces
$G/H$ of positive Euler characteristic with a simple Lie group $G$
are in one-to one correspondence with Lie subalgebras $\mathfrak{h}$, such
that $\mathfrak{t}\subset \mathfrak{h}\subset \mathfrak{g}$ and
$\mathfrak{h}\neq \mathfrak{g}$; we must identify
subalgebras, which are $\Ad(g)$-conjugate with respect to some
$g\in G$ such that $\Ad(g)(\mathfrak{t})=\mathfrak{t}$.
Any such Lie subalgebra $\mathfrak{h}$ is defined by a class of pairwise
$W$-isomorphic closed symmetric root
subsystems $A$ of $\Delta$, not equal to $\Delta$. By definition,
$A \subset \Delta$ is \textit{closed},
if $\alpha,\beta \in A$ and $\alpha \pm \beta \in \Delta$ imply
$\alpha \pm \beta \in A$,
and \textit{symmetric}, if $-\alpha \in A$ together with
$\alpha \in A.$ Then
\begin{equation}
\label{h} \mathfrak{h}=\mathfrak{t}\oplus \Lin \{ \cup_{\alpha \in
A}V_{\alpha}\},
\end{equation}
where $\Lin$ means a linear span.

\section{On the group $G_2$}

Let's describe all simply connected homogeneous spaces $G/H$ of
positive Euler characteristic for $G=G_2=\Aut(\mathbb{C}\,a)$. For
this we use the considerations from the previous section.

Let us give a description of the root system $\Delta$ of
the Lie algebra $g_2$. There are two simple roots $\alpha,\beta \in \Delta$
such that $\angle(\alpha,\beta)=\frac{5\pi}{6}$ and
$|\alpha|=\sqrt{3}|\beta|$. Then
$$\Delta=\{\pm\alpha, \pm\beta, \pm(\alpha+\beta), \pm(\alpha+2\beta),
\pm(\alpha+3\beta), \pm(2\alpha+3\beta)\}.$$ Under this,
$\pm\alpha, \pm(\alpha+3\beta), \pm(2\alpha+3\beta)$ are all long
roots. One can easily see that all non $W$-isomorphic closed
symmetric root subsystems of $\Delta_{G_2}$, not equal to
$\Delta_{G_2},$  are $\emptyset,$ $\{\pm\alpha\},$ $\{\pm\beta\},$
$\{\pm\beta,\pm(2\alpha+3\beta)\},$
$\{\pm\alpha,\pm(\alpha+3\beta), \pm(2\alpha+3\beta)\}.$

The first three cases give us respectively the following
(generalized) flag manifolds: $G_2/T^2,$ $G_2/SU(2)SO(2),$ and
$G_2/A_{1,3}SO(2),$ where $A_{1,3}$ is a Lie group with Lie
subalgebra of the type $A_1$ of index 3, see \cite{On}.
D.V.~Alekseevsky and A.~Arvanitoyeorgos proved in \cite{AA} that all
$G_2$-invariant Riemannian g.o. metrics on them with the full
connected isometry group $G_2$ are $G_2$-normal. The discussion in
Section \ref{HomManPECH} implies that any $G_2$-invariant metric
on these spaces, whose full connected isometry group is not $G_2,$
is $SO(7)$-normal (symmetric) metric on $G_2/SU(2)\cdot
SO(2)=Gr^{+}_{7,2}$.

The last two closed symmetric root subsystems are maximal, so they
correspond to maximal Lie subalgebras in $g_2,$ which are
respectively isomorphic to $su(2)\oplus su(2)$ and $su(3)$ with
the corresponding compact connected Lie subgroups $SO(4)$ and
$SU(3)$ and homogeneous spaces $G_2/SO(4)$ and $G_2/SU(3)=S^6$,
compare with \cite{On}. In the first (second) case
$$
\mathfrak{p}=V_{\alpha}\oplus V_{\alpha+\beta}\oplus V_{\alpha+2\beta}
\oplus V_{\alpha+3\beta}
$$
(respectively
$$
\mathfrak{p}=V_{\beta}\oplus V_{\alpha+\beta}\oplus V_{\alpha+2\beta}).
$$
It's well-known that irreducible components of a representation of
a compact Lie algebra are uniquely determined up to equivalence.
As a corollary, applying this to the adjoint representation of Lie
subalgebra $\mathfrak{t}\subset \mathfrak{h}$ on $\mathfrak{p}$,
one get that for any $\ad(\mathfrak{h})$-invariant subspace
$V\subset \mathfrak{p}$ there exists an equivalent
$\ad(\mathfrak{h})$-invariant subspace $V'\subset \mathfrak{p}$,
which is a direct sum of the given root vector subspaces
$V_{\gamma}, \gamma \in R$. One can easily see that in both cases
above there is no such $\ad(\mathfrak{h})$-invariant subspace
$V'\subset \mathfrak{p}$ besides $\mathfrak{p}$ and $\{0\}.$ Thus
the space $\mathfrak{p}$ is $\ad(\mathfrak{h})$-irreducible. This
means that the corresponding homogeneous spaces $G_2/H$ are
strongly isotropy irreducible. Then any $G_2$-invariant Riemannian
metric on $G_2/H$ is $G_2$-normal.

\begin{remark}
Note that $G_2/SO(4)$ is irreducible symmetric space, see \cite{Bes}.
\end{remark}

Therefore, we have the following

\begin{pred}\label{gog2}
Any g.o. (any $\delta$-homogeneous, in particular)
Riemannian homogeneous manifold $(G_2/H,\mu)$
of positive Euler characteristic is either $G_2$-normal or $SO(7)$-normal.
\end{pred}

Let us remark at the end that the very last root subsystem
contains only the long roots.

\section{Calculations with roots}

Let suppose that in the Notation of Section \ref{general},
$(M=G/H,\mu)$ is $G$-$\delta$-homogeneous simply connected
indecomposable Riemannian manifold with positive Euler
characteristic. Then $G$ is simple by Theorem \ref{KostN}, and we
have inclusions $T\subset H\subset G,$ where $T$ is a maximal
torus in $G$. Then we have some $\Ad(T)$-invariant $\langle
\cdot,\cdot \rangle$-orthogonal decomposition
$$
\mathfrak{g}=\mathfrak{t}\oplus \Lin \{ \cup_{\gamma\in C}V_{\gamma}\}
\oplus \Lin\{ \cup_{\alpha \in D}V_{\alpha}\},
$$
$C\cup D=\Delta$ is a set of all roots for Lie group $G$ with
respect to Lie algebra $\mathfrak{t}$ of $T$,
$V_{\alpha}=V_{-\alpha}$ and $V_{\gamma}=V_{-\gamma}$ are
two-dimensional "root spaces", and the first two summands give us a
decomposition of the Lie algebra $\mathfrak{h}$ of the Lie group
$H$, the last summand gives $\Ad(H)$-invariant vector subspace
$\mathfrak{p}$.

\begin{pred}
\label{inver} Let $\alpha_1, \dots, \alpha_k\in D$ are linearly
independent roots. Then there is a unique (up to multiplication by
constant) vector $t_c\in \Lin \{\alpha_1,\dots \alpha_k\}$ such that
for some real number $s$, $\Ad(\exp(st_c))=-\Id$ on
$\oplus_{i=1}^{k}V_{\alpha_i}$.
\end{pred}

\begin{proof}
One can easily prove  this by using the dual basis in Euclidean
space $\Lin \{\alpha_1,\dots \alpha_k\}$.
\end{proof}

\begin{pred}
\label{fix1} Let $\alpha_1, \dots, \alpha_k\in D$ are linearly
independent roots and $v=\sum_{i=1}^{k}v_{i},$ where $v_i\in
V_{\alpha_i},i=1,\dots,k,$ are non-zero vectors. Let $u(v)\neq 0,$
(see Remark \ref{notation})
and $C_{v}$ is the set of all $\gamma \in C$ such that
$V_{\gamma}$-component of $u(v)$ is not zero. Then
$$C_{v}\neq
\emptyset, \quad C_{v}\subset \Lin \{\alpha_1,\dots
\alpha_k\}-t_{c}^{\perp},
$$
where $t_{c}^{\perp}$ is the orthogonal compliment in $\Lin
\{\alpha_1,\dots \alpha_k\}$ to the vector $t_c$ from Proposition
\ref{inver}.
\end{pred}

\begin{proof}
Really, if $C_{v}=\emptyset,$ then $u(v):=u\in t$ and by
Proposition \ref{inver}
\begin{equation}
\label{vot}
\Ad(\exp(st_c))(w)=-w, \quad \Ad(\exp(st_c))(u(w))=u(w),
\end{equation}
since $[u,t_c]=0.$ This contradicts to Proposition \ref{fix}. So,
$C_{v}\neq \emptyset.$

Now, if some $\gamma \in C_{v}$ is not in $\Lin \{\alpha_1,\dots
\alpha_k\},$ then one can find a vector $w \in t,$ which is
orthogonal to all $\alpha_1,\dots \alpha_k$, but
$\langle w,\gamma \rangle\neq 0$. Then $[w,v]=0$,
while $[w,u(v)]\neq 0$ which contradicts to
Proposition \ref{fix}.

Finally, if $C_{v}\in t_c^{\perp},$ then one more we have
(\ref{vot}), which is impossible by Proposition \ref{fix}.
\end{proof}

Since roots $\alpha \in D, \gamma \in C$ are non-collinear, the
next proposition follows from Propositions \ref{fix} and
\ref{fix1}.

\begin{pred}
\label{root} If $v\in V_{\alpha}$, $\alpha \in D$,
 then $||v||=|v|=\sqrt{(v,v)}$, i.e. $v$ is a
$\delta$-vector.
\end{pred}

\begin{pred}
\label{two} We have at most two possibilities:
$(\cdot,\cdot)=x \langle \cdot,\cdot \rangle$ on $\mathfrak{p}$ or we have an
$\Ad(H)$-invariant $\langle \cdot,\cdot \rangle$-orthogonal direct
decomposition
$\mathfrak{p}=\mathfrak{p}_1\oplus \mathfrak{p}_2$ such that
$(\cdot,\cdot)=x_1 \langle\cdot,\cdot\rangle$ on
$\mathfrak{p}_1$ and $(\cdot,\cdot)=x_2 \langle \cdot,\cdot \rangle$ on
$\mathfrak{p}_2,$ where $x_1 \neq x_2$.
We have necessarily the first possibility, if all roots of
$G$ have one and the same length.
\end{pred}

\begin{proof}
The elements $\Ad(n)$, $n\in N(T)$, generate on $\mathfrak{t}$ a finite
Weyl group $W=W(T)$. It is known that $W$ is generated by
orthogonal reflections in hyper-planes in $\mathfrak{t}$, orthogonal to roots
in $\Delta\subset \mathfrak{t}$. From this and known classifications of roots
systems of compact simple Lie groups one can easily deduce that
$W$ acts transitively on every set of roots of equal lengths.
There are at most two such sets  in $\Delta$:
the set of all short roots $\Delta_s$ and the set of all long roots $\Delta_l$
(see Section \ref{sv}).
At the same time $\Ad(n)$,
$n\in N(T)$, acts transitively on the set of root vector spaces
$V_{\alpha}$, $\alpha \in \Delta_l \quad \mbox{or} \quad \alpha \in
\Delta_s$. Since $\|\cdot\|$ and $\langle \cdot,\cdot \rangle$ are
$\Ad(G)$-invariant,
we get by Proposition \ref{root} that
$$
(v_{\alpha},v_{\alpha})=\|v_{\alpha}\|^2=\|v_{\beta}\|^2=
(v_{\beta},v_{\beta})
$$
and
$$
\langle v_{\alpha},v_{\alpha} \rangle= \langle v_{\beta},v_{\beta} \rangle,
$$
if $\alpha, \beta \in \Delta_l$ or $\alpha, \beta \in \Delta_s$.
Here $v_{\alpha}\in V_{\alpha}$ mean special vectors from
Section \ref{sv}. From this follow the required statements.

\end{proof}

\begin{corollary}
\label{onel} Any $G$-$\delta$-homogeneous Riemannian manifold
$(G/H,\mu)$ of positive Euler characteristic
with $G=SU(l+1)$, $SO(2l)$, $E_6$, $E_7$, or $E_8$ is
$G$-normal.
\end{corollary}

Therefore, we should examine only the second case in Proposition \ref{two}.
Later on we shall use the following notation in this case:
\begin{equation}\label{not2}
\mathfrak{p}_1=\Lin\{\cup_{\beta \in B} V_{\beta}\},\quad
\mathfrak{p}_2=\Lin \{\cup_{\alpha \in A} V_{\alpha}\}, \quad
A=\Delta_l \cap D,\quad B=\Delta_s \cap D,
\end{equation}
where $\Delta_l$ ($\Delta_s$) means the set of long (respectively, short) roots of
Lie algebra $\mathfrak{g}$.

\begin{lemma}\label{pp}
Let $\mathfrak{g}=\mathfrak{h}\oplus \mathfrak{p}_1\oplus \mathfrak{p}_2$
as above, then
$[\mathfrak{p}_1,\mathfrak{p}_2]\neq 0$.
\end{lemma}

\begin{proof}
Let us suppose that $[\mathfrak{p}_1,\mathfrak{p}_2]=0$ and show that in this case
$\mathfrak{q}:=\mathfrak{p}_1+[\mathfrak{p}_1,\mathfrak{p}_1]$ is an proper
ideal of $\mathfrak{g}$. For this goal it is sufficient to show that
$[\mathfrak{h},\mathfrak{q}]\subset \mathfrak{q}$,
$[\mathfrak{p}_1,\mathfrak{q}]\subset \mathfrak{q}$ and
$[\mathfrak{p}_2,\mathfrak{q}]\subset \mathfrak{q}$.

Since $[\mathfrak{h},\mathfrak{p}_1]\subset \mathfrak{p}_1$ and
$[\mathfrak{p}_2,\mathfrak{p}_1]=0$, then
by the Jacobi identity we get
$[h,[\mathfrak{p}_1,\mathfrak{p}_1]] \subset [[h,\mathfrak{p}_1],\mathfrak{p}_1] \subset
[\mathfrak{p}_1,\mathfrak{p}_1] \subset \mathfrak{q}$ and
$[\mathfrak{p}_2,[\mathfrak{p}_1,\mathfrak{p}_1]]=0$.
Therefore, $[\mathfrak{h},\mathfrak{q}]\subset \mathfrak{q}$ and
$[\mathfrak{p}_2,\mathfrak{q}]=0$.

For any $X,Y \in \mathfrak{p}_1$ and $Z\in \mathfrak{p}_2$ we have
$\langle[X,Y],Z\rangle =- \langle Y, [X,Z]\rangle =0$, since
$[\mathfrak{p}_1,\mathfrak{p}_2]=0$. Hence,
$[\mathfrak{p}_1,\mathfrak{p}_1]\subset \mathfrak{p}_1\oplus \mathfrak{h}$,
and
$[\mathfrak{p}_1,\mathfrak{q}]\subset [\mathfrak{p}_1,\mathfrak{p}_1]+
[\mathfrak{p}_1,[\mathfrak{p}_1,\mathfrak{p}_1]] \subset
[\mathfrak{p}_1, \mathfrak{p}_1]+[\mathfrak{p}_1,\mathfrak{h}] \subset
\mathfrak{q}$.

Consequently, $\mathfrak{q}$ is an ideal of $\mathfrak{g}$. This ideal is proper, since
$\mathfrak{p}_2$ is $\langle \cdot,\cdot \rangle$-ortogonal to $\mathfrak{q}$ (see above).
On the other hand, $\mathfrak{g}$ is a simple Lie algebra
and contains no nontrivial ideal. This contradiction proves Lemma.
\end{proof}

\begin{lemma}
\label{pm} Let suppose that the root system $\Delta$ of a compact
simple Lie algebra $\mathfrak{g}\neq g_2$ contains two roots $\alpha \in
\Delta_l, \quad \beta \in \Delta_s$ of different lengths. Then at
most one of $\alpha +\beta$ or $\alpha-\beta$ is a root in
$\Delta.$
\end{lemma}

\begin{proof}
By previous description of $\Delta$, we have exactly three
possibilities for the angle between $\alpha$ and $\beta$:
$\frac{\pi}{4}, \frac{\pi}{2}, \frac{3\pi}{4}$. In the second case
no one of terms $\alpha +\beta$ or $\alpha-\beta$ is a root.
Otherwise there would be a root, longer than $\alpha,$ which is
impossible. In the first (respectively, third) case $\alpha
-\beta$ (respectively, $\alpha+\beta$) is a root, but not $\alpha
+\beta$ (respectively $\alpha-\beta$).
\end{proof}

\begin{lemma}
\label{lem} 1) The vector subspace
$$
\eta=\mathfrak{t}\oplus \Lin \{\cup_{\alpha \in \Delta_l}V_{\alpha}\}
$$
is a Lie subalgebra in $\mathfrak{g}$.

2) The vector subspace $\eta$ is a maximal subalgebra in
$\mathfrak{g}$, if $G\neq F_4$ and $G\neq Sp(l)$, $l\geq 3$.

3) If $G=Sp(l)$, then all non-collinear roots in
$\Delta_l$  are mutually orthogonal and
$[V_{\alpha_1},V_{\alpha_2}]=0$, if $\alpha_1\neq \pm\alpha_2$ are
roots in $\Delta_l$.

4) If $G=F_4$, then  $\eta$ is isomorphic to $so(8)=spin(8)$.
There is $\ad(\eta)$-invariant decomposition
$$
\mathfrak{g}=f_4=\eta \oplus \mathfrak{q}_1 \oplus \mathfrak{q}_2 \oplus \mathfrak{q}_3,
$$
where
$\mathfrak{q}_3=\Lin \{\cup_{\beta \in \Delta_{a}} V_{\beta}\}$,
and $\Delta_{a}$ consists of all roots in $\Delta_s$ of a form
$\pm e_{i},i=1,2,3,4$; $\mathfrak{q}_1$ ($\mathfrak{q}_2$) is spanned on the root spaces of
roots of the form $1/2(\pm e_1\pm e_2\pm e_3 \pm e_4)$ (see Section 11)
with the odd (respectively, even) number of signs "$-$" in this formula.
All modules $\mathfrak{q}_i$ are $\ad(\eta)$-irreducible,
and
$\mathfrak{r}_i=\eta \oplus \mathfrak{q}_i$, $1\leq i \leq 3$,
is a Lie algebra isomorphic to $so(9)=spin(9)$.
For $i\neq j$ there is an
automorphism of $f_4$ preserving $\eta$ and $\mathfrak{t}$, which
maps $\mathfrak{r}_i$ to $\mathfrak{r}_j$.
Any proper subalgebra of $\mathfrak{g}=f_4$, containing
$\eta$ and different from $\eta$, is one of the subalgebra $\mathfrak{r}_i$,
$1\leq i \leq 3$.
\end{lemma}

\begin{proof}
For $G=G_2$ all statements can be checked directly and easily.

Let $G$ be another simple group (with roots of different lengths),
and $\alpha, \beta \in \Delta_l.$ Then $<\alpha,\beta>=0$ or
$\angle(\alpha,\beta)=\frac{2\pi}{3}$ or
$\angle(\alpha,\beta)=\frac{\pi}{3}.$ In the first case $\alpha\pm
\beta$ cannot be a roots, so $[V_{\alpha},V_{\beta}]=0.$ In the
second (third) case orthogonal reflection of $\mathfrak{t}$ in the
hyperplane, $\langle \cdot,\cdot \rangle$-orthogonal to $\alpha$
(respectively, $-\alpha$), maps the root $\beta$ to the vector
$\alpha+\beta$ (respectively, to $\beta-\alpha$), so this vector
is a long root. At the same time, $\alpha-\beta$ (respectively,
$\beta+\alpha$) is not a root. So, we get
$[V_{\alpha},V_{\beta}]=V_{\alpha+\beta}$ (respectively,
$[V_{\alpha},V_{\beta}]=V_{\alpha-\beta}$). This finished the
proof of the first statement.

The second statement easily follows from the list of all roots of
a simple Lie algebra.

Let us remark that any maximal subalgebra $\theta$ in
$\mathfrak{g}=sp(l)$, $l\geq 3$, (with root system $C_l$),
containing $\eta$, has a form
$$
\theta=\eta\oplus \Lin \{\cup_{\alpha \in \Delta_s-\Delta_{i}} V_{\alpha}\},
$$
where $\Delta_{i}$ contains all roots in $\Delta_s$ of a form $\pm
e_{i}\pm e_{j}$ for a fixed $1\leq i\leq l,$ and all $j\neq i.$
All these Lie algebras $\theta_i$ are mutually isomorphic under
automorphisms of $\mathfrak{g}$ and are isomorphic to the Lie
algebra $\theta_1=sp(1)\oplus sp(l-1)$. So, if $\Theta$ is compact
connected Lie subgroup in $G=Sp(l)$ with Lie algebra $\theta_1$,
then we get the homogeneous space $G/\Theta=Sp(l)/Sp(1)\times
Sp(l-1)=\mathbb{H}P^{(l-1)}.$

All long roots for Lie algebra $sp(l)$ has the form $\pm 2e_i,
1\leq i\leq l,$ so we get the third statement.

One can check the first three statements of 4) directly.
All other statement are proved in \cite{Ad1}.
\end{proof}

\begin{lemma}
\label{import} The module
$\mathfrak{k}:=\mathfrak{h}\oplus \mathfrak{p}_2$ (see (\ref{not2})) is a Lie
subalgebra of $\mathfrak{g}$. As a corollary,
$[\mathfrak{p}_2,\mathfrak{p}_1]\subset \mathfrak{p}_1$.
\end{lemma}

\begin{proof}
It is clear that
$[\mathfrak{h}, \mathfrak{h}] \subset \mathfrak{h}\subset \mathfrak{k}$,
$[\mathfrak{h},\mathfrak{p}_2]\subset \mathfrak{p}_2 \subset \mathfrak{k}$ and
$\eta \subset \mathfrak{k}$ (see Lemma \ref{lem}). Note also that
$[\mathfrak{p}_2, \mathfrak{p}_2] \subset [\eta,\eta] \subset \eta \subset \mathfrak{k}$.
These considerations prove the first statement. The second statement is evident.
\end{proof}

The previous Lemma permits now to use all results of Section 10.

\begin{pred}
\label{ineq} Let suppose that we have the second possibility in
Proposition \ref{two} (so $\Delta$ has roots of two different
lengths), and $\mathfrak{g}\neq g_2$. There are
$\alpha \in A, \quad \beta \in B$ (see (\ref{not2})) such that
$[V_{\alpha},V_{\beta}]\neq 0$. For any such $\alpha, \beta,$
either $\alpha+2\beta \in C$ or $\alpha-2\beta \in C$. Moreover,
the following inequality holds:
\begin{equation}
\label{three} x_1<x_2\leq 2x_1.
\end{equation}
\end{pred}

\begin{proof}
The first statement follows from Lemma \ref{pp}.

If $[V_{\alpha},V_{\beta}]\neq 0,$ then by Lemma \ref{pm} we have
only two possible cases for the angle between $\alpha$ and
$\beta$: $\frac{\pi}{4}$ or $\frac{3\pi}{4}.$ Both cases are quite
similar, so let us consider the second one. In this case
\begin{equation}
\label{one} \alpha +\beta \in \Delta_s, \alpha + 2\beta \in
\Delta_l,\quad
|\alpha|_1=|\alpha+2\beta|_1=\sqrt{2}|\beta|_1=\sqrt{2}|\alpha+\beta|_1,
\end{equation}
where $|X|_1=\sqrt{\langle X,X \rangle}$ for $X\in \mathfrak{g}$.
Then
$$[[u_{\alpha},u_{\beta}],u_{\beta}]=[N_{\alpha,\beta}u_{\alpha+\beta},u_{\beta}]=
N_{\alpha,\beta}(N_{\alpha+\beta,-\beta}u_{\alpha}+N_{\alpha+\beta,\beta}u_{\alpha+2\beta}).$$
Here $N_{\alpha,\beta}=\pm(q+1),$ where $q=\max \{j:
\beta-j\alpha\in \Delta\}=0,$ so $N_{\alpha,\beta}=\pm 1.$
$N_{\alpha+\beta,-\beta}=\pm (p+1),$ where $p=\max \{j:
-\beta-j(\alpha+\beta)\in \Delta\}=1,$ so
$N_{\alpha+\beta,-\beta}=\pm 2.$ $N_{\alpha+\beta,\beta}=\pm
(l+1),$ where $l=\max \{j: \beta-j(\alpha+\beta)\in \Delta\}=1,$
so $N_{\alpha+\beta,\beta}=\pm 2.$ Hence we get
\begin{equation}
\label{two0} [[u_{\alpha},u_{\beta}],u_{\beta}]=2(\pm
u_{\alpha}\pm u_{\alpha+2\beta}),
\end{equation}
where one needs to take only one of four possible choices of
signs. Since $u_{\alpha}\in \mathfrak{p}_2,$ we see from (\ref{two0}) that
$[[u_{\alpha},u_{\beta}],u_{\beta}]_{\mathfrak{p}_2}\neq 0.$ Then
Proposition \ref{t31.9n} implies that
$[[u_{\alpha},u_{\beta}],u_{\beta}]_{\mathfrak{h}}\neq 0.$ It follows from
the formula (\ref{two0}) that
$[[u_{\alpha},u_{\beta}],u_{\beta}]_{\mathfrak{h}}=\pm 2u_{\alpha + 2\beta}.$
Hence $\alpha+ 2\beta \in C$.

In order to prove the inequality $x_1 <x_2$ take a $\delta$-vector
$Y\in \mathfrak{p}_2$ and some $U\in \mathfrak{p}_1$ such that $[U,Y]\neq 0$
(it is possible according to Proposition \ref{root} and Lemma \ref{pp}).
Using the inequality (\ref{th31.dop4}) of Corollary \ref{nado1}
in this case, we get $(x_1-x_2) \langle [U,Y],[U,Y] \rangle \leq 0$, therefore, $x_1 <x_2$.

It is possible to prove the inequality $x_2\leq 2x_1,$ using
Proposition \ref{t31.9n}. But we give a more clear proof. Let's
consider the ($\Ad(G)$-invariant) Chebyshev's norm $\|\cdot \|$ on
$\mathfrak{g}$, corresponding to $G$-$\delta$-homogeneous space
$(G/H,\mu)$ (see Theorem \ref{D}). According to Proposition
\ref{root}, for any root $\alpha \in A$ every $X\in V_{\alpha}$ is
a $\delta$-vector. Therefore, $\|X\|=\sqrt{(X,X)}=\sqrt{x_2}
|X|_1$. Similarly, for any root $\beta \in B$ every $Y\in
V_{\beta}$ is a $\delta$-vector and $\|Y\|=\sqrt{(Y,Y)}=\sqrt{x_1}
|Y|_1$. By above argument we can suppose that (\ref{one}) is
satisfied. Using the equations (\ref{N}), (\ref{one}) and
$\Ad(G)$-invariance of $\|\cdot \|$ and $|\cdot|_1$, we get that
$\|\alpha\|=\sqrt{x_2}|\alpha|_1=\sqrt{2x_2}|\beta|_1$ and
$\|\beta\|=\sqrt{x_1}|\beta|_1$. According to $\Ad(G)$-invariance
of $\|\cdot\|$ and $|\cdot|_1$ we get $\|\gamma\|=\|\beta\|$
($\|\gamma\|=\|\alpha\|$) for any $\gamma \in \Delta_s$
(respectively, for any $\gamma \in \Delta_l$), and by (\ref{one}),
$$
\sqrt{2x_2}|\beta|_1=\|\alpha+2\beta\| \leq \|\alpha
+\beta\|+\|\beta\|=2\|\beta\|=2\sqrt{x_1}|\beta|_1,
$$
which is equivalent to $x_2\leq 2x_1.$ Thus we get inequalities
(\ref{three}).
\end{proof}

\begin{remark}\label{stconv}
As it follows from the proof of inequalities (\ref{three}), for a
$G$-$\delta$-homogeneous Riemannian manifold $(G/H,\mu)$ with
$x_2=2x_1$ the restriction of the Chebyshev's norm $\|\cdot \|$ to
the Cartan subalgebra $\mathfrak{t}$ is not strictly convex norm.
\end{remark}

\begin{corollary}
\label{rflag} Every compact $G$-$\delta$-homogeneous Riemannian flag
manifold $M=G/T$ with a simple compact connected Lie group $G$ is
$G$-normal.
\end{corollary}

\begin{proof}
Let suppose that the space under consideration is not $G$-normal.
Then the first two statements in Proposition \ref{ineq} imply
that $C\neq \emptyset$. But this is impossible for
$\mathfrak{h}=\mathfrak{t}$.
\end{proof}

\begin{pred}
Every vector in $\mathfrak{p}_2$ for every $G$-$\delta$-homogeneous
Riemannian space $(G/H,\mu)$ is a $\delta$-vector.
\end{pred}

\begin{proof}
Let's take in the above notation $t:=x_1\leq x_2$ and an arbitrary
(non-zero) vector $v\in \mathfrak{p}_2.$ Let suppose at first that
$t=x_2.$ In this case the corresponding space $M_2=(G/H,\mu_2)$ is
$G$-normal. Then there is unique Killing vector field on $M_2$,
which as an element of Lie algebra of right-invariant vector
fields on $G$ can be naturally identified with $v\in
\mathfrak{p}_2\subset \mathfrak{g}.$ Then the Chebyshev's norm
$||X||_2=\sqrt{\mu_2(X(y),X(y))}$, where $y=H\in G/H=M$. Now, if
we take $t=x_1< x_2,$ leaving $x_2$ fixed, then for any point
$z\in M$ we will have
$$\sqrt{\mu(X(z),X(z))}\leq \sqrt{\mu_2(X(z),X(z))}\leq
\sqrt{\mu_2(X(y),X(y))}=\sqrt{\mu(X(y),X(y))}.$$

This means that $y=H$ is a point of maximal distortion of $X$
for $\mu$ also, which finishes the proof.
\end{proof}

The following proposition follows from $\Ad(G)$-invariance of
the Chebyshev's norm.

\begin{pred}
\label{inv} The set of all $\delta$-vectors in some vector
subspace $V_1\subset \mathfrak{p}_1$ is invariant under all
$\Ad(g)$, $g\in G$, which leave $V_1$ invariant.
\end{pred}

\section{The special second case}

Now we suppose that we have the second possibility in the
Proposition \ref{two}, hence $\Delta$ contain roots of different
length by Proposition \ref{two} and $G\neq G_2$ by Section
\ref{sv}. So we need to consider only the simple Lie groups $F_4,$
and $Sp(l), SO(2l+1),$ when $l\geq 1.$

If $l=1,$ then the center $C(Sp(1))$ is isomorphic to
$\mathbb{Z}_2$ and $Sp(1)/C(Sp(1))=SO(3).$ The unique nontrivial
Riemannian homogeneous space of positive Euler characteristic in
this case is the symmetric (irreducible) space
$Sp(1)/T=SO(3)/T=S^2$ of rank 1, which is $G$-normal, hence
$G$-$\delta$-homogeneous.

\begin{pred}
\label{decomp} In the notation above, the following statements hold:

1) If $G\neq Sp(l)$, $l\geq 3$, then $A\cup C=\Delta_l$, $B= \Delta_s$.

2) If $G=Sp(l)$, $l\geq 3$, then for every $\alpha \in A$ and
$\gamma \in C,$ $\langle \alpha,\gamma \rangle =0$ and
$[V_{\alpha},V_{\gamma}]=0$.

3) For every $\alpha \in A$ there is an $\beta \in B$ such that
$\langle \alpha,\beta \rangle \neq 0$. If $G\neq G_2,$ then one
(and only one) of the vectors $\alpha +\beta$ or $\alpha-\beta$ is
root in $B$, and $\alpha +2\beta$ (respectively, $\alpha-2\beta$)
is a root in $C$.
\end{pred}

\begin{proof}
The first statement in the case $G\neq F_4$ follows from the
statement 2) of Lemma \ref{lem} and from the inclusion $\eta\subset
\mathfrak{h}\oplus \mathfrak{p}_2$.

Suppose that $G=F_4$. By Lemma \ref{import}, $\mathfrak{p}_2\oplus
\mathfrak{h}$ is a proper Lie subalgebra in $f_4$, which contains $\eta$
by Lemma \ref{lem}. So, by the statement 4) in Lemma \ref{lem},
either $\mathfrak{p}_2\oplus \mathfrak{h}=\eta$, or $\mathfrak{p}_2\oplus
\mathfrak{h}=\mathfrak{r}_i$ for some $1\leq i \leq 3$.
The second case is impossible. Suppose the contrary. Since
$\mathfrak{r}_i=\eta \oplus \mathfrak{q}_i$, we get $\mathfrak{q}_i \subset \mathfrak{h}$.
On the other hand, the module $\mathfrak{q}_i$ generates the Lie algebra $\mathfrak{r}_i$
($(\mathfrak{r}_i,\eta)=(so(9),so(8))$). Since $\mathfrak{h}$ is a proper subalgebra in
$\mathfrak{r}_i$, this is impossible. Therefore, $\mathfrak{p}_2\oplus \mathfrak{h}=\eta$ and
$B$ coincides with the set $\Delta_s$.
This proves the first statement for $G=F_4$.

The second statement follows from the statement 3) of Lemma
\ref{lem}, if $\gamma \in \Delta_l$. The case $\gamma \in
\Delta_s$, can be considered as Lemma \ref{import} above.

Consider now the item 3). For any $\alpha \in A$ there is $\beta \in \Delta_s$ such that
$\gamma:=\alpha +\beta \in \Delta$ (otherwise an angle between $\alpha$ and any
$\beta \in \Delta_s$ is $\pi/2$, with using the Weyl group we get
the same for any root in $A$, but the latter contradicts to Lemma \ref{pp}).
It is clear that $\gamma \in \Delta_s$. Since $\gamma - \beta=\alpha \in A$,
then either $\beta$ or $\gamma$ is not in $C$, hence one of them is in $B$.
Other statements of this item
follow from Lemma \ref{import} and from Proposition \ref{ineq}.
\end{proof}

\begin{pred}
\label{sp} Up to change of indices, in the case of $G=Sp(l),$ we
must have $A=\{\pm 2e_1\}$, $\{\pm e_1 \pm e_i, 1< i \leq
l\}\subset B$.
\end{pred}

\begin{proof}
Let suppose that $A$ contains besides $\pm 2e_1$ (up to change of
indices) yet $\pm 2e_2$. Then by the statement 2) in Proposition
\ref{decomp}, $C$ cannot contain roots of the form $\pm e_i \pm
e_j$, $i<j$, where $i=1$ or $i=2$. So, $B$ contains all roots of
the form $\pm e_1 \pm e_i$, $1<i$, and $\pm e_2 \pm e_j$, $2<j$.
Let consider the root $-e_1+e_2 \in B$. Then
$[V_{2e_1},V_{-e_1+e_2}]=V_{e_1+e_2}$. Now by Lemma \ref{nontr}
$$
[V_{e_1+e_2},V_{-e_1+e_2}]=V_{2e_2}\oplus V_{2e_1}\subset \mathfrak{p}_2.
$$
So, in the previous notation
$$
\alpha:=2e_1, \quad \beta:=-e_1+e_2, \quad \alpha+\beta=e_1+e_2,
\quad \alpha +2\beta=2e_2\in A.
$$
We have got a contradiction with the second part of the second
statement in 3) of Proposition \ref{decomp}.

Now $A=\{\pm 2e_1\}$ and by the first part of the second statement
in 3) of Proposition \ref{decomp}, all roots of the form $\pm
e_1\pm e_i$, $1<i$, must lie in $B$.
\end{proof}

\begin{corollary}
In conditions of Proposition \ref{sp}, $\dim (\mathfrak{p}_{2})=2$.
\end{corollary}

\begin{pred}
\label{flag1} For the case $G=Sp(l), l\geq 2,$ the spaces under
consideration may have only one of the form $M=Sp(l)/U(1)\cdot
Sp(l-1)$ or $Sp(l)/U(1)\times Sp(k_2-1)\times \dots \times
Sp(l-k_m),$ where $1<k_2<\dots <k_m<l, m\geq 2.$
\end{pred}

\begin{proof}
In the Notation of Proposition \ref{sp}, let suppose also that all
other short roots (of the form $\pm e_i \pm e_j,\quad 2\leq i<
j\leq l,$) lie in $C$. In this case we get exactly the first case.
Here $U(1)\cdot Sp(l-1)$ is the centralizer of the root $2e_1\in
\mathfrak{t}$ and $\mathfrak{h}\oplus \mathfrak{p}_2=sp(1)\oplus
sp(l-1)\subset sp(l)$.

Let suppose that in the previous conditions $G=Sp(l)$ and $H\neq
U(1)\times Sp(l-1)$. From Propositions \ref{sp} and the first case
we get that $$U(1)\times Sp(1)^{l-1} \subset H \subset U(1)\times
Sp(l-1)\subset Sp(1)\times Sp(l-1).$$ Therefore, we obtain the
second case from the description of subgroups with maximal rank of
the group $Sp(l),$ obtained in Theorem II of \cite{Wang49}.
\end{proof}

\begin{theorem}
\label{flag}
For the case $G=Sp(l), l\geq 2,$ the spaces under
consideration may have only the form
$M=Sp(l)/U(1)\cdot Sp(l-1)$.
\end{theorem}

\begin{proof} Suppose the contrary, then according to Proposion \ref{flag1}
there is a $\delta$-homogeneous Riemannian manifold
$(G/H=Sp(l)/U(1)\times Sp(k_2-1)\times \cdots \times Sp(l-k_m),
\mu=\mu_{x_1,x_2})$, where $1<k_2<\cdots <k_m<l$, $m\geq 2$, and $x_1\neq x_2$.

Let $K=Sp(1)\times Sp(k_2-1)\times \cdots \times Sp(l-k_m)$,
$H\subset K\subset G$. Then
$\mathfrak{g}=\mathfrak{k}\oplus \mathfrak{p}_1$,
$sp(1)=u(1)\oplus \mathfrak{p}_2$.
We will use notation
$\mathfrak{h}_1=u(1)$,
$\mathfrak{h}_2=sp(k_2-1)\oplus \cdots \oplus sp(l-k_m)$, where
$\mathfrak{h}=\mathfrak{h}_1\oplus \mathfrak{h}_2$.
Let us consider $\Ad(H)$-invariant
submodules $\mathfrak{p}_{1,1}, \mathfrak{p}_{1,2} \subset \mathfrak{p}_1$
such that
$\mathfrak{g}=sp(l)=sp(1)\oplus sp(l-1) \oplus \mathfrak{p}_{1,1}$,
$sp(l-1)= \mathfrak{h}_2 \oplus \mathfrak{p}_{1,2}$, where all sums are
orthogonal with respect to $\langle \cdot,\cdot \rangle$, and
$\mathfrak{p}_1= \mathfrak{p}_{1,1}\oplus \mathfrak{p}_{1,2}$.

Take any $X \in \mathfrak{p}_{1,1}\subset \mathfrak{p}_1$ and any nontrivial
$Y\in \mathfrak{p}_2$. Then there is some $Z\in \mathfrak{h}$ such that
the vector $X+Y+Z$ is a $\delta$-vector. In particular, this vector is
geodesic for $(G/H,\mu)$. Then using Proposition \ref{t31.4}
we get that $[Z,Y]=0$. This means that $Z\in \mathfrak{h}_2$.

Take now any $U\in \mathfrak{p}_{1,2} \subset \mathfrak{p_1}$
and apply the inequality (\ref{th31.dop2})
from Proposition \ref{t31.5} in this situation.
It is clear that $[U,X]\in \mathfrak{p}_{1,1}\subset \mathfrak{p}_1$,
$[U,Y]=0$, $[Z,U]\subset \mathfrak{p}_{1,2}$ and
$\langle [U,X],[U,Z]\rangle=0$. Hence the inequality (10.10) take the form
$ x_1\langle [U,Z],[U,Z]\rangle \leq 0$, consequently, $[U,Z]=0$
for any $U\in \mathfrak{p}_{1,2}$. On the other hand, it is easy to see, that the submodule
$\mathfrak{p}_{1,2}$ generates the Lie algebra $sp(l-1)$
(the pair $(sp(l-1),\mathfrak{h}_2)$ is effective), therefore
$Z$ sits in the center of $sp(l-1)$ and $Z=0$.

Now we use Proposition \ref{t31.4} again. We get that
$[X,Y]=x_1/(x_2-x_1)[X,Z]=0$. Since $X \in \mathfrak{p}_{1,1}$ is arbitrary
we get $[Y, \mathfrak{p}_{1,1}]=0$.
This is impossible since $Y$ is nontrivial and the submodule
$\mathfrak{p}_{1,1}$ generates the Lie algebra $sp(l)$.
Therefore, $(G/H,\mu)$ is not $\delta$-homogeneous.
Theorem is proved.
\end{proof}

\begin{theorem}
\label{tam1} If $G=SO(2l+1)$, where $l\geq 2,$ then the space
$M:=G/H$ under consideration may have only one form
$M=SO(2l+1)/U(l)$.
\end{theorem}

\begin{proof}
The group $G=SO(2l+1)$ has the root system $B_l$. Then the Lie
algebra $\eta$ from Lemma \ref{lem} is isomorphic to the Lie
algebra $so(2l)$ of the Lie group $SO(2l)$ with the root system
$D_l$. In this case $\eta=\mathfrak{h}\oplus \mathfrak{p}_2$ and
$\mathfrak{p}_1=\Lin \{\cup_{\beta \in \Delta_s}V_{\beta}\}$ by the
statement 1) in Proposition \ref{decomp}. Therefore the
homogeneous space $(SO(2l+1)/H,\mu)$ under consideration is fibred
over rank 1 (hence irreducible) symmetric space
$SO(2l+1)/SO(2l)=S^{2l}$. So, the conditions of Theorem 4.1 in the
paper \cite{tam} are satisfied. Then by Table
I on the page 841 of this paper and by Theorem \ref{On} we must have $M=SO(2l+1)/U(l)$.
\end{proof}

\begin{remark}
The spaces in Theorem \ref{flag} and the spaces from Theorem
\ref{tam1} were appeared also in the paper \cite{AA} as
(generalized) flag manifolds, admitting non-normal invariant g.o.
Riemannian metrics. Earlier W.~Ziller proved that all invariant metrics on these
spaces are weakly symmetric (hence, g.o.) \cite{Zi96}.
\end{remark}

\begin{corollary}
\label{tr} For spaces in Theorem \ref{tam1}, every vector in
$\mathfrak{p}_1$ is a $\delta$-vector.
\end{corollary}

\begin{proof}
By the proof of Theorem \ref{tam1}, $p_1$ is naturally identified
with the tangent space at the initial point of a rank 1 symmetric
space $SO(2l+1)/SO(2l)=S^{2l}$, which is two-point homogeneous.
This implies that $\Ad(SO(2l+1))$ acts transitively on the unit
sphere in $(\mathfrak{p}_1,(\cdot,\cdot))$. The proof is finished
by applying of Propositions \ref{root} and \ref{inv}.
\end{proof}

\begin{theorem}
\label{tam2} There is no space $M:=G/H$ under consideration with $G=F_4$.
\end{theorem}

\begin{proof}
At first, we prove that $M=G/H$ with $G=F_4$
may have at most one form $M=F_4/\exp(u(4))$.

In this case $\mathfrak{h}\oplus \mathfrak{p}_2=\eta=so(8)$,
$\mathfrak{p}_1= \mathfrak{q}_1\oplus \mathfrak{q}_2\oplus \mathfrak{q}_3$
(see Lemma \ref{lem} and Proposition \ref{decomp}).
Let's consider a subalgebra
$\mathfrak{r}_3=\eta\oplus \mathfrak{q}_{3}=so(9)=spin(9)=\mathfrak{l}$.
By Proposition \ref{gotg}, the
Riemannian subspace $L/H=Spin(9)/H\subset F_4/H$ is totally
geodesic, hence $\delta$-homogeneous and g.o. space, and also has
positive Euler characteristic. Since $L=Spin(9)$ is a simple
group and the restriction of the Killing form of $f_4$ to
$\mathfrak{l}$ is $\Ad(L)$-invariant, then this restriction must be proportional to
the Killing form of $\mathfrak{l}$. We have $Spin(9)/H=(Spin(9)/C)/(H/C)=SO(9)/(H/C)$,
where $C$ is the common center of $Spin(9)$ and $H$.
Therefore, the Riemannian subspace
$SO(9)/(H/C)$ of $F_4/H$ is not $SO(9)$-normal, if $F_4/H$ is not
$F_4$-normal, because $\mathfrak{l}$ includes vector subspaces
$\mathfrak{p}_2$ and $\mathfrak{q}_3 \subset \mathfrak{p}_1$.

If $SO(9)/(H/C)$ is not $SO(9)$-$\delta$-homogeneous (being $\delta$-homogeneous),
then its full connected isometry group is not equal to $SO(9)$. Therefore, according to
Theorem \ref{On-Shch}, we must have $H/C=U(4)$ and $H=\exp(u(4))$. On the other hand,
if $SO(9)/(H/C)$ is $SO(9)$-$\delta$-homogeneous, then by Theorem
\ref{tam1}, we get again $H=\exp(u(4))$. Note that $\mathfrak{h}=u(4)$
is spanned on the Cartan subalgebra
$\mathfrak{t}$ and on the root spaces of the roots $\pm (e_i-e_j)$,
$1\leq i <j \leq 4$.

Now we shall prove that
the Riemannian manifold
$(G/H=F_4/\exp(u(4)), \mu=\mu_{x_1,x_2})$
is not g.o. for $x_1\neq x_2$.

Note that the submodule $\mathfrak{q}_2$ (see Lemma \ref{lem}) is not
$\ad(\mathfrak{h})$-irreducible. Really,
let us consider a two-dimensional submodule
$\mathfrak{q}\subset \mathfrak{q}_2$,
which is spanned
on the root space of the vectors $\pm 1/2(e_1+e_2+e_3+e_4)$.
It is clear that
$$
\left(\pm (e_i-e_j) \right)+ \left(\pm 1/2(e_1+e_2+e_3+e_4)\right)
$$
is not a root for any $1\leq i <j \leq 4$.
This means that $\mathfrak{q}$ commutes with every root spaces of the roots
$\pm (e_i-e_j)$. Therefore, $\mathfrak{q}$ is invariant under the action of
$\ad(\mathfrak{h})$.

Consider now any
$X\in \mathfrak{q} \subset \mathfrak{q}_2 \subset \mathfrak{p}_1$ and any
$Y\in \mathfrak{p}_2$. If $(F_4/\exp(u(4)), \mu)$ is a g.o. space,
then there is $Z \in \mathfrak{h}$ such that
$X+Y+Z$ is a geodesic vector. If we have $x_1 \neq x_2$ in addition, then
according to Proposition \ref{t31.4}, we get $[X,Y]=x_1/(x_2-x_1)[X,Z]$.
Since $[X,Z]\subset \mathfrak{q}$, we obtain that $[X,Y] \in \mathfrak{q}$
for any $X\in \mathfrak{q}$ and for any $Y\in \mathfrak{p}_2$.
Therefore, the module $\mathfrak{q}$ is $\ad(\eta)$-invariant which
is impossible, since
the module $\mathfrak{q}_2$ (containing $\mathfrak{q}$)
is $\ad(\eta)$-irreducible.
Therefore, $(F_4/\exp(u(4)),\mu)$ is not g.o. for $x_1\neq x_2$.
This finishes the proof.
\end{proof}

\section{On the space $SO(5)/U(2)=Sp(2)/U(1)\cdot Sp(1)=\mathbb{C}\,P^3$}

Here we find all $\delta$-homogeneous metrics on the space
$SO(5)/U(2)$, where $U(2)\subset SO(4) \subset SO(5)$, and the
pair $(SO(5),SO(4))$, $(SO(4),U(2))$ are irreducible symmetric.
Remind that the space $SO(5)/U(2)$ coincides with the space
$Sp(2)/ U(1)\cdot Sp(1)$.

For $A,B \in so(5)$ we define $\langle A,B\rangle
=-1/2\trace(A\cdot B)$. This is an $\Ad(SO(5))$-invariant inner
product on $so(5)$.
A matrix $A+\sqrt{-1}B\in u(2)$, where
$A=\left(
\begin{array}{rr}
0&c\\
-c&0\\
\end{array}
\right)$ and
$B=\left(
\begin{array}{rr}
a&d\\
d&b\\
\end{array}
\right)$ we embed into $so(4)$ via
$A+\sqrt{-1}B \mapsto
\left(
\begin{array}{rr}
A&B\\
-B&A\\
\end{array}
\right)$ in order to get the symmetric pair $(so(4),u(2))$ (see
e.g. \cite{Hel}). Also we use the standard embedding $so(4)$ into
$so(5)$: $A \mapsto \diag(A,0)$.

It is known the following
$\langle\cdot,\cdot\rangle$-orthogonal decomposition:
$$
\mathfrak{g}=so(5)=so(4)\oplus \mathfrak{p}_1 =u(2)\oplus \mathfrak{p}_2
\oplus \mathfrak{p}_1, \quad \mathfrak{p}=\mathfrak{p}_1\oplus \mathfrak{p}_2,
$$
where
$$
u(2)=\left\{ \left(
\begin{array}{rrrrrr}
0&c&a&d&0\\
-c&0&d&b&0\\
-a&-d&0&c&0\\
-d&-b&-c&0&0\\
0&0&0&0&0\\
\end{array}
\right) \,;\quad a,b,c,d \in \mathbb{R}\, \right\},
$$
$$
\mathfrak{p}_1=\left\{ X= \left(
\begin{array}{rrrrrr}
0&0&0&0&k\\
0&0&0&0&l\\
0&0&0&0&m\\
0&0&0&0&n\\
-k&-l&-m&-n&0\\
\end{array}
\right) \,;\quad k,l,m,n \in \mathbb{R}\, \right\},
$$
$$
\mathfrak{p}_2=\left\{ Y=\left(
\begin{array}{rrrrrr}
0&e&0&f&0\\
-e&0&-f&0&0\\
0&f&0&-e&0\\
-f&0&e&0&0\\
0&0&0&0&0\\
\end{array}
\right) \,;\quad e,f \in \mathbb{R}\, \right\},
$$
and the modules $\mathfrak{p}_1$ and $\mathfrak{p}_2$ are
$\Ad(U(2))$-invariant and $\Ad(U(2))$-irreducible. Note that for
vectors $X$ from $\mathfrak{p}_1$ as above we have $\langle
X,X\rangle =k^2+l^2+m^2+n^2$, and for vectors $Y\in
\mathfrak{p}_2$ we have $\langle Y,Y\rangle =2e^2+2f^2$.

Let us consider the invariant metric $\mu={\mu}_{x_1,x_2}$ on
$SO(5)/U(2)$, corresponding to the inner product
$$
(\cdot,\cdot)=x_1\langle \cdot,\cdot\rangle|_{\mathfrak{p}_1}+ x_2\langle
\cdot,\cdot\rangle|_{\mathfrak{p}_2}
$$
for some positive $x_1$ and $x_2$. We know that every such metric
is a g.o.-metric \cite{Zi96},\cite{tam}.
From the discussion in Section
\ref{HomManPECH} we get the following

\begin{pred}\label{vspom0}
The full connected isometry group of $(SO(5)/U(2),\mu)$ is
$SO(5)$, excepting the case $x_2=2x_1$, where the full connected
isometry group is $SO(6)/\{\pm I\}$, and the metric $\mu$ is
$SO(6)$-normal (in the last case $(SO(5)/U(2),\mu)$ is isometric
to the complex projective space $\mathbb{C}P^3$ with the standard
Fubini metric).
\end{pred}

Let $E_{i,j}$ be a $(5\times 5)$-matrix, whose $(i,j)$-th entry is
equal to $1$, and all other entries are zero. For any $1\leq i <j
\leq 5$ put $F_{i,j}=E_{i,j}-E_{j,i}$. Let consider the following
subspace of $\mathfrak{p}=\mathfrak{p}_1\oplus \mathfrak{p}_2$:
$$
\mathfrak{q}=\mathbb{R} \cdot F_{1,5} \oplus \mathbb{R}\cdot
(F_{1,4}-F_{2,3}).
$$

\begin{pred}\label{vspom1}
For any vector $V\in \mathfrak{p}$ there is $a \in H=U(2)$ such that
$\Ad(a)(V)\in \mathfrak{q}$.
\end{pred}

\begin{proof}
Let $V=X+Y$, where $X\in \mathfrak{p}_1$ and $Y\in \mathfrak{p}_2$.
We know by (the proof of) Corollary \ref{tr} that $\Ad(U(2))$
acts transitively on
the unit sphere in $\mathfrak{p}_1$. Therefore, we may assume that
$X=bF_{1,5}$ for some $b\in \mathbb{R}$. We have
$$
Y=c_1(F_{1,2}-F_{3,4})+ c_2(F_{1,4}-F_{2,3})
$$
for some real $c_1$ and $c_2$. Note that $[F_{2,4},X]=0$.
Therefore, $X$ is invariant under $\Ad(a)$, where
$a=\exp(tF_{2,4})$. On the other hand,
$$
\Ad(a)(Y)=\widetilde{c}_1(F_{1,2}-F_{3,4})+
\widetilde{c}_2(F_{1,4}-F_{2,3})\in \mathfrak{p}_2,
$$
where
$$
\widetilde{c}_1=c_1\cos(t)+c_2\sin(t),\quad
\widetilde{c}_2=c_2\cos(t)-c_1\sin(t).
$$
For some suitable $t\in \mathbb{R}$ we get that
$\widetilde{c}_1=0$. Therefore,
$\Ad(a)(V)=bF_{1,5}+\widetilde{c}_2(F_{1,4}-F_{2,3}) \in \mathfrak{q}$.
\end{proof}

\begin{pred}\label{vspom2}
Let $W=X+Y+Z$, where $X+Y \in \mathfrak{q}$ and $Z\in
\mathfrak{h}=u(2)$, be a non-trivial geodesic vector on
$(SO(5)/U(2),\mu)$, $x_2 \neq x_1$, $x_2\neq 2x_1$. Then we have
one of the following possibilities:

1)
$W=bF_{1,5}+\frac{x_2}{x_1}cF_{1,4}+\frac{x_2-2x_1}{x_1}cF_{2,3}$
for some $b\neq 0$, $c\neq 0$;

2) $W=
d(F_{1,4}-F_{2,3})+a_1(F_{1,2}+F_{3,4})+a_2(F_{1,4}+F_{2,3})+
a_3(F_{1,3}-F_{2,4})$ for some $d\neq 0$, $a_1,a_2,a_3 \in
\mathbb{R}$;

3) $W= eF_{1,5}+fF_{2,4}$ for some $e\neq 0$ and $f\in
\mathbb{R}$.
\end{pred}

\begin{proof}
Let $W=X+Y+Z$, where $X=bF_{1,5}\in p_1$, $Y=c(F_{1,4}-F_{2,3})$,
and
$Z=b_1(F_{1,2}+F_{3,4})+b_2(F_{1,4}+F_{2,3})+b_3F_{1,3}+b_4F_{2,4}$.
Since $W$ is geodesic vector, then from Proposition \ref{t31.4} we
have
$$
[Z,Y]=0, \quad [X,Y]=\frac{x_1}{x_2-x_1}[X,Z].
$$
Direct calculations show that
$$[Z,Y]=c(b_3+b_4)(F_{1,2}-F_{3,4}), [X,Y]=bcF_{4,5},
[X,Z]=b(b_1F_{2,5}+b_3F_{3,5}+b_2F_{4,5}).$$

If $b\neq 0$ and $c\neq 0$, then $b_1=b_3=b_4=0$ and
$b_2=\frac{x_2-x_1}{x_1}c$.

If $b=0$ and $c \neq 0$, then $b_4=-b_3$.

If $b\neq 0$ and $c=0$, then we have $b_1=b_2=b_3=0$.

The proposition is proved.
\end{proof}

\begin{pred}\label{vspom3}
The Riemannian manifold $(SO(5)/U(2),\mu)$ is
$SO(5)-\delta$-homogeneous if and only if for every $b\neq 0$ and
every $c\neq 0$ the vector
$$
W=\left(
\begin{array}{ccccr}
0 & 0 & 0 & \frac{x_2}{x_1}c & b \\
0 & 0 & \frac{x_2-2x_1}{x_1}c & 0 & 0 \\
0 & \frac{2x_1-x_2}{x_1}c & 0 &  0 & 0 \\
-\frac{x_2}{x_1}c & 0 & 0 & 0 & 0 \\
-b & 0 & 0 & 0 & 0
\end{array}
 \right) =b F_{1,5}+ \frac{x_2}{x_1}c F_{1,4}+ \frac{x_2-2x_1}{x_1}c F_{2,3}
$$
is $\delta$-vector on $(SO(5)/U(2),\mu)$.
\end{pred}

\begin{proof}
If $(SO(5)/U(2),\mu)$ is $SO(5)-\delta$-homogeneous, then for
every vector of the form $V=X+Y$, where $X=bF_{1,5}\in
\mathfrak{p}_1$, $Y=c(F_{1,4}-F_{2,3}) \in \mathfrak{p}_2$, $b\neq
0$, $c\neq 0$, there is $Z\in \mathfrak{h}$ such that the vector
$W=X+Y+Z$ is $\delta$-vector. In particular, $W$ is geodesic
vector. According to Proposition \ref{vspom2}, we get that
$$
W=b F_{1,5}+ \frac{x_2}{x_1}c F_{1,4}+ \frac{x_2-2x_1}{x_1}cF_{2,3}.
$$
Therefore, this $W$ is a $\delta$-vector.

Let us suppose now that all vectors of the form
$$
W=b F_{1,5}+ \frac{x_2}{x_1}c F_{1,4}+ \frac{x_2-2x_1}{x_1}cF_{2,3},
$$
where $b\neq 0$ and $c\neq 0$, are $\delta$-vectors. Since the
limit of any sequence of $\delta$-vectors is a $\delta$-vector
itself, we get that the vectors $W$ as above are $\delta$-vectors
for $b=0$ or $c=0$ also.

Therefore, for any vector $X+Y\in \mathfrak{q}$ there is $Z\in \mathfrak{h}$
such that the vector $X+Y+Z$ is $\delta$-vector. Using Proposition
\ref{vspom1}, we get that $(SO(5)/U(2),\mu)$ is
$SO(5)-\delta$-homogeneous in this case.
\end{proof}

\begin{lemma}\label{sim}
For every $b,c,x_1,x_2 \in \mathbb{R}$ with the properties
$$
b\neq 0, \quad x_1\neq 0, \quad 2x_1> x_2,
$$
the following inequality is fulfilled:
$$
\left(|c|(2x_1-x_2)+\sqrt{b^2x_1^2+c^2x_2^2}\right)^2x_2 <
2x_1^2(x_1b^2+2x_2c^2).
$$
\end{lemma}

\begin{proof} It is enough to consider the case $x_2>0$. In this case
we have the following chain of equivalent inequalities.
$$
\left(c^2(2x_1-x_2)^2+b^2x_1^2+c^2x_2^2+
2|c|(2x_1-x_2)\sqrt{b^2x_1^2+c^2x_2^2}\right)x_2 <
2x_1^3b^2+4x_1^2x_2c^2;
$$
$$
2|c|(2x_1-x_2)\sqrt{b^2x_1^2+c^2x_2^2}x_2 <
2x_1^3b^2+4x_1^2x_2c^2- c^2(2x_1-x_2)^2x_2-b^2x_1^2x_2-c^2x_2^3 =
$$
$$
(2x_1-x_2)x_1^2b^2+2x_2^2(2x_1-x_2)c^2;
$$
$$
2|c|\sqrt{b^2x_1^2+c^2x_2^2}x_2 < x_1^2b^2+2x_2^2c^2;
$$
$$
4c^2(b^2x_1^2+c^2x_2^2)x_2^2 =4x_1^2x_2^2b^2c^2+4x_2^4c^4<
x_1^4b^4+4x_1^2x_2^2b^2c^2+4x_2^4c^4=(x_1^2b^2+2x_2^2c^2)^2.
$$
\end{proof}

\begin{pred}\label{vspom5}
If $2x_1\geq x_2 \geq x_1$, then the Riemannian manifold
$(SO(5)/U(2),\mu)$ is $SO(5)-\delta$-homogeneous.
\end{pred}

\begin{proof}
We may assume by continuity, that $x_1<x_2<2x_1$.

According to Proposition \ref{vspom3}, we only need to prove that
every vector of the form
$$
W =b F_{1,5}+ \frac{x_2}{x_1}c F_{1,4}+ \frac{x_2-2x_1}{x_1}cF_{2,3},
$$
where $b\neq 0$ and $c\neq 0$, is $\delta$-vector on
$(SO(5)/U(2),\mu)$.

Let us consider the orbit $O(W)$ of $W$ under the action of
$\Ad(G)=\Ad(SO(5))$. Since $O(W)$ is compact, there is
$\widetilde{W} \in O(W)$ such that
$$
(\widetilde{W}|_{\mathfrak{p}},\widetilde{W}|_{\mathfrak{p}})
\geq (V|_{\mathfrak{p}},V|_{\mathfrak{p}})
$$
for every $V\in O(W)$.

Therefore, $\widetilde{W}$ is a $\delta$-vector. According to
Proposition \ref{vspom1} we may assume, that $\widetilde{W}|_{\mathfrak{p}}\in
\mathfrak{q}$. Now it is sufficient to show that
$$
(\widetilde{W}|_{\mathfrak{p}},\widetilde{W}|_{\mathfrak{p}}) \leq
(W|_{\mathfrak{p}},W|_{\mathfrak{p}})
$$

We shall use the following idea. Since $\widetilde{W} \in O(W)$,
then the matrices $-W^2$ and $-\widetilde{W}^2$ has one and the
same set of eigenvalues. The eigenvalues of $-W^2$ are the
following:
$$
0,\quad \frac{c^2(2x_1-x_2)^2}{x_1^2},\quad
\frac{b^2x_1^2+c^2x_2^2}{x_1^2},
$$
where two last eigenvalues are of multiplicity $2$. Since $x_2>
x_1$, we obviously get
$$
b^2x_1^2+c^2x_2^2 > c^2(2x_1-x_2)^2.
$$
Note also that $(W|_{\mathfrak{p}},W|_{\mathfrak{p}})=x_1b^2+2x_2c^2$.

Since $\widetilde{W}$ is geodesic vector and
$\widetilde{W}|_{\mathfrak{p}} \in \mathfrak{q}$, then by
Proposition \ref{vspom2} we have one of the following
possibilities:

1) $\widetilde{W}=\widetilde{b}F_{1,5}+
\frac{x_2}{x_1}\widetilde{c}F_{1,4}+\frac{x_2-2x_1}{x_1}\widetilde{c}F_{2,3}$
for some $\widetilde{b}\neq 0$, $\widetilde{c}\neq 0$;

2) $\widetilde{W}=
d(F_{1,4}-F_{2,3})+a_1(F_{1,2}+F_{3,4})+a_2(F_{1,4}+F_{2,3})+
a_3(F_{1,3}-F_{2,4})$ for some $d\neq 0$, $a_1,a_2,a_3 \in
\mathbb{R}$;

3) $\widetilde{W}= eF_{1,5}+fF_{2,4}$ for some $e\neq 0$ and $f\in
\mathbb{R}$.

Let us consider these cases separately.

{\bf Case 1)}. Eigenvalues of $-\widetilde{W}^2$ in this case are
the following:
$$
0,\quad \frac{\widetilde{c}^2(2x_1-x_2)^2}{x_1^2},\quad
\frac{\widetilde{b}^2x_1^2+\widetilde{c}^2x_2^2}{x_1^2},
$$
where two last eigenvalues are of multiplicity $2$. Since
$\widetilde{b}^2x_1^2+\widetilde{c}^2x_2^2 >
\widetilde{c}^2(2x_1-x_2)^2$ (remind that $x_2 > x_1$) and
$\widetilde{W}\in O(W)$, we get that
$$
\widetilde{b}^2x_1^2+\widetilde{c}^2x_2^2 =
b^2x_1^2+c^2x_2^2,\quad \widetilde{c}^2(2x_1-x_2)^2=
c^2(2x_1-x_2)^2,
$$
which implies $c^2=\widetilde{c}^2$ and $b^2=\widetilde{b}^2$
(since $2x_1>x_2$). Therefore
$$
(\widetilde{W}|_{\mathfrak{p}},\widetilde{W}|_{\mathfrak{p}})=
x_1\widetilde{b}^2+2x_2\widetilde{c}^2=
x_1b^2+2x_2c^2= (W|_{\mathfrak{p}},W|_{\mathfrak{p}}).
$$

{\bf Case 2).} In this case the eigenvalues of $-\widetilde{W}^2$
are the following:
$$
0,\quad d^2+a_1^2+a_2^2+a_3^2-2\sqrt{d^2(a_1^2+a_2^2+a_3^2)},\quad
d^2+a_1^2+a_2^2+a_3^2+2\sqrt{d^2(a_1^2+a_2^2+a_3^2)},
$$
where two last eigenvalues are of multiplicity $2$.

Since $\widetilde{W}\in O(W)$, we obtain
$$
(|d|-|s|)^2=d^2+s^2-2\sqrt{d^2s^2}=
\frac{c^2(2x_1-x_2)^2}{x_1^2},\quad
(|d|+|s|)^2=d^2+s^2+2\sqrt{d^2s^2}=
\frac{b^2x_1^2+c^2x_2^2}{x_1^2},
$$
where $s^2=a_1^2+a_2^2+a_3^2$. We get from these equations
$$
2|d|=(|d|-|s|)+(|d|+|s|)\leq \frac{|c|(2x_1-x_2)}{x_1}+
\frac{\sqrt{b^2x_1^2+c^2x_2^2}}{x_1}.
$$
Using Lemma \ref{sim}, we get
$$
4d^2x_1^2x_2 \leq
\left(|c|(2x_1-x_2)+\sqrt{b^2x_1^2+c^2x_2^2}\right)^2x_2 <
2x_1^2(x_1b^2+2x_2c^2).
$$
Therefore
$$
(\widetilde{W}|_{\mathfrak{p}},\widetilde{W}|_{\mathfrak{p}})=
2x_2d^2 < x_1b^2+2x_2c^2=(W|_{\mathfrak{p}},W|_{\mathfrak{p}}).
$$

{\bf Case 3).} In this case the eigenvalues of $-\widetilde{W}^2$
are the following:
$$
0,\quad e^2, e^2, f^2, f^2.
$$
Therefore
$$
e^2= \frac{c^2(2x_1-x_2)^2}{x_1^2}\quad\mbox{or}\quad e^2=
\frac{b^2x_1^2+c^2x_2^2}{x_1^2},
$$
Since $2x_1 >x_2>x_1$, we get
$$
x_1b^2+2x_2c^2> \frac{b^2x_1^2+c^2x_2^2}{x_1}>
\frac{c^2(2x_1-x_2)^2}{x_1},
$$
which implies $x_1b^2+2x_2c^2>x_1e^2$. Therefore
$$
(\widetilde{W}|_{\mathfrak{p}},\widetilde{W}|_{\mathfrak{p}})=x_1e^2
<x_1b^2+2x_2c^2=(W|_{\mathfrak{p}},W|_{\mathfrak{p}}).
$$

The above considerations prove that $W$ is a $\delta$-vector on
$(SO(5)/U(2),\mu)$. This proves the proposition.
\end{proof}

\begin{theorem}\label{main}
The Riemannian manifold $(SO(5)/U(2),\mu={\mu}_{x_1,x_2})$ is
$\delta$-homogeneous if and only if $x_1\leq x_2 \leq 2x_1$. For
$x_2=x_1$ it is $SO(5)$-normal homogeneous; for $x_2=2x_1$ it is
$SO(6)$-normal homogeneous; for $x_2\in (x_1,2x_1)$ it is not
normal homogeneous with respect to any its isometry group, but
$SO(5)$-$\delta$-homogeneous.
\end{theorem}

\begin{proof}
If $(SO(5)/U(2),\mu={\mu}_{x_1,x_2})$ is $\delta$-homogeneous,
then it is $SO(6)$-$\delta$-homogeneous or
$SO(5)$-$\delta$-homogeneous, see Theorem \ref{Shchet}. In the
first case it is $SO(6)$-homogeneous. Then by Example \ref{2pec},
we have $x_2=2x_1$. In the second case, by Proposition \ref{ineq}
we get $x_1\leq x_2 \leq 2x_1$. On the other hand, for $x_2=x_1$
and for $x_2=2x_1$ the metric $\mu$ is $SO(5)$-normal homogeneous
and $SO(6)$-normal homogeneous respectively (see Example
\ref{2pec}). From Proposition \ref{vspom5} we get that the
Riemannian manifold $(SO(5)/U(2),\mu)$ is $\delta$-homogeneous for
$2x_1> x_2> x_1$. The theorem is proved.
\end{proof}

\begin{remark}
According to Theorem \ref{Kon57}, the metrics in Theorem \ref{main}  with
the condition $x_2\in (x_1, 2x_2)$ are not naturally
reductive (with respect to any isometry group) in spite of the fact that
they are $\delta$-homogeneous.
\end{remark}

\begin{remark}
It follows from \cite{V} that the Riemannian manifolds in Theorem
\ref{main} have positive sectional curvatures and their (exact)
pinch constant is $\varepsilon=(\frac{x_2}{4x_1})^2.$ This means
that if we scale them so that their maximal sectional curvature
will be 1, then minimal sectional curvature will be $\varepsilon.$
\end{remark}





\begin{thebibliography}{99}

\bibitem{Ad}
{\sl Adams, J. Frank.} Lectures on Lie groups. W.A.Benjamin, Inc.,
New York-Amsterdam, 1969.

\bibitem{Ad1}
{\sl Adams, J. Frank.} Lectures on exceptional Lie groups.
The University of Chicago Press, Chicago-London, 1996.

\bibitem{Al68}
{\sl Aleksevskii~D.V.} Compact quaternion spaces.
Funk. Anal. Pril., {\bf 2(2)} (1968), 11--20.

\bibitem{AA}
{\sl Alekseevsky, D.V. and Arvanitoyeorgos, A.} Riemannian flag
manifolds with homogeneous geodesics, to appear in Trans. AMS.


\bibitem{AkhVinb}
{\sl Akhiezer~D.N and Vinberg~E.B.}
Weakly symmetric spaces and spherical varieties. Transformation
Groups, {\bf 4} (1999), 3--24.


\bibitem{B}
{\sl Berestovski\u\i, V.N.} Homogeneous Busemann $G$-spaces.
(Russian) Sibirsk. Mat. Zh., {\bf 23 (2)} (1982), 3--15.


\bibitem{B2}
{\sl Berestovski\u\i, V.N.} Homogeneous Riemannian manifolds of
positive Ricci curvature. (Russian) Mat. Zametki, {\bf 58 (3)}(1995),
334--340; translation in Math. Notes, {\bf 58 (3-4)} (1996), 905--909.


\bibitem{BerP}
{\sl Berestovskii, Valera, and Plaut, Conrad}. Homogeneous spaces
of curvature bounded below.  J. Geom. Anal., {\bf 9 (2)} (1999), 203--219.


\bibitem{BG}
{\sl Berestovskii, V.N. and Guijarro, Luis}. A metric
characterization of Riemannian submersions. Ann. Global Anal.
Geom., {\bf 18 (6)} (2000), 577--588.

\bibitem{Berger}
{\sl Berger, M.} Les varietes riemanniennes homogenes normales a
courbure strictement positive. Ann.\ Sc.\ Norm.\ Sup.\ Pisa,
{\bf 15} (1961), 179--246.

\bibitem{Berger66}
{\sl Berger, M.} Trois remarques sur les vari\'{e}t\'{e}s
riemanniennes \`{a} courbure positive. C. R. Acad. Sci. Paris
S\'{e}r. A-B, {\bf 263} (1966), 76--78.


\bibitem{Bes}
{Besse, A.L.} Einstein Manifolds. Springer-Verlag, Berlin,
Heidelberg, New York, London, Paris, Tokyo, 1987.


\bibitem{BorSieb}
{\sl Borel, A. and de Siebenthal, J.},
Les sous-groups ferm\'{e}s de rang maximum des groups de Lie clos.
Comment. Math. Helv., {\bf 23} (1949), 200--221.


\bibitem{Burb4}
{\sl Bourbaki, N.} Groupes et alg\'{e}bras de Lie. Chapter 9.
Groupes de Lie r\'{e}els compacts. MASSON, Paris, 1982.


\bibitem{Bot}
{Bott, R.} Vector fields and characteristic numbers. Michigan
Math. J., {\bf 14} (1967), 231--244.

\bibitem{Bus}
{\sl Busemann, Herbert.} The geometry of geodesics. Academic Press
Inc., New York, N.Y., 1955.

\bibitem{CG}
{Cheeger, J. and Gromoll, D.} The splitting theorem for manifolds
of nonnegative Ricci curvature. J. Dif. Geometry, {\bf 6(6)}
(1971), 119--128.


\bibitem{KF}
{\sl Cohn-Vossen, S.} Existenz k\"urzester Wege. Compositio Math.,
{\bf 3} (1936), 441--452.

\bibitem{DZ}
{\sl D'Atri, J.E. and Ziller, W.} Naturally reductive metrics and
Einstein metrics on compact Lie groups. Memoirs  Amer. Math. Soc.,
{\bf 18} (1979), no. 215, 1--72.

\bibitem{DKN}
{\sl Du\v{s}ek~Z., Kowalski~O., Nik\v{c}evi\'c~S.~\v{Z}.~}
New examples of Riemannian g.o. manifolds in dimension 7. Differential
Geom. Appl., {\bf 21} (2004), 65--78.

\bibitem{F}
{\sl Freudenthal, H.} Clifford-Wolf-Isometrien symmetrischer
R\"{a}ume, Math. Ann., {\bf 150} (1963), 136--149.


\bibitem{GO}
{\sl Gorbatsevich, V.V. and Onishchik, A.L.} Lie transformation
groups. Lie groups and Lie algebras, I, 95--235. Encyclopedia
Math. Sci., 20, Springer, Berlin, 1993.

\bibitem{Hel}
{\sl Helgason, S.} Differential geometry and symmetric spaces.
Academic Press Inc., New-York, 1962.

\bibitem{HS}
{\sl Hopf, H. and Samelson, H.} Ein Satz \"uber die Wirkungr\"aume
geschlossener Liescher Gruppen. (German). Comment. Math. Helv.,
{\bf 13} (1940-41), 240--251.

\bibitem{J}
{\sl Jacobson, Nathan.} Lie algebras. Interscience Tracts in Pure
and Applied Mathematics, No. 10. Interscience Publishers (a
division of John Wiley and Sons), New York-London, 1962.

\bibitem{Ker}
{\sl Kerr, M.} Some new homogeneous Einstein metrics on symmetric
spaces. Trans. Amer. Math. Soc., {\bf 348} (1996), 153--171.


\bibitem{K}
{\sl Kobayashi, Sh.} Transformation groups in differential
geometry. Springer, Berlin, 1972.


\bibitem{KN}
{\sl Kobayashi, S. and Nomizu, K.}
Foundations of differential
geometry. Vol.~I -- A Wiley-Interscience Publication, New York,
1963; Vol.~II -- A Wiley-Interscience Publication, New York, 1969.


\bibitem{Kost56}
{\sl Kostant, B.} On differential geometry and homogeneous spaces.
Proc. Math. Acad. Sci. USA, {\bf 42} (1956), 258--261, 354--357.

\bibitem{Kost57}
{\sl Kostant, B.} On holonomy and homogeneous spaces.
Nagoya Math. J., {\bf 12} (1957), 31--54.


\bibitem{KV}
{\sl Kowalski, O. and  Vanhecke, L.} Riemannian manifolds with
homogeneous geodesics. Boll. Unione Mat. Ital. vii.Ser. B, {\bf 5 (1)}
(1991), 189--246.

\bibitem{Kr}
Kr\"amer~M.~ Sph\"arische Untergruppen in kompakten zusammenh\"angenden
Liegruppen. Compositio Math., {\bf 38(2)} (1979), 129--153.

\bibitem{Laur}
{\sl Lauret~J.}
Commutative spaces which are not weakly symmetric. Bull. London Math. Soc.,
{\bf 30} (1998), 29--36.

\bibitem{MT}
{\sl Mimura, M. and Toda, H.} Topology of Lie groups, I and II.
AMS Translations 91, Providence, Rhode Island, 1991.


\bibitem{Neill}
{\sl O'Neill, B.} The fundamental equations of a submersion,
Michigan Math. J., {\bf 13(4)} (1966), 459--469.

\bibitem{Ng}
{\sl Nguy\v{e}\~n~H.~D.} Compact weakly symmetric spaces and spherical pairs.
Proc.\ AMS, {\bf 128(11)} 2000, 3425--3433.

\bibitem{On}
{\sl Onishchik, A.L.} Topology of Transitive Transformation
Groups. Johann Ambrosius Barth: Leipzig, Berlin, Heidelberg, 1994.

\bibitem{On92}
{\sl Onishchik, A.L.} The group of isometries of a compact
Riemannian homogeneous space. Differential geometry and its
applications (Eger, 1989), 597--616, Colloq. Math. Soc. Jano\'s
Bolyai, 56, North-Holland, Amsterdam, 1992.

\bibitem{Selb}
{\sl Selberg~A.}
Harmonic analysis and discontinuous groups in weakly symmetric spaces with
applications to Dirichlet series. J. Indian Math. Soc. N.S., {\bf 20} (1956), 47--87.

\bibitem{Shchet1}
{\sl Shchetinin, A.N.} On a class of compact homogeneous spaces.
I. Ann. Global Anal. Geom., {\bf 6} (1988), 119--140.

\bibitem{Shchet2}
{\sl Shchetinin, A.N.} On a class of compact homogeneous spaces.
II. Ann. Global Anal. Geom., {\bf 8} (1990), 227--247.


\bibitem{tam}
{\sl Tamaru, H.} Riemannian g.o. spaces fibered over irreducible
symmetric spaces. Osaka J. Math., {\bf 36} (1999), 835--851.


\bibitem{tam1}
{\sl Tamaru, H.} Riemannian geodesic orbit metrics on fiber
bundles. Algebra, Groups and Geometries, {\bf 15} (1998), 55--67.

\bibitem{T}
{\sl Toponogov, V.A.} Metric structure of Riemannian spaces of
nonnegative curvature containing direct lines (Russian). Sib.
Math. Zh., {\bf 5(5)} (1964), 1358--1369.

\bibitem{V}
{\sl Vol'per, D.E.} Sectional curvatures of nonstandard metrics on
$CP^{2n+1}$ (Russian). Sibirsk. Mat. Zh., {\bf 40(1)} (1999), 49--56;
translation in Siberian Math. J., {\bf 40(1)} (1999), 39--45.


\bibitem{Wallach72}
{\sl Wallach, N.R.} Compact homogeneous Riemannian manifolds with
strictly positive curvature, Annals of Math., {\bf 96(2)} (1972),
277--295.


\bibitem{Wang49}
{\sl Wang, H.C.}
Homogeneous spaces with non-vanishing Euler characteristic.
Annals of Math., 2nd Ser., {\bf 50(4)} (1949), 925--953.


\bibitem{Wa-Zi1}
{\sl Wang, M. and Ziller, W.}
On isotropy irreducible Riemannian manifolds.
Acta Math., {\bf 166} (1991), 223--261.

\bibitem{Wa}
{\sl Warner, Frank W.} Foundations of differentiable manifolds and
Lie groups. Corrected reprint of the 1971 edition. Graduate Texts
in Mathematics, 94. Springer-Verlag, New York-Berlin, 1983.

\bibitem{W}
{\sl J.A.~Wolf}, Spaces of constant curvature. University of
California, Berkley, 1972.


\bibitem{Zi1}
{\sl Ziller, W.} Homogeneous Einstein Metrics on Spheres and
Projective Spaces.  Math. Ann., {\bf 259} (1982), 351-358.


\bibitem{Zi96}
{\sl Ziller, W.} Weakly symmetric spaces, 355--368. In: Progress
in Nonlinear Differential Equations. V.~20. Topics in geometry: in
memory of Joseph D'Atri. Birkh{\"a}user, 1996.


\end{thebibliography}
\end{document}